\documentclass[a4paper,11pt]{article}
\usepackage{amssymb,amscd,amsfonts,amsbsy}
\usepackage{enumerate}
\usepackage{epsfig}
\input amssymb.sty
\input euscript.sty

\def\ring{\mathaccent"0017 }

\newcommand{\RR}{{\mathbb{R}}}
\newcommand{\NN}{{\mathbb{N}}}
\newcommand{\ZZ}{{\mathbb{Z}}}
\newcommand{\CC}{{\mathbb{C}}}
\newcommand{\meanint}{{\int{\mkern-19mu}-}}

\newtheorem{proposition}{Proposition}[section]
\newtheorem{theorem}[proposition]{Theorem}
\newtheorem{lemma}[proposition]{Lemma}
\newtheorem{corollary}[proposition]{Corollary}

\newtheorem{definition}{Definition}[section]

\begin{document}

\title{{The Dirichlet problem in Lipschitz domains for higher order 
elliptic systems with rough coefficients}
\thanks{2000 {\it Math Subject Classification.} Primary: 35G15, 42B20, 
35J55, 35J40. Secondary 35J67, 35E05, 46E39.
\newline
{\it Key words}: higher order elliptic systems, Besov spaces, 
weighted Sobolev spaces, mean oscillations, BMO, VMO, Lipschitz domains, 
Dirichlet problem, traces, extension operators
\newline
${}$\qquad The work of authors was supported in part by from NSF, DMS and 
FRG grants as well as from the Swedish National Science Research Council.}} 
 
\author{V.\, Maz'ya, M.\, Mitrea and T.\, Shaposhnikova}

\date{~}

\maketitle

\begin{abstract}
We study the Dirichlet problem in Lipschitz domains and with boundary data 
in Besov spaces, for divergence form strongly elliptic systems of arbitrary 
order, with bounded, complex-valued coefficients. Our main result gives a 
sharp condition on the local mean oscillations of the coefficients of the 
differential operator and the unit normal to the boundary (which is 
automatically satisfied if these functions belong to space {\rm VMO}) 
guaranteeing that the solution operator associated with this problem 
is an isomorphism. 
\end{abstract}

\section{Introduction}
\setcounter{equation}{0}

\subsection{Formulation of the main result} 

A fundamental theme in the theory of partial differential equations, which 
has profound and intriguing connections with many other subareas of analysis, 
is the well-posedness of various classes of boundary value problems under 
minimal smoothness assumptions on the boundary of the domain and on the 
coefficients of the corresponding differential operator. 
The main result of this paper is 
{\it the solution of the Dirichlet problem for higher order, 
strongly elliptic systems in divergence form, with complex-valued, bounded, 
measurable coefficients in Lipschitz domains, and for boundary data in 
Besov spaces, under sharp smoothness assumptions}. In order to be more 
specific we need to introduce some notation. 

Let $m,l\in\NN$ be two fixed integers and, for a bounded Lipschitz domain
$\Omega$ in $\mathbb{R}^n$, 
$n\geq 2$, with outward unit normal $\nu=(\nu_1,...,\nu_n)$ consider 
the Dirichlet problem for the operator 
\begin{equation}\label{LOL}
{\mathcal L}(X,D_X)\,{\mathcal U}
:=\sum_{|\alpha|=|\beta|=m}D^\alpha(A_{\alpha\beta}(X)D^\beta{\mathcal U}),
\end{equation}
\noindent i.e., 
\begin{equation}\label{e0}
\left\{
\begin{array}{l}
\displaystyle{\sum_{|\alpha|=|\beta|=m}
D^\alpha(A_{\alpha\beta}(X)\,D^\beta\,{\mathcal U})}=0
\qquad\mbox{for}\,\,X\in\Omega,
\\[28pt] 
{\displaystyle\frac{\partial^k{\mathcal U}}{\partial\nu^k}}=g_k
\,\,\quad\mbox{on}\,\,\partial\Omega,\qquad 0\leq k\leq m-1.
\end{array}
\right.
\end{equation}
 \noindent Here and elsewhere, 
$D^\alpha=(-i\partial/\partial x_1)^{\alpha_1}\cdots
(-i\partial/\partial x_n)^{\alpha_n}$ if $\alpha=(\alpha_1,...,\alpha_n)$. 
The coefficients $A_{\alpha\beta}$ are $l\times l$ matrix-valued functions 
with measurable, complex-valued entries, for which there exists some finite 
constant $\kappa>0$ (referred to in the sequel as the ellipticity constant 
of ${\mathcal L}$) such that 
\begin{equation}\label{A-bdd}
\sum_{|\alpha|=|\beta|=m}\|A_{\alpha\beta}\|_{L_\infty(\Omega)}\leq\kappa^{-1}
\end{equation}
\noindent and such that the coercivity condition  
\begin{equation}\label{coercive}
\Re\,\int_\Omega\sum_{|\alpha|=|\beta|=m}\langle A_{\alpha\beta}(X) 
D^\beta\,{\mathcal V}(X),\,D^\alpha\,{\mathcal V}(X)\rangle\,dX 
\geq\kappa\sum_{|\alpha|=m}\|D^\alpha\,{\mathcal V}\|^2_{L_2(\Omega)}
\end{equation}
\noindent holds for all $\CC^l$-valued 
functions ${\mathcal V}\in C^\infty_0(\Omega)$. Throughout the paper, 
$\Re\,z$ denotes the real part of $z\in\CC$ and $\langle\cdot,\cdot\rangle$ 
stands for the canonical inner product in $\mathbb{C}^l$. 
Since, generally speaking, $\nu$ is merely bounded and measurable,
care should be exercised when defining iterated normal derivatives. For the 
setting we have in mind it is natural to take
$\partial^k/\partial\nu^k:=(\sum_{j=1}^n\xi_j\partial/\partial x_j)^k
\mid_{\xi=\nu}$ or, more precisely, 
\begin{equation}\label{nuk}
\frac{\partial^k{\mathcal U}}{\partial\nu^k}
:=i^k\sum_{|\alpha|=k}\frac{k!}{\alpha!}\,
\nu^\alpha\,{\rm Tr}\,[D^\alpha{\mathcal U}],\qquad 0\leq k\leq m-1,
\end{equation}
\noindent where ${\rm Tr}$ is the boundary trace operator and
$\nu^\alpha:=\nu_1^{\alpha_1}\cdots\nu_n^{\alpha_n}$ if 
$\alpha=(\alpha_1,...,\alpha_n)$.
Now, if $p\in(1,\infty)$, $a\in(-1/p,1-1/p)$ are fixed and 
$\rho(X):={\rm dist}\,(X,\partial\Omega)$, a solution for (\ref{e0}) is 
sought in $W^{m,a}_p(\Omega)$, defined as the space of vector-valued 
functions for which 
\begin{equation}\label{W-Nr}
\Bigl(\sum_{|\alpha|\leq m}\int_\Omega|D^\alpha{\mathcal U}(X)|^p 
\rho(X)^{pa}\,dX\Bigr)^{1/p}<\infty.
\end{equation}
\noindent In particular, as explained later on, the traces in (\ref{nuk}) 
exist in the Besov space $B_p^{s}(\partial\Omega)$, where $s:=1-a-1/p\in(0,1)$,
for any ${\mathcal U}\in W^{m,a}_p(\Omega)$. 
Recall that, with $d\sigma$ denoting the area element on $\partial\Omega$,  
\begin{equation}\label{Bes-xxx}
f\in B_p^s(\partial\Omega)\Leftrightarrow 
\|f\|_{B_p^s(\partial\Omega)}:=\|f\|_{L_p(\partial\Omega)}
+\Bigl(\int_{\partial\Omega}\int_{\partial\Omega}
\frac{|f(X)-f(Y)|^p}{|X-Y|^{n-1+sp}}\,d\sigma_Xd\sigma_Y\Bigr)^{1/p}<\infty.
\end{equation}

The above definition takes advantage of the Lipschitz manifold structure of 
$\partial\Omega$. On such manifolds, smoothness spaces of index $s\in(0,1)$ 
can be defined in an intrinsic, invariant fashion by lifting their Euclidean 
counterparts onto the manifold itself via local charts. Nonetheless, the very 
nature of the problem investigated in this paper requires the consideration of
{\it higher order smoothness spaces on $\partial\Omega$}, in which case the 
above approach is no longer effective. An alternative point of view has been 
introduced by H.\,Whitney in {\bf\cite{Wh}} where he considered higher order 
Lipschitz spaces on arbitrary closed sets (see also C.\,Fefferman's
article {\bf{\cite{Fe}}} for related issues). An extension of this 
circle of ideas pertaining to the full scale of Besov and Sobolev spaces on 
irregular subsets of $\RR^n$ can be found in the book {\bf\cite{JW}} by 
A.\,Jonsson and H.\,Wallin. Here we further refine this theory in the context 
of Lipschitz domains. For the purpose of this introduction we note that one 
possible description of these higher order Besov spaces on the boundary of 
a Lipschitz domain $\Omega\subset\RR^n$, and for $m\in\NN$, $p\in(1,\infty)$, 
$s\in(0,1)$, reads
\begin{equation}\label{Bes-X}
\dot{B}^{m-1+s}_p(\partial\Omega)=\,\mbox{the closure of}\,\,
\Bigl\{(i^{|\alpha|}D^\alpha\,{\mathcal V}|_{\partial\Omega})
_{|\alpha|\leq m-1}:\,
{\mathcal V}\in C^\infty_0(\RR^n)\Bigr\}\mbox{ in }
B_p^s(\partial\Omega)
\end{equation}
\noindent (making no notational distinction between a Banach space 
${\mathfrak X}$ and 
${\mathfrak X}^N={\mathfrak X}\oplus\cdots\oplus{\mathfrak X}$).
A formal definition, which involves higher order Taylor remainders in place of
$f(X)-f(Y)$ in (\ref{Bes-xxx}), along with other equivalent characterizations 
of $\dot{B}^{m-1+s}_p(\partial\Omega)$ can be found in \S{7.1}. 
Given (\ref{nuk})-(\ref{W-Nr}) and (\ref{Bes-X}), a necessary condition for 
the boundary data $\{g_k\}_{0\leq k\leq m-1}$ in (\ref{e0}) is that 
\begin{equation}\label{data-B}
g_k=\sum_{|\alpha|=k}\frac{k!}{\alpha!}\,\nu^\alpha\,f_\alpha,\quad
0\leq k\leq m-1,\quad\mbox{for some }\dot{f}=\{f_\alpha\}_{|\alpha|\leq m-1}
\in\dot{B}^{m-1+s}_p(\partial\Omega). 
\end{equation}

Let ${\rm BMO}$ and ${\rm VMO}$ stand, respectively, for the John-Nirenberg 
space of functions of bounded mean oscillations and the Sarason space 
of functions of vanishing mean oscillations (considered either on $\Omega$, 
or on $\partial\Omega$). Our main result, pertaining to the well-posedness 
of the problem (\ref{e0}), then reads as follows. 
\begin{theorem}\label{Theorem}
Let $\Omega\subset\RR^n$ be a bounded Lipschitz domain whose Lipschitz 
constant is $\leq M$, and assume that the operator ${\mathcal L}$ in 
{\rm (\ref{LOL})} satisfies {\rm (\ref{A-bdd})}-{\rm (\ref{coercive})}. 
Then there exists $c>0$, depending only on $M$ and the ellipticity constant 
of ${\mathcal L}$, with the following significance. 

Given $p\in (1,\infty)$ and $s\in (0,1)$ set $a:=1-s-1/p$. Then the Dirichlet 
problem {\rm (\ref{e0})} with boundary data as in {\rm (\ref{data-B})} 
has a unique solution ${\mathcal U}$ for which {\rm (\ref{W-Nr})} holds 
provided the coefficient matrices $A_{\alpha\beta}$ and the 
exterior normal vector $\nu$ to $\partial\Omega$ satisfy
\begin{equation}\label{axxx}
{\rm dist}\,(\nu,{\rm VMO}(\partial\Omega))
+\sum_{|\alpha|=|\beta|=m}{\rm dist}\,(A_{\alpha\beta},{\rm VMO}(\Omega))
\leq\,c\,s(1-s)\Bigl(p^2(p-1)^{-1}+s^{-1}(1-s)^{-1}\Bigr)^{-1}. 
\end{equation}
In particular, the above Dirichlet problem is well-posed 
for each $p\in (1,\infty)$ and $s\in (0,1)$ granted that 
$A_{\alpha\beta}\in {\rm VMO}(\Omega)$ and $\nu\in{\rm VMO}(\partial\Omega)$.
\end{theorem}
\noindent The above result (along with its inhomogeneous version, 
presented in \S{8.1}) is sharp. See \S{8.2} for a discussion. 
Note that $[\nu]_{{\rm BMO}(\partial\Omega)}
+\sum_{|\alpha|=|\beta|=m}[A_{\alpha\beta}]_{{\rm BMO}(\Omega)}
\leq\,c\,s(1-s)\bigl(p^2(p-1)^{-1}+s^{-1}(1-s)^{-1}\bigr)^{-1}$
suffices for (\ref{axxx}) to hold.
Other corollaries of independent interest are presented below.

\subsection{Some consequences of the main result and of its proof}

In the proof of Theorem~\ref{Theorem} we shall actually work with a related
condition in place of (\ref{axxx}), which quantifies the local 
oscillations of the coefficient matrices and the unit normal. 
To state this formally, we need one final piece of terminology. 
By the {\it BMO mod VMO character} 
of a function $F\in L_1(\Omega)$ we shall understand the quantity
\begin{equation}\label{e60}
\{F\}_{\ast,\Omega}:=\mathop{\hbox{lim\,sup}}_{\varepsilon\to 0} 
\left(\mathop{\hbox{sup}}_{(B_\varepsilon)_\Omega}
\meanint_{\!\!\!B_\varepsilon\cap\Omega}\,\,
\meanint_{\!\!\!B_\varepsilon\cap\Omega}\,
\Bigl|\,F(X)-F(Y)\,\Bigr|\,dXdY\right), 
\end{equation}
\noindent where $(B_\varepsilon)_\Omega$ stands for the collection of 
arbitrary balls centered at points of $\Omega$ and of radius $\varepsilon$, 
and the barred integral is the mean value. In a similar fashion, if 
$(B_\varepsilon)_{\partial\Omega}$ is the collection of $n$-dimensional 
balls with centers on $\partial\Omega$ and of radius $\varepsilon$, 
and if $f\in L_1(\partial\Omega)$, we set 
\begin{equation}\label{e61}
\{f\}_{\ast,\partial\Omega}
:=\mathop{\hbox{lim\,sup}}_{\varepsilon\to 0}
\left(\mathop{\hbox{sup}}_{(B_\varepsilon)_{\partial\Omega}}
\meanint_{\!\!\!B_\varepsilon\cap\partial\Omega}\,\,
\meanint_{\!\!\!B_\varepsilon\cap\partial\Omega}\,
\Bigl|\,f(X)-f(Y)\,\Bigr|\,d\sigma_Xd\sigma_Y\right), 
\end{equation}
For an arbitrary function $F$ we obviously have
(with the dependence on the domain dropped) 
$\{F\}_{\ast}\leq 2\,{\rm dist}\,(F,{\rm VMO})$, where the 
distance is taken in ${\rm BMO}$. Moreover, as a consequence of a result due
to D.\,Sarason (cf. Lemma~2 on p.\,393 of {\bf{\cite{Sar}}}), there 
exists $C>0$ such that ${\rm dist}\,(F,{\rm VMO})\leq C\{F\}_{\ast}$.
Thus, all together, $\{F\}_{\ast}\sim{\rm dist}\,(F,{\rm VMO})$ so that 
condition (\ref{axxx}) becomes equivalent to 
\begin{equation}\label{a0}
\{\nu\}_{\ast,\partial\Omega}
+\sum_{|\alpha|=|\beta|=m}\{ A_{\alpha\beta}\}_{\ast,\Omega}
\leq\,c\,s(1-s)\Bigl(p^2(p-1)^{-1}+s^{-1}(1-s)^{-1}\Bigr)^{-1}.
\end{equation}

A corollary of our main result is that, under the hypotheses of 
Theorem~\ref{Theorem}, the problem 
\begin{equation}\label{e0-bis}
\left\{
\begin{array}{l}
\displaystyle{\sum_{|\alpha|=|\beta|=m}
D^\alpha(A_{\alpha\beta}(X)\,D^\beta\,{\mathcal U})}=0
\qquad\mbox{in}\,\,\,\Omega,
\\[24pt] 
i^{|\gamma|}\,{\rm Tr}\,[D^\gamma\,{\mathcal U}]=f_\gamma
\,\,\quad\mbox{on}\,\,\partial\Omega,\qquad |\gamma|\leq m-1, 
\end{array}
\right.
\end{equation}
\noindent has a unique solution ${\mathcal U}$ which satisfies (\ref{W-Nr}) 
whenever 
\begin{equation}\label{data-G}
\dot{f}:=\{f_\gamma\}_{|\gamma|\leq m-1}\in\dot{B}^{m-1+s}_p(\partial\Omega).
\end{equation}
\noindent However, an advantage of the classical formulation 
(\ref{e0}) over (\ref{e0-bis}) is that, in the former case, the number of 
boundary conditions is {\it minimal}. For a domain $\Omega\subset\RR^2$ 
with boundary of class $C^{1+r}$, ${\textstyle\frac 1{2}}<r<1$, and for real,
constant coefficient, scalar operators, the limiting case $p=\infty$ of the 
Dirichlet problem (\ref{e0-bis}) has been considered by S.\,Agmon 
in {\bf\cite{Ag1}}. Exploiting the special nature of the layer potentials
associated with the equation in the two-dimensional setting, he has proved 
that there exists a unique solution 
${\mathcal U}\in C^{m-1+s}(\overline{\Omega})$, $0<s<r$, whenever 
$f_\gamma=i^{|\gamma|}D^\gamma{\mathcal V}|_{\partial\Omega}$, 
$|\gamma|\leq m-1$, for some function 
${\mathcal V}\in C^{m-1+s}(\overline{\Omega})$.
See also {\bf\cite{Ag2}} for a related problem. 

The innovation that allows us to consider boundary data in Besov spaces 
as in (\ref{data-G}) is the systematic use of {\it weighted Sobolev spaces} 
such as those associated with the norm in (\ref{W-Nr}). In relation to the 
standard Besov scale in $\RR^n$, we would like to point out that, thanks to 
Theorem~4.1 in {\bf\cite{JK}} on the one hand, and Theorem~1.4.2.4 and 
Theorem~1.4.4.4 in {\bf\cite{Gr}} on the other, we have
\begin{equation}\label{incls}
\begin{array}{l}
a=1-s-\frac{1}{p}\in (0,1-1/p)\Longrightarrow
W^{m,a}_p(\Omega)\hookrightarrow B^{m-1+s+1/p}_p(\Omega),
\\[10pt]
a=1-s-\frac{1}{p}\in (-1/p,0)\Longrightarrow
B^{m-1+s+1/p}_p(\Omega)\hookrightarrow W^{m,a}_p(\Omega).
\end{array}
\end{equation}
\noindent Of course, $W^{m,a}_p(\Omega)$ is just the classical Sobolev 
space $W^{m}_p(\Omega)$ when $a=0$. 

Remarkably, the classical trace theory for ordinary Sobolev spaces
in domains with smooth boundaries turns out to have a most satisfactory 
analogue in this weighted context and for Lipschitz domains. One of our 
main results in this regard is identifying the correct class of boundary data
for higher order Dirichlet problems for functions in $W^{m,a}_p(\Omega)$.
In the process, we establish that 
\begin{equation}\label{newtrace-ker}
\begin{array}{l}
{\mathcal U}\in W^{m,a}_p(\Omega)\mbox{ and }
\frac{\partial^k{\mathcal U}}{\partial\nu^k}=0\mbox{ on }\partial\Omega
\mbox{ for }0\leq k\leq m-1
\\[6pt]
\Longleftrightarrow 
{\mathcal U}\mbox{ belongs to the closure of $C^\infty_0(\Omega)$ in }
W^{m,a}_p(\Omega),
\end{array}
\end{equation}
\noindent which provides an answer to the question raised by J.\,Ne\v{c}as 
in Problem~4.1 on page 91 of his 1967 book {\bf{\cite{Nec}}}. In the context 
of unweighted Sobolev spaces and for smoother domains, such a result has been 
known for a long time (cf., e.g., P.\,Grisvard, S.M.\,Nikol'ski\u{\i} and 
H.\,Triebel's monographs {\bf{\cite{Gr}}}, {\bf{\cite{Ni}}}, {\bf{\cite{Tr}}} 
and the references therein). 

As a consequence of the trace theory developed in \S{7}
we have that, in the context of Theorem~\ref{Theorem}, 
\begin{equation}\label{Trace}
\sum_{|\alpha|\leq m-1}
\|{\rm Tr}\,[D^\alpha\,{\mathcal U}]\|_{B_p^{s}(\partial\Omega)}
\sim \left(\sum_{|\alpha|\leq m}\int_{\Omega}\rho(X)^{p(1-s)-1}\,
|D^\alpha{\mathcal U}(X)|^p\,dX\right)^{1/p},
\end{equation}
\noindent uniformly in ${\mathcal U}$ satisfying 
${\mathcal L}(X,D_X)\,{\mathcal U}=0$ in $\Omega$. The estimate (\ref{Trace}) 
can be viewed as a significant generalization of a well-known 
characterization of the membership of a function to a Besov space in 
$\RR^{n-1}$ in terms of weighted Sobolev norm estimates for 
its harmonic extension to $\RR^n_+$ (see, e.g., Proposition~$7'$ 
on p.\,151 of E.\,Stein's book {\bf\cite{St}}).

\subsection{A brief overview of related work}

Broadly speaking, there are two types of questions pertaining to the
well-posedness of the Dirichlet problem in a Lipschitz domain 
$\Omega$ for a divergence form, strongly elliptic system (\ref{LOL}) 
of order $2m$ with boundary data in Besov spaces indexed by $s$ and $p$. 
\vskip 0.08in
\noindent {\it Question I.} Granted that the coefficients of ${\mathcal L}$ 
exhibit a certain amount of smoothness, identify the indices $p$, $s$  
for which this boundary value problem is well-posed. 
\vskip 0.08in
\noindent {\it Question II.} Alternatively, having fixed the indices $s$ 
and $p$, characterize the smoothness of $\partial\Omega$ and 
of the coefficients of ${\mathcal L}$ for which the aforementioned problem 
is well-posed. 
\vskip 0.08in 
\noindent These, as well as other related issues, have been a driving 
force behind many exciting, recent developments in partial differential 
equations and allied fields. An authoritative account of their impact is given 
by C.\,Kenig in {\bf\cite{Ke}} where he describes the state of 
the art in this field of research up to mid 1990's. 
One generic problem which falls under the scope of {\it Question I} is to 
determine the optimal scale of spaces on which the Dirichlet problem for a
strongly elliptic system of order $2m$ is solvable in an 
{\it arbitrary Lipschitz domain} $\Omega$ in $\RR^n$. The most basic case, 
that of the constant coefficient Laplacian in arbitrary Lipschitz domains in 
$\RR^n$, is now well-understood thanks to the work of B.\,Dahlberg and 
C.\,Kenig {\bf\cite{DK}}, in the case of $L_p$-data, and D.\,Jerison and 
C.\,Kenig {\bf\cite{JK}}, in the case of Besov data. The case of (\ref{LapJK})
for boundary data exhibiting higher regularity (i.e., $s>1$) has been recently
dealt with by V.\,Maz'ya and T.\,Shaposhnikova in {\bf\cite{MS2}} where 
optimal smoothness conditions for $\partial\Omega$ are found in terms of the 
properties of $\nu$ as a Sobolev space multiplier. Generalizations of 
(\ref{LapJK}) to the case of variable-coefficient, single, second order 
elliptic equations have been obtained in 
{\bf\cite{MT1}}, {\bf\cite{MT2}}, {\bf\cite{MT3}}. 

In spite of substantial progress in recent years, there remain many basic 
open questions, particularly for $l>1$ and/or $m>1$, even in the case of 
{\it constant coefficient} operators in Lipschitz domains. In this context, 
one significant problem is to determine the sharp range of $p$'s for which 
the Dirichlet problem for strongly elliptic systems with $L_p$-boundary data 
is well-posed. In {\bf\cite{PV}}, 
J.\,Pipher and G.\,Verchota have developed a $L_p$-theory for real, constant 
coefficient, higher order systems $L=\sum_{|\alpha|=2m}A_\alpha D^\alpha$ when 
$p$ is near $2$, i.e., $2-\varepsilon<p<2+\varepsilon$ with $\varepsilon>0$ 
depending on the Lipschitz character of $\Omega$. On p.\,2 of {\bf\cite{PV}} 
the authors ask whether the $L_p$-Dirichlet problem for these operators is 
solvable in a given Lipschitz domain for 
$p\in (2-\varepsilon,\frac{2(n-1)}{n-3}+\varepsilon)$, and a positive answer 
has been recently given by Z.\,Shen in {\bf\cite{Sh}}. Let us also mention 
here the work {\bf\cite{AP}} of V.\,Adolfsson and J.\,Pipher who have dealt 
with the Dirichlet problem for the biharmonic operator in arbitrary Lipschitz 
domains and with data in Besov spaces, {\bf\cite{Ve}} where G.\,Verchota 
formulates and solves a Neumann-type problem for the bi-Laplacian in Lipschitz
domains and with boundary data in $L_2$, {\bf\cite{MMT}} where the authors
treat the Dirichlet problem for strongly elliptic systems 
of second order in an arbitrary Lipschitz subdomain $\Omega$ of a (smooth) 
Riemannian manifold and with boundary data in $B^s_p(\partial\Omega)$, 
when $2-\varepsilon<p<2+\varepsilon$ and $0<s<1$, as well as the paper 
{\bf\cite{KM}} by V.\,Kozlov and V.\,Maz'ya, which contains an explicit 
description of the asymptotic behavior of null-solutions of constant 
coefficient, higher order, elliptic operators near points on the boundary 
of a domain with a sufficiently small Lipschitz constant. 

A successful strategy for dealing with {\it Question II} consists of  
formulating and solving the analogue of the original problem in a 
standard case, typically when $\Omega=\RR^n_+$ and ${\mathcal L}$ has 
constant coefficients, and then deviating from this most standard setting 
by allowing perturbations of a certain magnitude. A paradigm result in this 
regard, going back to the work of S.\,Agmon, A.\,Douglis, L.\,Nirenberg and 
V.A.\,Solonnikov in the 50's and 60's, is that the Dirichlet problem 
is solvable in the context of Sobolev-Besov spaces if $\partial\Omega$ 
is sufficiently smooth and if ${\mathcal L}$ has continuous coefficients. 
The latter requirement is an artifact of the method of proof (based 
on Korn's trick of freezing the coefficients) which requires measuring 
the size of the oscillations of the coefficients in a {\it pointwise sense} 
(as opposed to integral sense, as in (\ref{e60})). 
For a version of {\it Question II}, corresponding to boundary data 
selected from $\prod_{k=0}^{m-1}B^{m-1-k+s}_p(\partial\Omega)$, 
optimal results have been obtained by V.\,Maz'ya and 
T.\,Shaposhnikova in {\bf\cite{MS}}. In this context, the natural language 
for describing the smoothness of the domain $\Omega$ is that of Sobolev 
space multipliers. 

In the smooth context, problems such as (\ref{e0}) 
have been investigated by many authors, including 
L.\,G{\aa}rding {\bf{\cite{Gar}}}, M.I.\,Vi\v{s}ik {\bf{\cite{Vi}}}, 
F.E.\,Browder {\bf{\cite{Br}}}, S.\,Agmon, A.\,Douglis and L.\,Nirenberg 
{\bf\cite{Ag1}}, {\bf\cite{Ag2}}, {\bf\cite{ADN}}, V.\,Solonnikov 
{\bf\cite{Sol1}}, {\bf\cite{Sol2}}, L.\,H\"ormander {\bf{\cite{Ho}}}, 
G.\,Grubb and N.J.\,Kokholm {\bf{\cite{GK}}}. A related result is as follows.
If $\Omega\subset{\mathbb{R}}^n$ is an {\it arbitrary bounded open set}, and
$g\in W^m_2(\Omega)$, $m\in{\mathbb{N}}$, is given, then the problem 
\begin{eqnarray}\label{Gen-Dir1}
\Delta^m u=0\,\,\mbox{ in }\,\,\Omega,\quad
u\in W^m_2(\Omega),\quad D^{\alpha}(u-g)|_{\partial\Omega}=0
\,\,\mbox{ for }\,\,|\alpha|\leq m-1,
\end{eqnarray}
\noindent where the boundary traces are taken in a generalized sense, 
has a unique solution. Building on some earlier work of K.\,Friedrichs, 
S.L.\,Sobolev has considered this problem in {\bf{\cite{Sob}}} in the case 
when $\partial\Omega$ consists of a finite union of submanifolds of 
${\mathbb{R}}^n$ of arbitrary codimension.  
This result also appears in Sobolev's 1950 monograph {\bf{\cite{Sob2}}}.
For arbitrary domains, the well-posedness of (\ref{Gen-Dir1}) has been
established by L.I.\,Hedberg in {\bf{\cite{Hed-1}}}, {\bf{\cite{Hed-2}}}.
The issue of continuity of the variational solution for higher-order equations
at boundary points has been studied by V.\,Maz'ya in {\bf{\cite{Mazz-4}}}.

\subsection{Comments on the proof of Theorem~\ref{Theorem} and the 
layout of the paper}

While the study of boundary value problems for elliptic differential 
operators with rough coefficients goes a long way back (it suffices to point 
to the connections with Hilbert's 19-th problem and De Giorgi-Nash-Moser 
theory), a lot of attention has been devoted lately to the class of 
operators with coefficients in ${\rm VMO}$. Part of the impetus for the 
recent surge of interest in this particular line of work stems from a key 
observation made by F.\,Chiarenza, M.\,Frasca and P.\,Longo in the 
early 1990's. More specifically, while investigating interior estimates for 
the solution of a scalar, second-order elliptic differential equation of the 
form ${\mathcal L}\,{\mathcal U}=F$, these authors have noticed in that 
${\mathcal U}$ can be related to $F$ via a potential theoretic representation 
formula in which the residual terms are commutators between operators of 
Mikhlin-Calder\'on-Zygmund type, on the one 
hand, and operators of multiplication by the coefficients of ${\mathcal L}$, 
on the other hand. This made it possible to control these terms by invoking 
the commutator estimate of Coifman-Rochberg-Weiss ({\bf\cite{CRW}}).  
An alternative method, based on maximal operators and good-$\lambda$ 
inequalities, has been developed by L.\,Caffarelli and I.\,Peral in 
{\bf\cite{CaPe}}, whereas when $\Omega=\RR^n$, an approach based 
on estimates for the Riesz transforms has been devised by T.\,Iwaniec 
and C.\,Sbordone in {\bf\cite{IS}}. Further related results can be found 
in {\bf\cite{AQ}}, {\bf\cite{By1}}, {\bf\cite{CFL1}}, {\bf\cite{Faz}}, 
{\bf\cite{MPS}}. 

Compared to the aforementioned works, our approach is more akin to that  
of F.\,Chiarenza and collaborators {\bf\cite{CFL1}} 
though there are fundamental differences between solving boundary problems 
for higher order and for second order operators. 
One difficulty inherently linked with the case $m>1$ arises from the 
way the norm in (\ref{W-Nr}) behaves under a change of variables
$\varkappa:\Omega=\{(X',X_n):\,X_n>\varphi(X')\}\to\RR^n_+$ designed to 
flatten the Lipschitz surface $\partial\Omega$. When $m=1$, a simple 
bi-Lipschitz changes of variables such as 
$\Omega\ni (X',X_n)\mapsto(X',X_n-\varphi(X'))\in\RR^n_+$ will do, but 
matters are considerable more subtle in the case $m>1$. In this latter 
situation, we employ a special global flattening map first introduced by 
J.\,Ne\v{c}as (in a different context; cf. p.\,188 in {\bf\cite{Nec}}) 
and then independently rediscovered and/or further adapted to new 
 settings by several authors, including V.\,Maz'ya and 
T.\,Shaposhnikova in {\bf\cite{MS}}, 
B.\,Dahlberg, C.\,Kenig J.\,Pipher, E.\,Stein and G.\,Verchota 
(cf. {\bf\cite{Dah}} and the discussion in {\bf\cite{DKPV}}), 
and S.\,Hofmann and J.\,Lewis in {\bf\cite{HL}}. Our main novel contribution 
in this regard is adapting this circle of ideas to the context when one 
seeks pointwise estimates for higher order derivatives of $\varkappa$ 
and $\lambda:=\varkappa^{-1}$ in terms of 
$[\nabla\varphi]_{{\rm BMO}(\RR^{n-1})}$. 

Another ingredient of independent interest is deriving estimates for 
$D_x^\alpha D_y^\beta G(x,y)$ where $G$ is the Green function associated 
with a constant (complex) coefficient system $L(D)$ of order $2m$ in the 
upper half space, which are sufficiently well-suited for deriving commutator 
estimates in the spirit of {\bf\cite{CRW}}. The methods employed in earlier 
works are largely based on explicit representation formulas for $G(x,y)$ and, 
hence, cannot be adapted easily to the case of general, non-symmetric, complex 
coefficient, higher order systems. By way of contrast, our approach 
consists of proving directly that the residual part 
$R(x,y):=G(x,y)-\Phi(x-y)$, where $\Phi$ is a fundamental solution 
for $L(D)$, has the property that $D_x^\alpha D_y^\beta R(x,y)$ 
is a Hardy-type kernel whenever $|\alpha|=|\beta|=m$. See also
{\bf{\cite{AQ}}} for a discussion of the difficulties encountered when 
estimating the residual part $R(x,y)$ in the case when ${\mathcal L}$ 
is not necessarily symmetric and has complex coefficients. 

The layout of the paper is as follows. Section~2 contains estimates 
for the Green function in the upper-half space. Section~3 deals with 
integral operators (of Mikhlin-Calder\'on-Zygmund and Hardy type) as well as
commutator estimates on weighted Lebesgue spaces. In the last part 
of this section we also revisit Gagliardo's extension operator and
establish estimates in the context of ${\rm BMO}$. 
Section~4 contains a discussion of the Dirichlet problem for higher 
order, variable coefficient, strongly elliptic systems in the upper-half space.
The adjustments necessary to treat the case of an unbounded domain 
lying above the graph of a Lipschitz function are presented in Section~5,
whereas in Section~6 we explain how to handle the case of a bounded
Lipschitz domain. In Section~7 we study traces and extension operators 
for higher order smoothness spaces on Lipschitz domains. Finally, in 
Section~8, we deal with the inhomogeneous version of (\ref{e0}); 
cf. Theorem~\ref{Theorem1} from which Theorem~\ref{Theorem} follows.

\section{Green's matrix estimates in the half-space}
\setcounter{equation}{0}

\subsection{Statement of the main result}
Fix two nonnegative integers $m,l$ and let ${L}(D_x)$ be a matrix-valued 
differential operator
\begin{equation}\label{eq1.1}
{L}(D_x)=\sum_{|\alpha|=2m}A_{\alpha} D^{\alpha}_x,
\end{equation}
\noindent
\noindent where the $A_\alpha$'s are constant $l\times l$ matrices with
complex entries. Throughout the paper, $D^\alpha_x:= i^{-|\alpha|}
\partial_{x_1}^{\alpha_1}\partial_{x_2}^{\alpha_2}\cdots
\partial_{x_n}^{\alpha_n}$ if
$\alpha=(\alpha_1,\alpha_2,...,\alpha_n)\in\NN_0^n$.
Here and elsewhere, $\NN$ stands for the collection of all
positive integers and $\NN_0:=\NN\cup\{0\}$. 
Assume that ${L}$ is strongly elliptic, i.e., there exists $\kappa>0$ such 
that $\sum_{|\alpha|=m}\|A_\alpha\|_{\CC^{l\times l}}\leq\kappa^{-1}$ and
\begin{equation}\label{eq1.2}
\Re\,\langle{L}(\xi)\eta,\eta\rangle_{\mathbb{C}^l}\geq\kappa\,
|\xi|^{2m}\,\|\eta\|^2_{\mathbb{C}^l},\qquad 
\forall\,\xi\in\RR^n,\,\,\,\forall\,\eta\in\mathbb{C}^l.
\end{equation}
\noindent In what follows, in order to simplify notations, we shall denote the 
norms in different finite-dimensional real Euclidean spaces by $|\cdot|$ 
irrespective of their dimensions. Also, quite frequently, we shall make 
no notational distinction between a space of scalar functions, call 
it ${\mathfrak X}$, and the space of vector-valued functions 
(of a fixed, finite dimension) whose components are in ${\mathfrak X}$.
We denote by $F(x)$ a  fundamental matrix of the  operator ${L}(D_x)$, i.e., 
an $l\times l$ matrix solution of the system
\begin{equation}\label{eq1.3}
{L}(D_x)F(x)=\delta(x)I_l\quad\mbox{in}\,\,\mathbb{R}^n,
\end{equation}
\noindent where $I_l$ is the $l\times l$ identity matrix and $\delta$ is the
Dirac function. We consider the Dirichlet problem
\begin{equation}\label{eq1.5}
\left\{
\begin{array}{l}
{L}(D_x)u=f\qquad\qquad\qquad \mbox{in}\,\,\mathbb{R}^n_+,
\\[6pt] 
{\rm Tr}\,[\partial^j u/\partial x_n^j]=f_j\quad j=0,1,\ldots,m-1,   
\qquad\quad\,\,\mbox{on}\,\,\mathbb{R}^{n-1},
\end{array}
\right.
\end{equation}
\noindent where 
$\mathbb{R}^n_+:=\{x=(x',x_n):\,x'\in\mathbb{R}^{n-1},\,x_n>0\}$
and ${\rm Tr}$ is the boundary trace operator. 
Hereafter, we shall identify $\partial\RR^n_+$ with $\RR^{n-1}$ in 
a canonical fashion.

For each $y'\in \mathbb{R}^{n-1}$ we introduce the Poisson matrices 
$P_0,\ldots, P_{m-1}$ for the problem (\ref{eq1.5}), i.e., the solutions of 
the boundary-value problems
\begin{equation}\label{eq1.6}
\left\{
\begin{array}{l}
{ L}(D_x)P_j(x,y')= 0\,I_l
\qquad\qquad\qquad\mbox{in}\,\,\mathbb{R}^n_+,
\\[10pt]
\displaystyle{\left(\frac{\partial^k}{\partial x_n^k}P_j\right)
(\,(x',0),y'\,)}
=\delta_{jk}\,\delta(x'-y')I_l\,\,\,{\rm for}\,\,\,x'\in\mathbb{R}^{n-1},\,\,
0\leq k\leq m-1,
\end{array}
\right.
\end{equation}
\noindent where $\delta_{jk}$ is the usual Kronecker symbol and 
$0\leq j\leq m-1$. The matrix-valued function $P_j(x,0')$ is positive 
homogeneous of degree $j+1-n$, i.e.,
\begin{equation}\label{eq1.7}
P_j(x,0')=|x|^{j+1-n}\,P_j(x/|x|,0'),\qquad x\in\mathbb{R}^n,
\end{equation}
\noindent where $0'$ denotes the origin of $\mathbb{R}^{n-1}$.
The restriction of $P_j(\cdot,0')$ to the upper half-sphere
$S^{n-1}_+$ is smooth and vanishes on the equator along with all of
its derivatives up to order $m-1$ (see for example, \S{10.3} in 
{\bf\cite{KMR2}}). Hence,
\begin{equation}\label{eq1.8}
\|P_j(x,0')\|_{\mathbb{C}^{l\times l}}
\leq C\,\frac{x_n^m}{|x|^{n+m-1-j}},\qquad x\in\RR^n_+,
\end{equation}
\noindent and, consequently,
\begin{equation}\label{eq1.9}
\|P_j(x,y')\|_{\mathbb{C}^{l\times l}}
\leq C\,\frac{x_n^m}{|x-(y',0)|^{n+m-1-j}},\qquad
x\in\RR^n_+,\,\,\,\,y'\in\RR^{n-1}.
\end{equation}

By $G(x,y)$ we shall denote the Green's matrix of the problem (\ref{eq1.5}), 
i.e., the unique solution of the boundary-value problem
\begin{equation}\label{eq1.10}
\left\{
\begin{array}{l}
{L}(D_x)G(x,y)=\delta(x-y)I_l\quad\mbox{for}\,\,x\in\mathbb{R}^n,
\\[6pt]
\displaystyle{\left(\frac{\partial ^j}{\partial x_n^j}G\right)((x',0),y)}
=0\,I_l\qquad
\mbox{for}\,\,x'\in\mathbb{R}^{n-1}, \,\,\, 0\leq j\leq m-1,
\end{array}
\right.
\end{equation}
\noindent where $y\in\mathbb{R}^n_+$ is regarded as a parameter.
We now introduce the matrix 
\begin{equation}\label{defRRR}
R(x,y):=F(x-y)-G(x,y),\qquad x,y\in\RR^n_+,
\end{equation}
\noindent so that, for each fixed $y\in\mathbb{R}^n_+$,
\begin{equation}\label{eq1.12}
\left\{
\begin{array}{l}
{L}(D_x)\,R(x,y)=0\qquad\qquad\qquad\qquad\qquad\quad\,\,\,\,
\mbox{for}\,\,x\in\mathbb{R}^n,
\\[6pt]
\displaystyle{\left(\frac{\partial^j}{\partial x_n^j}R\right)((x',0),y)}
=\left(\frac{\partial^j}{\partial x_n^j}F\right)((x',0)-y)
\quad\mbox{for}\,\,x'\in\mathbb{R}^{n-1},\,\, 0\leq j\leq m-1.
\end{array}
\right.
\end{equation}
\noindent Our goal is to prove the following result.
\begin{theorem}\label{th1}
For all multi-indices $\alpha,\beta$ of length $m$
\begin{equation}\label{mainest}
\|D^\alpha_xD^\beta_y R(x,y)\|_{\mathbb{C}^{l\times l}}
\leq C\,|x-\bar{y}|^{-n},
\end{equation}
\noindent for $x,y\in\RR^n_+$, where $\bar{y}:=(y',-y_n)$ is
the reflection of the point $y\in\RR^n_+$ with respect to $\partial\RR^n_+$.
\end{theorem}
In the proof of Theorem~\ref{th1} we distinguish two cases, $n>2m$ and 
$n\leq 2m$, which we shall treat separately. 
Our argument pertaining to
the situation when $n>2m$ is based on the following useful estimate 
for a parameter dependent integral. 
\begin{lemma}\label{lem1}
Let $a$ and $b$ be two non-negative numbers and assume that 
$\zeta\in \mathbb{R}^N$. Then for every $\varepsilon>0$ and $0<\delta<N$
there exists a constant $c(N,\varepsilon,\delta)>0$ such that 
\begin{equation}\label{E1}
\int_{\mathbb{R}^N}
\frac{d\eta}{(|\eta|+a)^{N+\varepsilon}(|\eta-\zeta|+ b)^{N-\delta}} 
\leq\frac{c(N,\varepsilon,\delta)}{a^\varepsilon(|\zeta|+a+b)^{N-\delta}}.
\end{equation}
\end{lemma}

\noindent The proof is postponed for \S{2.4}, as to prevent disrupting the
flow of the presentation.

\subsection{Proof of Theorem~\ref{th1} for $n>2m$}

In the case when $n>2m$ there exists a unique fundamental matrix $F(x)$ 
for the operator (\ref{eq1.1}) which is positive homogeneous of 
degree $2m-n$. We shall use the integral representation formula 
\begin{equation}\label{IntRRR}
R(x,y)=R_0(x,y)+\ldots+ R_{m-1}(x,y),\qquad x,y\in\RR^n_+,
\end{equation}
\noindent where $R(x,y)$ has been introduced in (\ref{defRRR}) and, 
with $P_j$ as in (\ref{eq1.6}), we set  
\begin{equation}\label{eq1.13}
R_j(x,y):=\int_{\mathbb{R}^{n-1}}P_j(x,\xi')\,
\left(\frac{\partial^j}{\partial x_n^j}F\right)((\xi',0)-y)\,d\xi',
\qquad 0\leq j\leq m-1.  
\end{equation}
\noindent Then, thanks to (\ref{eq1.8}) we have
\begin{equation}\label{eq1.14}
\|R_j(x,y)\|_{\mathbb{C}^{l\times l}}\leq C\,\int_{\mathbb{R}^{n-1}}
\frac{x_n^m}{|x-(\xi',0)|^{n+m-1-j}}\cdot\frac{d\xi'}{|(\xi',0)-y|^{n-2m+j}}.
\end{equation}
Next, using Lemma~\ref{lem1} with 
\begin{eqnarray}\label{n=n-1}
N=n-1 &, & \quad a=x_n\,\,,\quad \delta=2m-j-1,
\\
\varepsilon=m-j &, & \quad b=y_n\,\,\,,\quad \zeta=y'-x',
\nonumber
\end{eqnarray}
\noindent we obtain from (\ref{eq1.14})
\begin{equation}\label{Rj}
\|R_j(x,y)\|_{\mathbb{C}^{l\times l}}
\leq\frac{C\,x_n^j}{(|y'-x'|+x_n+y_n)^{n-2m+j}},\qquad 0\leq j\leq m-1.
\end{equation}
\noindent Summing up over $j=0,\ldots,m-1$ gives, by virtue of (\ref{IntRRR}), 
the estimate
\begin{equation}\label{eq1.16}
\|R(x,y)\|_{\mathbb{C}^{l\times l}}\leq C\,|x-{\bar y}|^{2m-n},
\qquad x,y\in\RR^n_+.
\end{equation}

To obtain pointwise estimates for derivatives of $R(x,y)$, we make 
use of the following local estimate for a solution of problem (\ref{eq1.5}) 
with $f=0$. Recall that $W^s_p$ stands for the classical $L_p$-based Sobolev 
space of order $s$. The particle {\it loc} is used to brand the local versions 
of these (and other) spaces.
\begin{lemma}\label{l1.1}{\rm [see {\bf\cite{ADN}}]}
Let $\zeta$ and $\zeta_0$ be functions in $C_0^\infty(\mathbb{R}^n)$ 
such that $\zeta_0 =1$ in a neighborhood of $\mbox{supp}\,\zeta$. Then the 
solution $u\in W_2^m(\mathbb{R}_+^n,loc)$ of problem {\rm (\ref{eq1.5})} 
with $f=0$ and $f_j\in W_p^{k+1-j-1/p}(\mathbb{R}^{n-1},loc)$, 
where $k\geq m$ and $p\in(1,\infty)$, belongs to 
$W_p^{k+1}(\mathbb{R}^n_+,loc)$ and satisfies the estimate
\begin{equation}\label{eq1.17}
\|\zeta u\|_{W_p^{k+1}(\mathbb{R}^n_+)}\leq C\,
\Bigl(\sum_{j=0}^{m-1}\|\zeta_0 f_j\|_{W_p^{k+1-j-1/p}(\mathbb{R}^{n-1})}
+\|\zeta_0 u\|_{L_p(\mathbb{R}^n_+)}\Bigr),
\end{equation}
\noindent where $C$ is independent of $u$ and $f_j$.
\end{lemma}
\noindent Let $B(x,r)$ denote the ball of radius $r>0$ centered at $x$.
\begin{corollary}\label{c1.2}
Assume that $u\in W_2^m(\mathbb{R}_+^n,loc)$ is a solution of problem 
{\rm (\ref{eq1.5})} with $f=0$ and 
$f_j\in C^{k+1-j}(\mathbb{R}^{n-1},loc)$. 
Then for any $z\in\overline{\mathbb{R}^n_+}$ and $\rho>0$,
\begin{equation}\label{eq1.18}
\sup_{\mathbb{R}^n_+\cap B(z,\rho)}|\nabla _k u|\leq C\,
\Bigl(\,\rho^{-k}\sup_{\mathbb{R}^n_+\cap B(z,2\rho)}|u|
+\sum_{j=0}^{m-1}\sum_{s=0}^{k+1-j}
\rho^{s+j-k}\sup_{\mathbb{R}^{n-1}\cap B(z,2\rho)}|\nabla'_{s}f_j|\Bigr),
\end{equation}
\noindent where $\nabla'_s$ is the gradient of order $s$ in $\mathbb{R}^{n-1}$.
Here $C$ is a  constant independent of $\rho$, $z$, $u$ and $f_j$.
\end{corollary}
\noindent{\bf Proof.} Given the dilation invariant nature of the estimate
we seek, it suffices to assume that $\rho=1$. 
Given $\phi\in C^{k+1-j}(\RR^{n-1})$ supported in $\RR^{n-1}\cap B(z,2)$,
we observe that, for a suitable $\theta\in (0,1)$,
\begin{equation}\label{simple}
\|\phi\|_{W_p^{k+1-j-1/p}(\mathbb{R}^{n-1})}
\leq C\|\phi\|_{L_p(\mathbb{R}^{n-1})}^\theta
\|\phi\|_{W_p^{k+1-j}(\mathbb{R}^{n-1})}^{1-\theta}
\leq C\sum_{s=0}^{k+1-j}\sup_{\mathbb{R}^{n-1}\cap B(z,2)}|\nabla'_{s}\phi|.
\end{equation}
\noindent Also, if $p>n$, 
\begin{equation}\label{eq1.19}
\sup\limits_{\mathbb{R}^n_+}|\nabla_k v|\leq C\,
\|v\|_{W_p^{k+1}(\mathbb{R}^n_+)},
\end{equation}
\noindent by Sobolev's inequality. Now, (\ref{eq1.18}) follows
by combining (\ref{simple}), (\ref{eq1.19}) with Lemma~\ref{l1.1}. 
\hfill$\Box$
\vskip 0.08in
Given $x,y\in\RR^n_+$, set $\rho:=|x-\bar{y}|/5$ and pick
$z\in\partial\RR^n_+$ such that $|x-z|=\rho/2$. It follows that for any
$w\in\RR^n_+\cap B(z,2\rho)$ we have
$|x-\bar{y}|\leq |x-z|+|z-w|+|w-\bar{y}|\leq \rho/2+2\rho+|w-\bar{y}|
\leq |x-\bar{y}|/2+|w-\bar{y}|$. Consequently,
$|x-\bar{y}|/2\leq |w-\bar{y}|$  for every $w\in\RR^n_+\cap B(z,2\rho)$,
so that, ultimately,
\begin{equation}\label{nablaPhi}
\rho^{\nu-k}\sup_{w\in\mathbb{R}^{n-1}\cap B(z,2\rho)}
\|\nabla'_{\nu}F(w-y)\|_{\CC^{l\times l}}
\leq \frac{C}{|x-\bar{y}|^{n-2m+k}},
\end{equation}
\noindent for each $\nu\in\NN_0$.
Granted (\ref{eq1.16}) and our choice of $\rho$, we altogether obtain that
\begin{equation}\label{eq1.20}
\|D^\alpha_x R(x,y)\|_{\mathbb{C}^{l\times l}}
\leq C_{k}\,|x-\bar{y}|^{2m-n-k},\quad x,y\in\RR^n_+,\quad\forall\,
\alpha\in\NN_0^n,\,\,|\alpha|\leq k.
\end{equation}
In the following two formulas, it will be convenient to use the notation 
$R_{\mathcal L}$ for the matrix $R$ associated with the operator 
${\mathcal L}(D_x)$ as in (\ref{defRRR}). By Green's formula
\begin{equation}\label{eq1.11}
R_{\mathcal L}(y,x)=\Bigl[R_{{\mathcal L}^*}(x,y)\Bigr]^*,\qquad x,y\in\RR^n_+,
\end{equation}
\noindent where the superscript star indicates adjunction.
In order to estimate {\it mixed} partial derivatives, we observe 
that (\ref{eq1.11}) entails
\begin{equation}\label{eq1.21}
(D^\beta_y R_{\mathcal L})(x,y)
=\Bigl[(D^\beta_x R_{{\mathcal L}^*})(y,x)\Bigr]^*
\end{equation}
\noindent and remark that ${\mathcal L}^*$ has properties similar 
to ${\mathcal L}$. This, in concert with (\ref{eq1.20}) and 
the fact that $|x-\bar{y}|=|\bar{x}-y|$ for $x, y\in\mathbb{R}^n_+$, yields
\begin{equation}\label{eq1.22}
\|D^\beta_y R(x,y)\|_{\mathbb{C}^{l\times l}}
\leq C_{\beta}\,|x-\bar{y}|^{2m-n-|\beta|}.
\end{equation}
\noindent Let us also point out that by formally differentiating
(\ref{eq1.12}) with respect to $y$ we obtain
\begin{equation}\label{eq1.23}
\left\{
\begin{array}{l}
{\mathcal L}(D_x)\,[D^\beta_yR_{\mathcal L}(x,y)]=0
\qquad\qquad\qquad\qquad\qquad 
\mbox{for}\,\,x\in\mathbb{R}^n,
\\[10pt]
\displaystyle{\left(\frac{\partial^j}{\partial x_n^j}
D^\beta_y R\right)((x',0),y)
=\left(\frac{\partial ^j}{\partial x_n^j}(-D)^\beta F\right)((x',0)-y)},
\,\,x'\in\mathbb{R}^{n-1},\,\,0\leq j\leq m-1.
\end{array}
\right.
\end{equation}
\noindent With (\ref{eq1.22}) and (\ref{eq1.23}) in place of (\ref{eq1.16}) 
and (\ref{eq1.12}), respectively, we can now run the same program as above 
and obtain the estimate 
\begin{equation}\label{Eq3}
\|D^\alpha_xD^\beta_y R(x,y)\|_{\mathbb{C}^{l\times l}}
\leq C_{\alpha\beta}\,|x-\bar{y}|^{2m-n-|\alpha|-|\beta|},\quad
\forall\,\alpha,\beta\in{\mathbb{N}}_0.
\end{equation}

\subsection{Proof of Theorem~\ref{th1} for $n\leq 2m$}

When $n\leq 2m$ we shall use the method of descent. To get started, fix an 
integer $N$ such that $N>2m$ and let $(x,z)\mapsto {\mathcal G}(x,y,z-\zeta)$ 
denote the Green matrix with singularity at $(y,\zeta)\in\RR^n\times\RR^{N-n}$
of the Dirichlet problem for the operator ${\mathcal L}(D_x)+(-\Delta_z)^m$ 
in the $N$-dimensional half-space 
\begin{equation}\label{RN}
\mathbb{R}^N_+:=
\{(x,z):\,z\in\mathbb{R}^{N-n},\,x=(x',x_n),\,x'\in\mathbb{R}^{n-1},\,x_n>0\}.
\end{equation}
\noindent Also, recall that $G(x,y)$ stands for the Green matrix of 
the problem (\ref{eq1.5}). 
\begin{lemma}\label{lem3}
For all multi-indices $\alpha$ and $\beta$ of order $m$ and for all 
$x$ and $y$ in $\mathbb{R}^n_+$
\begin{equation}\label{Eq12}
D^\alpha_x D^\beta_y G(x',y)
=\int_{\mathbb{R}^{N-n}}D^\alpha_x D^\beta_y{\mathcal G}(x,y,-\zeta)\,d\zeta.
\end{equation}
\end{lemma}
\noindent{\bf Proof.} The strategy is to show that 
\begin{equation}\label{pairGf}
\int_{\mathbb{R}^n_+}D^\alpha_xD_y^\beta G(x,y)\,f_\beta(y)\,dy
=\int_{\mathbb{R}^n_+}\int_{\mathbb{R}^{N-n}}D^\alpha_xD_y^\beta 
{\mathcal G}(x, y,-\zeta)\,d\zeta\,f_\beta(y)\,dy
\end{equation}
\noindent for each $f_\beta\in C^\infty_0(\mathbb{R}^n_+)$, from which 
(\ref{Eq12}) clearly follows. To justify (\ref{pairGf}) for a fixed,
arbitrary $f_\beta\in C^\infty_0(\mathbb{R}^n_+)$, we let $u$ be the unique 
vector-valued function satisfying $D^\alpha u\in L^2(\mathbb{R}^n_+)$ for 
all $\alpha$ with $|\alpha|=m$, and such that 
\begin{equation}\label{Eq9}
\left\{
\begin{array}{l}
{\mathcal L}(D_x)u=D^\beta_x f_\beta \qquad{\rm in}\,\,\mathbb{R}^n_+,
\\[6pt]
\displaystyle{\left(\frac{\partial^j u}{\partial x_n^j}\right)(x',0)=0
\qquad {\rm on}\,\,\mathbb{R}^{n-1}},\,\,0\leq j\leq m-1. 
\end{array}
\right.
\end{equation}
\noindent It is well-known that for each $\gamma\in\NN_0^n$ 
\begin{equation}\label{Eq10}
|D^\gamma u(x)|\leq C_\gamma\,|x|^{m-n-|\gamma|}\qquad{\rm for}\,\,|x|>1.
\end{equation}
\noindent This follows, for instance, from Theorem~6.1.4 {\bf\cite{KMR1}} 
combined with Theorem~10.3.2 {\bf\cite{KMR2}}. Also, as a consequence of 
Green's formula, the solution of the problem (\ref{Eq9}) satisfies
\begin{equation}\label{Eq14}
D^\alpha_x u(x)
=\int_{\mathbb{R}^n_+}D^\alpha_x(-D_y)^\beta G(x,y)\,f_\beta(y)\,dy.
\end{equation}

We shall now derive yet another integral representation formula for 
$D^\alpha_x u$ in terms of (derivatives of) ${\mathcal G}$ which is 
similar in spirit to (\ref{Eq14}). Since $N>2m$, (\ref{Eq3}) implies
\begin{equation}\label{Eq13}
\|D^\alpha_x D^\beta_y {\mathcal G}(x,y,-\zeta)\|_{\mathbb{C}^{l\times l}} 
\leq c\,(|x-y|+|\zeta|)^{-N}.
\end{equation}
\noindent Let us now fix $x\in\mathbb{R}^n_+$, $\rho>0$ and introduce 
a cut-off function $H\in C^\infty(\mathbb{R}^{N-n})$ which satisfies
$H(z)=1$ for $|z|\leq 1$ and $H(z)=0$ for $|z|\geq 2$. We may then write
\begin{equation}\label{u=G}
u(x)=\int_{\mathbb{R}^N} {\mathcal G}(x,y,-\zeta) 
\Bigl[H\bigl(\zeta/\rho\bigr)D^\beta f_\beta(y)+(-\Delta_\zeta)^m 
\bigl(H\bigl(\zeta/\rho\bigr)\,u(y)\bigr)\Bigr]\,dy\,d\zeta,
\end{equation}
\noindent which further implies
\begin{eqnarray}\label{est1}
&& \Bigl|D^\alpha_x u(x)-\int_{\mathbb{R}^N} D^\alpha_x(-D_y)^\beta\,
{\mathcal G}(x,y,-\zeta)\,H\bigl(\zeta/\rho\bigr)\,f_\beta(y)\,dy\,d\zeta\Bigr|
\nonumber\\[6pt]
&& \qquad\leq c\,\sum_{|\gamma|=m}\int_{\mathbb{R}^N_+}
\|D^\alpha_x D_\zeta^\gamma\,{\mathcal G}(x,y,-\zeta)\|
_{\mathbb{C}^{l\times l}}\,
\bigl|u(y)\,D^\gamma_\zeta\bigl(H\bigl(\zeta/\rho\bigr)\bigr)\bigr|\,d\zeta.
\end{eqnarray}
\noindent By (\ref{Eq10}) and (\ref{Eq13}), the expression in the 
right-hand side of (\ref{est1}) does not exceed
\begin{eqnarray*}
&& c\,\rho^{-m}\int_{\rho<|\zeta|<2\rho}d\zeta
\int_{\mathbb{R}^{n-1}}(|x-y|+|\zeta|)^{-N}\,|y|^{m-n}\,dy
\\[6pt]
&& \qquad\qquad
\leq c\,\rho^{N-n-m}\int_{\mathbb{R}^{n-1}}(|y|+\rho)^{-N}|y|^{m-n}\,dy
=c\,\rho^{-n}.
\end{eqnarray*}
\noindent This estimate, in concert with (\ref{Eq13}), allows us 
to obtain, after making $\rho\to\infty$, that 
\begin{equation}\label{DalphaU}
D^\alpha_x u(x)
=\int_{\mathbb{R}^n_+}\int_{\mathbb{R}^{N-n}}D^\alpha_x(-D_y)^\beta 
{\mathcal G}(x, y,-\zeta)\,d\zeta\,f_\beta(y)\,dy.
\end{equation}
\noindent Now (\ref{pairGf}) follows readily from this and (\ref{Eq14}). 
\hfill$\Box$
\vskip 0.08in

Having disposed of Lemma~\ref{lem3}, we are ready to discuss the

\vskip 0.08in
\noindent{\bf End of Proof of Theorem~\ref{th1}.} Assume that $2m\geq n$ 
and let $N$ be  again an integer such that $N>2m$. Denote by 
${\mathcal F}(x,z)$ the fundamental solution of the operator 
${\mathcal L}(D_x)+(-\Delta_z)^m$, which is positive homogeneous of 
degree $2m-N$ and is singular at $(0,0)\in\RR^n\times\RR^{N-n}$. 
Then the identity 
\begin{equation}\label{Eq17}
D^{\alpha+\beta}_x F(x)=\int_{\mathbb{R}^{N-n}}D^{\alpha+\beta}_x
{\mathcal F}(x,-\zeta)\,d\zeta
\end{equation}
\noindent can be established as in the proof of Lemma~\ref{lem3}.
Combining (\ref{Eq17}) with Lemma~\ref{lem3}, we arrive at
\begin{equation}\label{Eq18}
D^\alpha_x D^\beta_y R(x',y)=\int_{\mathbb{R}^{N-n}}D^\alpha_x D^\beta_y 
{\mathcal R}(x,y,-\zeta)\,d\zeta,
\end{equation}
\noindent where ${\mathcal R}(x,y,z):={\mathcal G}(x,y,z)-{\mathcal F}(x-y,z)$.
Consequently, 
\begin{equation}\label{DalDbetU}
\|D^\alpha_x D^\beta_y {\mathcal R}(x,y,-\zeta)\|_{\mathbb{C}^{l\times l}} 
\leq C(|x-\bar{y}|+|\zeta|)^{-N}.
\end{equation}
\noindent by (\ref{eq1.20})  with $k=0$ and $N$ in place of $n$. 
This estimate, together with (\ref{Eq18}), then yields (\ref{mainest}).
\hfill$\Box$

\subsection{Proof of Lemma~\ref{lem1}}

Write ${\mathcal J}={\mathcal J}_1+{\mathcal J}_2$ 
where ${\mathcal J}$ stands for the integral in the left side of (\ref{E1}), 
whereas  ${\mathcal J}_1$ and ${\mathcal J}_2$ denote the integrals
obtained by splitting the domain of integration in ${\mathcal J}$ into 
the ball $B_a=\{\eta\in \mathbb{R}^N:\,|\eta|<a\}$ and $\RR^n\setminus B_a$, 
respectively. If $|\zeta|<2a$, then
\begin{equation}\label{I1}
{\mathcal J}_1\leq a^{-N-\varepsilon}
\int_{B_a}\frac{d\eta}{(|\eta-\zeta|+b)^{N-\delta}}
\leq c\,a^{-N-\varepsilon}\int_{B_{4a}}\frac{d\xi}{(|\xi|+b)^{N-\delta}}.
\end{equation}

\noindent Hence
\begin{equation}\label{II1}
{\mathcal J}_1\leq 
\left\{
\begin{array}{l}
c\,a^{-N-\varepsilon}\,a^{N}/b^{N-\delta}\qquad{\rm if}\,\,a<b,
\\[6pt]
c\,a^{-N-\varepsilon+\delta}\qquad\quad{\rm if}\,\,a>b,
\end{array}
\right.
\end{equation}
\noindent so that, in particular, 
\begin{equation}\label{I1est}
|\zeta|<2a\Longrightarrow 
{\mathcal J}_1\leq c\,a^{-\varepsilon}(|\zeta|+a+b)^{\delta-N}.
\end{equation}

Let us now assume that $|\zeta|>2a$. Then
\begin{equation}\label{I1est2}
{\mathcal J}_1\leq\int_{B_a}\frac{d\eta}{(|\eta|+a)^{N+\varepsilon}}\, 
\frac{c}{(|\zeta|+b)^{N-\delta}}
\leq c\,a^{-\varepsilon}(|\zeta|+a+b)^{\delta-N},
\end{equation}
\noindent which is of the right order. As for ${\mathcal J}_2$, we write 
\begin{equation}\label{I2}
{\mathcal J}_2\leq\int_{\mathbb{R}^n\backslash B_a}
\frac{d\eta}{|\eta|^{N+\varepsilon}(|\eta-\zeta|+b)^{N-\delta}} 
={\mathcal J}_{2,1}+{\mathcal J}_{2,2}.
\end{equation}
\noindent where ${\mathcal J}_{2,1}$, ${\mathcal J}_{2,2}$ are obtained by 
splitting the domain of integration in the above integral into the set 
$\{\eta:|\eta|>\max\{a,2|\zeta|\}\}$ and its complement in 
$\mathbb{R}^n\backslash B_a$. We have
\begin{eqnarray}\label{Eq4}
{\mathcal J}_{2,1} & \leq & \int_{|\eta|>\max\{a,b,2|\zeta|\}}
\frac{d\eta}{|\eta|^{N+\varepsilon}(|\eta|+b)^{N-\delta}} 
+\int_{b>|\eta|>\max\{a,b,2|\zeta|\}}
\frac{d\eta}{|\eta|^{N+\varepsilon}(|\eta|+b)^{N-\delta}}
\nonumber\\[6pt]
& \leq & c\Biggl(\int_{|\eta|>\max\{a,b,2|\zeta|\}}
\frac{d\eta}{|\eta|^{2N+\varepsilon-\delta}}
+\frac{1}{b^{N-\delta}}\int_{b>|\eta|>\max\{a,b,2|\zeta|\}}
\frac{d\eta}{|\eta|^{N+\varepsilon}}\Biggr)
\nonumber\\[6pt]
& \leq & \frac{c}{(a+b+|\zeta|)^{N+\varepsilon-\delta}}
+\frac{c}{a^\varepsilon(a+b+|\zeta|)^{N-\delta}}
\nonumber\\[6pt]
& \leq & \frac{c}{a^\varepsilon(a+b+|\zeta|)^{N-\delta}}.
\end{eqnarray}

There remains to estimate the integral
\begin{equation}\label{Eq5}
{\mathcal J}_{2,2}=\int_{B_{2|\zeta|}\backslash B_a}
\frac{d\eta}{|\eta|^{N+\varepsilon}(|\eta-\zeta|+b)^{N-\delta}} 
={\mathcal J}_{2,2}^{(1)}+{\mathcal J}_{2,2}^{(2)},
\end{equation}
\noindent where ${\mathcal J}_{2,2}^{(1)}$ and ${\mathcal J}_{2,2}^{(2)}$ are 
obtained by splitting the domain of integration in ${\mathcal J}_{2,2}$ into 
$B_{|\zeta|/2}\backslash B_a$ and its complement (relative to 
$B_{2|\zeta|}\backslash B_a$). On the one hand, 
\begin{equation}\label{Eq6}
{\mathcal J}_{2,2}^{(1)}\leq\frac{c}{(|\zeta|+b)^{N-\delta}}
\int_{B_{|\zeta|/2}\backslash B_a}\frac{d\eta}{|\eta|^{N+\varepsilon}} 
\leq\frac{c}{a^\varepsilon(|\zeta|+a+b)^{N-\delta}}.
\end{equation}
\noindent On the other hand, whenever $|\zeta|>a/2$, the integral 
${\mathcal J}_{2,2}^{(2)}$, which extends over all $\eta$'s such 
that $|\eta|>a$, $2|\zeta|>|\eta|>|\zeta|/2$, can be estimated as
\begin{eqnarray}\label{cv-p}
{\mathcal J}_{2,2}^{(2)} & \leq & \frac{c}{|\zeta|^{N+\varepsilon}}
\int_{B_{2|\zeta|}\backslash B_a}\frac{d\eta}{(|\eta-\zeta|+b)^{N-\delta}}
\leq\frac{c}{|\zeta|^{N+\varepsilon}}
\int_{B_{4|\zeta|}}\frac{d\xi}{(|\xi|+b)^{N-\delta}}
\nonumber\\[6pt]
& \leq & \frac{c}{|\zeta|^{N+\varepsilon}}
\Biggl(\int_{{|\xi|<4|\zeta|}\atop{|\xi|<b}}
\frac{d\xi}{(|\xi|+b)^{N-\delta}}
+\int_{{|\xi|<4|\zeta|}\atop{|\xi|>b}}\frac{d\xi}{(|\xi|+b)^{N-\delta}}\Biggr).
\end{eqnarray}
\noindent Consequently, 
\begin{equation}\label{J22}
{\mathcal J}_{2,2}^{(2)}\leq c\,
\frac{\min\{|\zeta|,b\}^N}{|\zeta|^{N+\varepsilon} b^{N-\delta}}.
\end{equation}
\noindent Using $|\zeta|>a/2$ and the obvious inequality
\begin{equation}\label{trivial}
\min\{|\zeta|,b\}^N\cdot\max\{|\zeta|,b\}^{N-\delta}\leq|\zeta|^N\,b^{N-\delta}
\end{equation}
\noindent we arrive at
\begin{equation}\label{Eq7}
{\mathcal J}_{2,2}^{(2)}\leq c\,a^{-\varepsilon}(|\zeta|+a+b)^{\delta-N}.
\end{equation}
\noindent The estimate (\ref{Eq7}), along with (\ref{Eq6}) and (\ref{Eq5}), 
gives the upper bound $c\,a^{-\varepsilon}(|\zeta|+ a+b)^{\delta-N}$ for 
${\mathcal J}_{2,2}$. Combining this with (\ref{Eq4}) we obtain the same 
majorant for ${\mathcal J}_{2}$ which, together with a similar result 
for ${\mathcal J}_{1}$ already obtained, leads to (\ref{E1}). 
The proof of the lemma is therefore complete.

\section{Properties of integral operators in a half-space}
\setcounter{equation}{0}

\noindent In \S{3.1} and \S{3.2} we prove estimates for commutators 
(and certain commutator-like operators) between 
integral operators in $\mathbb{R}^n_+$ and multiplication operators with 
functions of bounded mean oscillations, in weighted Lebesgue spaces on 
$\mathbb{R}^n_+$. Subsection~3.3 contains ${\rm BMO}$ and pointwise 
estimates for extension operators from $\mathbb{R}^{n-1}$ onto 
$\mathbb{R}^n_+$. Throughout, given two Banach spaces $E,F$, we let 
${\mathfrak L}(E,F)$ stand for the space of bounded linear operators 
from $E$ into $F$, and abbreviate ${\mathfrak L}(E):={\mathfrak L}(E,E)$.
Also, given $p\in[1,\infty]$, an open set ${\mathcal O}\subset\RR^n$ 
and a measurable nonnegative function $w$ on ${\mathcal O}$, we let 
$L_p({\mathcal O},w(x)\,dx)$ denote the Lebesgue space of (classes of) 
functions which are $p$-th power integrable with respect to the weighted 
measure $w(x)\,dx$ on ${\mathcal O}$. Finally, following a well-established 
custom, $A(r)\sim B(r)$ will mean that each quantity is 
bounded by a fixed multiple of the other, uniformly in the parameter $r$.

\subsection{Kernels with singularities along $\partial\mathbb{R}^{n}_+$}

Recall that $L_p(\RR^n_+,\,x_n^{ap}\,dx)$ stands for the weighted Lebesgue 
space of $p$-th power integrable functions in $\RR^n_+$ corresponding to the 
weight $w(x):=x_n^{ap}$, $x=(x',x_n)\in\RR^n_+$. 
\begin{proposition}\label{tp3}
Let $a\in\RR$, $1<p<\infty$, and assume that ${\mathcal Q}$ is a 
non-negative measurable function on 
$\{\zeta=(\zeta',\zeta_n)\in\RR^{n-1}\times\RR:\,\zeta_n>-1\}$, 
which also satisfies
\begin{equation}\label{CC-1}
\int_{\mathbb{R}^n_+}
{\mathcal Q}(\zeta',\zeta_n-1)\,\zeta_n^{-a-1/p}\,d\zeta<\infty.
\end{equation}
\noindent Then the operator
\begin{equation}\label{CC-2}
Qf(x):=x_n^{-n}\int_{\mathbb{R}^n_+}{\mathcal Q}
\Bigl(\frac{y-x}{x_n}\Bigr)f(y)\,dy,\qquad x=(x',x_n)\in\RR^n_+,
\end{equation}
\noindent initially defined on functions $f\in L_p(\mathbb{R}^n_+)$ with 
compact support in $\mathbb{R}^n_+$, can be extended by continuity to an 
operator acting from $L_p(\RR^n_+,\,x_n^{a p}\,dx)$ into itself, with the 
norm satisfying
\begin{equation}\label{CC-3}
\|Q\|_{{\mathfrak L}(L_p(\RR^n_+,\,x_n^{ap}dx))}\leq\int_{\RR^n_+}
{\mathcal Q}(\zeta',\zeta_n-1)\,\zeta_n^{-a-1/p}\,d\zeta.
\end{equation}
\end{proposition}
\noindent{\bf Proof.} Introducing the new variable 
$\zeta:=(x_n^{-1}(y'-x'), x_n^{-1}y_n)\in\RR^n_+$, we may write 
\begin{equation}\label{CC-4}
|Qf(x)|\leq\int_{\RR^n_+}{\mathcal Q}(\zeta',\zeta_n-1) 
|f(x' +x_n\zeta', x_n\zeta_n)|d\zeta,\qquad\forall\,x\in\RR^n_+.
\end{equation}
\noindent Then, by Minkowski's inequality,
\begin{eqnarray}\label{CC-5}
\|Qf\|_{L_p(\RR^n_+,x_n^{a p}\,dx)} & \leq & \int_{\RR^n_+}
{\mathcal Q}(\zeta',\zeta_n-1)\Bigl(\int_{\RR^n_+} 
x_n^{a p}|f(x'+x_n\zeta',x_n\zeta_n)|^p\,dx\Bigr)^{1/p}d\zeta
\nonumber\\[6pt]
& = & \Bigl(\int_{\RR^n_+}{\mathcal Q}(\zeta',\zeta_n-1)\,
\zeta_n^{-a-1/p}\,d\zeta\Bigr)\|f\|_{L_p(\RR^n_+,x_n^{a p}\,dx)},
\end{eqnarray}

\noindent as desired. 
\hfill$\Box$
\vskip 0.08in

Recall that $\bar{y}:=(y',-y_n)$ if $y=(y',y_n)\in\RR^{n-1}\times\RR$. 
\begin{corollary}\label{Cor1}
Consider 
\begin{equation}\label{1a}
Rf(x):=\int_{\mathbb{R}^n_+}\frac{\log\,\bigl(\frac{|x-y|}{x_n}+2\bigr)}
{|x-\bar{y}|^n} f(y)\,dy,\qquad x=(x',x_n)\in\RR^n_+.
\end{equation}
\noindent Then for each $1<p<\infty$ and each $a\in (-1/p,1-1/p)$ the 
operator $R$ is bounded from $L_p(\RR^n_+,\,x_n^{a p}\,dx)$ into itself.
Moreover, there exists $c(n)$, independent of $a$, $p$ and $s:=1-a-1/p$, for which 
\begin{equation}\label{70a}
\|R\|_{{\mathfrak L}(L_p(\RR^n_+,\,x_n^{a p}\,dx))}
\leq\frac{c(n)\,p^2}{(pa+1)(p(1-a)-1)}=\frac{c(n)}{s(1-s)}.
\end{equation}
\end{corollary}
\noindent{\bf Proof.} The result follows from Proposition~\ref{tp3} with
\begin{equation}\label{CC-6}
{\mathcal Q}(\zeta):=\frac{\log\,(|\zeta|+2)}{(|\zeta|^2+1)^{n/2}},
\end{equation}
\noindent and from the obvious inequality $2|x-\bar{y}|^2 \geq |x-y|^2 +x_n^2$.
\hfill$\Box$
\vskip 0.08in

Let us note here that Corollary~\ref{Cor1} immediately yields the following. 
\begin{corollary}\label{Cor2}
Consider 
\begin{equation}\label{2a}
Kf(x):=\int_{\RR^n_+}\frac{f(y)}{|x-\bar{y}|^n}\,dy,\qquad x\in\RR^n_+.
\end{equation}
\noindent Then for each $1<p<\infty$ and $a\in (-1/p,1-1/p)$ the 
operator $K$ is bounded from $L_p(\RR^n_+,\,x_n^{a p}\,dx)$ into itself. 
In addition, there exists $c(n)$, independent of $a$, $p$ and $s:=1-a-1/p$, 
such that 
\begin{equation}\label{71a}
\|K\|_{{\mathfrak L}(L_p(\RR^n_+,\,x_n^{a p}\,dx))}
\leq \frac{c(n)\,p^2}{(pa+1)(p(1-a)-1)}=\frac{c(n)}{s(1-s)}. 
\end{equation}
\end{corollary}
Recall that the barred integral stands for the mean-value (taken in the 
integral sense). 
\begin{lemma}\label{Lem1}
Assume that $1<p<\infty$, $a\in (-1/p,1-1/p)$, and recall the operator 
$K$ introduced in {\rm (\ref{2a})}. Further, consider a non-negative, 
measurable function $w$ defined on $\RR^n_+$ and fix a family of 
balls ${\mathcal F}$ which form a Whitney covering of $\RR^n_+$. 
Then the norm of $wK$ as an operator from $L_p(\RR^n_+,\,x_n^{a p}\,dx)$ 
into itself is equivalent to
\begin{equation}\label{CC-7}
\sup\limits_{B\in{\mathcal F}}\meanint_{\!\!\!\!B}w(y)^p\,dy.
\end{equation}
\noindent Furthermore, there exists $c(n)$, independent 
of $w$, $p$, $a$ and $s:=1-a-1/p$, such that 
\begin{equation}\label{CC-70}
\|w\,K\|_{{\mathfrak L}(L_p(\RR^n_+,\,x_n^{a p}\,dx))}\leq\frac{c(n)}{s(1-s)}
\sup\limits_{B\in{\mathcal F}}\Bigl(\meanint_{\!\!\!\!B}w(y)^p\,dy\Bigr)^{1/p}.
\end{equation}
\end{lemma}
\noindent{\bf Proof.} Fix $f\geq 0$ and denote by $|B|$ the Euclidean 
volume of $B$. Sobolev's embedding theorem allows us to write
\begin{equation}\label{CC-8}
\|Kf\|^p_{L_\infty(B)}\leq c(n)\,|B|^{-1}\sum_{j=0}^n |B|^{jp/n}
\|\nabla_j Kf\|^p_{L_p(B)},\qquad\forall\,B\in{\mathcal F}. 
\end{equation}
\noindent Hence,
\begin{equation}\label{N1-bis}
\int_{\RR^n_+}|x_n^{a}w(x)(Kf)(x)|^p\,dx\leq c(n)\,
\sup\limits_{B\in{\mathcal F}}\meanint_{\!\!\!\!B}w(y)^p\,dy
\int_{\RR^n_+}x_n^{pa}\sum_{0\leq j\leq l}x_n^{jp}|\nabla_j Kf|^p\,dx. 
\end{equation}

\noindent Observing that $x_n^j|\nabla_j\, Kf|\leq c(n)\,Kf$ and referring to 
Corollary~\ref{Cor2}, we arrive at the required upper estimate for the norm 
of $wK$. The lower estimate is obvious.
\hfill$\Box$
\vskip 0.08in
We momentarily pause in order to collect some definitions and set up 
basic notation pertaining to functions with bounded mean oscillations.
Let $f$ be a locally integrable function defined on $\mathbb{R}^n$ and 
define the seminorm
\begin{equation}\label{semi1}
[f]_{{\rm BMO}(\mathbb{R}^n)}:=\sup_{B}
\meanint_{\!\!\!B}\,\Bigl|f(x)-\meanint_{\!\!\!B}f(y)\,dy\Bigr|\,dx, 
\end{equation}
\noindent where the supremum is taken over all balls $B$ in ${\mathbb{R}^n}$. 
If $f$ is a locally integrable function on $\mathbb{R}^n_+$, set 
\begin{equation}\label{semi2}
[f]_{{\rm BMO}(\mathbb{R}^n_+)}:=\mathop{\hbox{sup}}_{(B)}
\meanint_{\!\!\!B\cap\mathbb{R}^n_+}\,
\Bigl|f(x)-\meanint_{\!\!\!B\cap\mathbb{R}^n_+}f(y)\,dy\Bigr|\,dx, 
\end{equation}
\noindent where, this time, the supremum is taken over the collection 
$(B)$ of all balls $B$ with centers in $\overline{\mathbb{R}^n_+}$.
Then the following inequalities are straightforward
\begin{equation}\label{N4-bis}
[f]_{{\rm BMO}(\mathbb{R}^n_+)}\leq\mathop{\hbox{sup}}_{(B)}
\meanint_{\!\!\!B\cap\mathbb{R}^n_+}\,\meanint_{\!\!\!B\cap\mathbb{R}^n_+}
\Bigl|f(x)-f(y)\,\Bigr|\,dxdy\leq 2\,[f]_{{\rm BMO}(\mathbb{R}^n_+)}.
\end{equation}
\noindent We also record here the equivalence relation 
\begin{equation}\label{semi}
[f]_{{\rm BMO}(\mathbb{R}^n_+)}\sim [{\rm Ext}\,f]_{{\rm BMO}(\mathbb{R}^n)},
\end{equation}
\noindent where ${\rm Ext}\,f$ is the extension of $f$ onto 
$\mathbb{R}^n$ as an even function in $x_n$.
Finally, by ${\rm BMO}({\mathbb{R}^n_+})$ we denote the collection of 
equivalence classes, mod constants, of functions $f$ on $\mathbb{R}^n_+$ 
for which $[f]_{{\rm BMO}({\mathbb{R}^n_+})}<\infty$. 
\begin{proposition}\label{tp3-bis}
Let $b\in{\rm BMO}(\RR^n_+)$ and consider the operator 
\begin{equation}\label{eqp8-bis}
Tf(x):=\int_{\RR^n_+}\frac{|b(x)-b(y)|}{|x-\bar{y}|^n}f(y)\,dy,
\qquad x\in\RR^n_+.
\end{equation}
\noindent Then for each $p\in(1,\infty)$ and $a\in (-1/p,1-1/p)$
\begin{equation}\label{eqp9bis}
T:L_p(\RR^n_+,\,x_n^{a p}\,dx)\longrightarrow
L_p(\RR^n_+,\,x_n^{a p}\,dx)
\end{equation}
\noindent is a well-defined, bounded operator, such that if $s:=1-a-1/p$ then 
\begin{equation}\label{CC-71}
\|T\|_{{\mathfrak L}(L_p(\RR^n_+,\,x_n^{a p}\,dx))}
\leq\frac{c(n)}{s(1-s)}\,[b]_{{\rm BMO}(\RR^n_+)}.
\end{equation}
\end{proposition}
\noindent{\bf Proof.} Given $x\in\RR^n_+$ and $r>0$, we shall use 
the abbreviations
\begin{equation}\label{CC-9}
\bar{b}_r(x):=\meanint_{\!\!\!B(x,r)\cap\RR^n_+}b(y)\,dy,
\qquad\quad D_r(x):=|b(x)-\bar{b}_{r}(x)|,
\end{equation}
\noindent and make use of the integral operator
\begin{equation}\label{CC-10}
Sf(x):=\int_{\mathbb{R}^n_+}\frac{D_{|x-\bar{y}|}(x)}{|x-\bar{y}|^n}
\,f(y)\,dy,\qquad x\in\RR^n_+,
\end{equation}
\noindent as well as its adjoint $S^*$. For each nonnegative, 
measurable function $f$ on $\RR^n_+$ and each $x\in\RR^n_+$, 
\begin{eqnarray}\label{CC-11}
Tf(x) & \leq & Sf(x)+S^*f(x)+\int_{\RR^n_+}
\frac{|\bar{b}_{|x-\bar{y}|}(x)-\bar{b}_{|x-\bar{y}|}(y)|}{|x-\bar{y}|^n}\,
f(y)\,dy
\nonumber\\[6pt]
& \leq & Sf(x)+S^*f(x)+c(n)\,[b]_{{\rm BMO}(\RR^n_+)}Kf(x),
\end{eqnarray}
\noindent where $K$ has been introduced in (\ref{2a}). Making use of 
Corollary~\ref{Cor2}, we need to estimate only the norm of $S$. Obviously,
\begin{equation}\label{CC-12}
Sf(x)\leq D_{x_n}(x)Kf(x)+\int_{\RR^n_+}
\frac{|\bar{b}_{x_n}(x)-\bar{b}_{|x-\bar{y}|}(x)|}{|x-\bar{y}|^n}\,f(y)\,dy.
\end{equation}
\noindent Setting $r=|x-\bar{y}|$ and $\rho=x_n$ in the standard inequality
\begin{equation}\label{CC-13}
|\bar{b}_\rho(x)-\bar{b}_r(x)|\leq c(n)\,\log\,\Bigl(\frac{r}{\rho}+1\Bigr)
[b]_{{\rm BMO}(\RR^n_+)},
\end{equation}
\noindent where $r>\rho$, we arrive at
\begin{equation}\label{CC-14}
Sf(x)\leq D_{x_n}(x)Kf(x)+ c(n)\,[b]_{BMO(\mathbb{R}^n_+)}\, Rf(x),
\end{equation}
\noindent where $R$ is defined in (\ref{1a}). Let ${\mathcal F}$ 
be a Whitney covering of $\RR^n_+$ with open balls. For an arbitrary
$B\in{\mathcal F}$, denote by $\delta$ the radius of $B$. 
By Lemma~\ref{Lem1} with $w(x):=D_{x_n}(x)$, the norm of the operator 
$D_{x_n}(x)K$ does not exceed
\begin{eqnarray}\label{CC-15}
\sup\limits_{B\in{\mathcal F}}\Bigl(\meanint_{\!\!\!B}|D_{x_n}(x)|^p\,dx
\Bigr)^{1/p} 
& \leq & c(n)\,\sup\limits_{B\in{\mathcal F}}\Bigl(\meanint_{\!\!\!B}
|b(x)-\bar{b}_\delta(x)|^p\,dx\Bigr)^{1/p}
+c(n)\,[b]_{{\rm BMO}(\mathbb{R}^n_+)}
\nonumber\\[6pt]
&\leq & c(n)\,[b]_{{\rm BMO}(\mathbb{R}^n_+)},
\end{eqnarray}
\noindent by the John-Nirenberg inequality. Here we have also used the 
triangle inequality and the estimate (\ref{CC-13}) in order to replace 
$\bar{b}_{x_n}(x)$ in the definition of $D_{x_n}(x)$ by $\bar{b}_{\delta}(x)$.
The intervening logarithmic factor is bounded independently of $x$ since $x_n$
is comparable with $\delta$, uniformly for $x\in B$. With this estimate
in hand, a reference to Corollary~\ref{Cor1} gives that 
\begin{equation}\label{CC-16}
\begin{array}{l}
S:L_p(\RR^n_+,\,x_n^{a p}\,dx)\to L_p(\RR^n_+,\,x_n^{a p}\,dx)
\mbox{ boundedly} 
\\[6pt]
\mbox{for each }p\in(1,\infty)\mbox{ and each }a\in (-1/p,1-1/p).
\nonumber
\end{array}
\end{equation}
\noindent The corresponding estimate for the norm of $S$ is implicit in the
above argument. 
By duality, it follows that $S^*$ enjoys the same property and,
hence, the operator $T$ is bounded on $L_p(\RR^n_+,\,x_n^{a p}\,dx)$
for each $p\in(1,\infty)$ and $a\in(-1/p,1-1/p)$, thanks to (\ref{CC-11})
and Corollary~\ref{Cor2}. The fact that the operator norm of $T$ can be 
estimated in the desired fashion is implicit in the above reasoning. 
\hfill$\Box$

\subsection{Preliminary estimates for singular integrals on weighted 
Lebesgue spaces}

We need the analogue of Proposition~\ref{tp3-bis} for the class 
of Mikhlin-Calder\'on-Zygmund singular integral operators. 
Let $S^{n-1}$ stand for the unit sphere in $\RR^n$ and recall that
\begin{equation}\label{CZ-op}
{\mathcal S}f(x)=p.v.\int_{\RR^n}k(x,x-y)f(y)\,dy,\qquad x\in\RR^n,
\end{equation}
\noindent (where $p.v.$ indicates that the integral is taken in the
principal value sense) is called a Mikhlin-Calder\'on-Zygmund operator 
provided the function $k:\RR^n\times(\RR^n\setminus\{0\})\to\RR$ satisfies:
\begin{itemize}
\item[(i)] $k(x,\cdot)\in C^\infty(\RR^n\setminus\{0\})$
and, for almost each $x\in\RR^n$,
\begin{equation}\label{kk-est}
\max_{|\alpha|\leq 2n}\|D_z^\alpha k(x,z)\|_{L_\infty(\RR^n\times S^{n-1})}
<\infty.
\end{equation}
\item[(ii)]  $k(x,\lambda z)=\lambda ^{-n}k(x,z)$ for each $z\in\RR^n$
and each $\lambda\in\RR$, $\lambda>0$;
\item[(iii)] $\int_{S^{n-1}} k(x,\omega)\,d\omega=0$,
where $d\omega$ indicates integration with respect to $\omega\in S^{n-1}$.
\end{itemize}
\vskip 0.08in
It is well-known that the Mikhlin-Calder\'on-Zygmund operator ${\mathcal S}$ 
and its commutator $[{\mathcal S},b]$ with the operator of multiplication by 
a function $b\in{\rm BMO}(\mathbb{R}^n_+)$ are bounded operators in 
$L_p(\mathbb{R}^n_+)$ for each $1<p<\infty$. Then 
\begin{equation}\label{b1}
\|{\mathcal S}\|_{{\mathfrak L}(L_p(\mathbb{R}^n_+))}\leq c(n)\,p\,p',\qquad 
\|\,[{\mathcal S},b]\,\|_{{\mathfrak L}(L_p(\mathbb{R}^n_+))}
\leq c(n)\,p\,p'\,[b]_{{\rm BMO}(\RR^n_+)},
\end{equation}
\noindent where $1/p+1/p'=1$ and $c(n)$ depends only on $n$ and the quantity 
in (\ref{kk-est}). The first estimate in (\ref{b1}) goes back to the work 
of A.\,Calder\'on and A.\,Zygmund (see also the comment on 
p.\,22 of {\bf\cite{St}} regarding the dependence on the parameter 
$p$ of the constants involved).
The second estimate in (\ref{b1}) was originally proved for 
convolution type operators by R.\,Coifman, R.\,Rochberg and G.\,Weiss in 
{\bf\cite{CRW}} and a standard expansion in spherical harmonics allows
to extend this result to the case of operators with variable-kernels 
of the type considered above. 

We wish to extend (\ref{b1}) to the case when the Lebesgue measure
is replaced by $x_n^{ap}\,dx$, with $1<p<\infty$ and $a\in(-1/p,1-1/p)$. 
Incidentally, $a\in(-1/p,1-1/p)$ corresponds precisely to the range 
of $a$'s for which $w(x):=x_n^{ap}$ is a weight in Muckenhoupt's $A_p$ 
class, although here we prefer to give a direct, elementary proof. 
\begin{proposition}\label{Prop4}
Retain the above conventions and hypotheses. Then the operator ${\mathcal S}$ 
and its commutator $[{\mathcal S},b]$ with a function 
$b\in{\rm BMO}(\mathbb{R}^n_+)$ are bounded when acting from 
$L_p(\RR^n_+,\,x_n^{ap}\,dx)$ into itself for each $p\in(1,\infty)$ 
and $a\in (-1/p, 1-1/p)$. Then, with $s:=1-a-1/p$ and $1/p+1/p'=1$,   
\begin{eqnarray}\label{b2}
&&\|{\mathcal S}\|_{{\mathfrak L}(L_p(\mathbb{R}^n_+,\,x_n^{ap}\,dx))}
\leq c(n)\Bigl(p\,p'+\frac{1}{s(1-s)}\Bigr),
\\[6pt]
&&\|\,[{\mathcal S},b]\,\|_{{\mathfrak L}(L_p(\mathbb{R}^n_+,\,x_n^{ap}\,dx))} 
\leq c(n)\Bigl(p\,p'+\frac{1}{s(1-s)}\Bigr)\,[b]_{{\rm BMO}(\RR^n_+)}.
\label{b2'}
\end{eqnarray}
\end{proposition}
\noindent{\bf Proof.} Let $\chi_j$ be the characteristic function of the 
layer $2^{j/2}<x_n\leq 2^{1+j/2}$, $j=0,\pm 1,\ldots$, so that
$\sum_{j\in\ZZ}\chi_j=2$. We then write ${\mathcal S}$ as the sum 
${\mathcal S}_1+{\mathcal S}_2$, where
\begin{equation}\label{CC-17}
{\mathcal S}_1:=\frac{1}{4}\sum_{|j-k|\leq 3}\chi _j{\mathcal S}\chi_k.
\end{equation}
\noindent The following chain of inequalities is evident
\begin{eqnarray}\label{CC-18}
\|{\mathcal S}_1\, f\|_{L_p(\mathbb{R}^n_+,\,x_n^{ap}\,dx)}
& \leq & \Bigl(\sum_j\int_{\RR^n_+}\chi_j(x)\,
\Bigl|{\mathcal S}\Bigl(\sum_{|k-j|\leq 3}\chi_k f\Bigr)(x)\Bigr|^p\,
x_n^{ap}\,dx\Bigr)^{1/p}
\nonumber\\[6pt]
& \leq & c(n)\Bigl(\sum_j\int_{\mathbb{R}^n_+}
\Bigl|{\mathcal S}\Bigl(\sum_{|k-j|\leq 3} 
\chi_k 2^{ja/2} f\Bigr)(x)\Bigr|^p\,dx\Bigr)^{1/p}.
\end{eqnarray}
\noindent In concert with the first estimate in (\ref{b1}), this entails
\begin{eqnarray}\label{CC-19}
\|{\mathcal S}_1\,f\|_{L_p(\mathbb{R}^n_+,\,x_n^{ap}\,dx)}
& \leq & c(n)\,p\, p'\Bigl(\sum_j 
\int_{\RR^n_+}\Bigl(\sum_{|k-j|\leq 3}
\chi_k 2^{ja/2}|f|\Bigr)^p\,dx\Bigr)^{1/p}
\nonumber\\[6pt]
& \leq & c(n)\,p\,p'\Bigl(\int_{\RR^n_+}|f(x)|^p\,x_n^{ap}\,dx\Bigr)^{1/p},
\end{eqnarray}
\noindent which is further equivalent to
\begin{equation}\label{b3}
\|{\mathcal S}_1\|_{{\mathfrak L}(L_p(\mathbb{R}^n_+,\,x_n^{ap}\,dx))}
\leq c(n)\,p\,p'.
\end{equation}
\noindent Applying the same  argument to $[{\mathcal S}_1,b]$ and referring 
to (\ref{b1}), we arrive at
\begin{equation}\label{b4}
\|\,[{\mathcal S}_1,b]\,\|_{{\mathfrak L}(L_p(\mathbb{R}^n_+,\,x_n^{ap}\,dx))} 
\leq c(n)\,p\,p'\,[b]_{{\rm BMO}(\RR^n_+)}.
\end{equation}

It remains to obtain the analogues of (\ref{b3}) and ({\ref{b4}) with 
${\mathcal S}_2$ in place of ${\mathcal S}_1$. One can check directly that the 
modulus of the kernel of ${\mathcal S}_2$ does not exceed 
$c(n)\,|x-\bar{y}|^{-n}$ and that the modulus of the kernel of 
$[{\mathcal S}_2,b]$ is majorized by $c(n)\,|b(x)-b(y)|\,|x-\bar{y}|^{-n}$. 
Then the desired conclusions follow from 
Corollary~\ref{Cor2} and Proposition~\ref{tp3-bis}. 
\hfill$\Box$

\subsection{BMO type estimates for Gagliardo's extension operator}  

Here we shall revisit a certain operator $T$, extending functions 
defined on $\mathbb{R}^{n-1}$ into functions defined on $\mathbb{R}^n_+$,
first introduced by E.\,Gagliardo in {\bf\cite{Ga}}. Fix a smooth, radial, 
decreasing, even, non-negative function $\zeta$ in $\mathbb{R}^{n-1}$ 
such that $\zeta(t)=0$ for $|t|\geq 1$ and 
\begin{equation}\label{zeta-int}
\int\limits_{\mathbb{R}^{n-1}}\zeta(t)\,dt=1.
\end{equation}
\noindent (A standard choice is 
$\zeta(t):=c\,{\rm exp}\,(-1/(1-|t|^2)_+)$ for a suitable $c$.) 
Following {\bf\cite{Ga}} we then define 
\begin{equation}\label{10.1.20}
(T\varphi)(x',x_n):=\int\limits_{\mathbb{R}^{n-1}}\zeta(t)\varphi(x'+x_nt)\,dt,
\qquad(x',x_n)\in\RR^n_+,
\end{equation}
\noindent acting on functions $\varphi$ from $L_1(\mathbb{R}^{n-1},loc)$. 
To get started, we note that 
\begin{eqnarray}\label{10.1.29}
\nabla_{x'}(T\varphi)(x',x_n)
& = & \int\limits_{\mathbb{R}^{n-1}}\zeta(t)\nabla\varphi(x'+tx_n)\,dt,
\\[6pt]
{\partial\over\partial x_n}(T\varphi)(x',x_n)
& = & \int\limits_{\mathbb{R}^{n-1}}\zeta(t)\,t\,\nabla\varphi(x'+t x_n)\,dt,
\label{10.1.30}
\end{eqnarray}
\noindent and, hence, we have the estimate
\begin{equation}\label{10.1.21}
\|\nabla_x\,(T\varphi)\|_{L_\infty(\mathbb{R}^n_+)} 
\leq c\,\|\nabla_{x'}\,\varphi\|_{L_\infty(\mathbb{R}^{n-1})}.
\end{equation}
\noindent Refinements of (\ref{10.1.21}) are contained in the 
Lemmas~\ref{lem7}-\ref{lem8} below. 
\begin{lemma}\label{lem7}
{\rm (i)} For each multi-indices $\alpha$ with $|\alpha|>1$
there exists $c>0$ such that 
\begin{equation}\label{1.2}
\Bigl|D^\alpha_{x}(T\varphi)(x)\Bigr|
\leq c\,x_n ^{1-|\alpha|}[\nabla\varphi]_{\rm BMO(\mathbb{R}^{n-1})},
\qquad\forall\,x=(x',x_n)\in\mathbb{R}^n_+. 
\end{equation}
{\rm (ii)} There exists $c>0$ such that 
\begin{equation}\label{Tfi}
\Bigl|(T\varphi)(x)-\varphi(x')\Bigr|
\leq c\,x_n[\nabla\varphi]_{\rm BMO(\mathbb{R}^{n-1})},\qquad
\forall\,x=(x',x_n)\in\mathbb{R}^n_+.
\end{equation}
\end{lemma}
\noindent{\bf Proof.} Rewriting (\ref{10.1.30}) as 
\begin{equation}\label{1.3}
{\partial\over\partial x_n}(T\varphi)(x',x_n)
=x_n^{1-n}\int\limits_{\mathbb{R}^{n-1}}\zeta\Bigl(\frac{\xi-x'}{x_n}\Bigr)
\frac{\xi-x'}{x_n}\Bigl(\nabla\varphi(\xi)-\meanint_{\!\!\!|z-x'|<x_n} 
\nabla\varphi(z)dz\Bigr)d\xi
\end{equation}
\noindent we obtain
\begin{equation}\label{1.12}
\Bigl|D^\gamma_{x}{\partial\over\partial x_n}(T\varphi)(x)\Bigr|
\leq c\,x_n^{-|\gamma|}[\nabla\varphi]_{\rm BMO(\mathbb{R}^{n-1})}
\end{equation}
\noindent for every non-zero multi-index $\gamma$. Furthermore, for 
$i=1,\ldots n-1$, by (\ref{10.1.29})
\begin{equation}\label{TT1}
\frac{\partial}{\partial x_i}\nabla_{x'}(T\varphi)(x) 
=x_n^{1-n}\int\limits_{\mathbb{R}^{n-1}}\partial_i\zeta
\Bigl(\frac{\xi-x'}{x_n}\Bigr)\Bigl(\nabla\varphi(\xi)
-\meanint_{\!\!\!|z-x'|<x_n}\nabla\varphi(z)dz\Bigr)d\xi.
\end{equation}
\noindent Hence, once again
\begin{equation}\label{TT2}
\Bigl|D^\gamma_x\frac{\partial}{\partial x_i}\nabla_{x'}(T\varphi)(x)\Bigr| 
\leq c\,x_n^{-|\gamma|-1}[\nabla\varphi]_{\rm BMO(\mathbb{R}^{n-1})},
\end{equation}
\noindent and the estimate claimed in (i) follows.
Finally, (ii) is a simple consequence of (i) and the fact that 
$(T\varphi)|_{\RR^{n-1}}=\varphi$. 
\hfill$\Box$
\vskip 0.08in

\noindent{\bf Remark.} In concert with Theorem~2 on p.\,62-63 in 
{\bf\cite{St}}, formula (\ref{10.1.29}) yields the pointwise estimate
\begin{equation}\label{HL-Max}
|\nabla\,(T\varphi)(x)|
\leq c\,{\mathcal M}(\nabla\varphi)(x'),\qquad x=(x',x_n)\in\RR^{n}_+,
\end{equation}
\noindent where ${\mathcal M}$ is the classical Hardy-Littlewood 
maximal function (cf., e.g., Chapter~I in {\bf\cite{St}}). 
As for higher order derivatives, an inspection of the above proof 
reveals that 
\begin{equation}\label{sharp}
\Bigl|D_x^\alpha(T\varphi)(x)\Bigr|
\leq c\,x_n^{1-|\alpha|}(\nabla\varphi)^\#(x'),\qquad (x',x_n)\in\RR^{n},
\end{equation}
\noindent holds for each multi-index $\alpha$ with $|\alpha|>1$, 
where $(\cdot)^\#$ is the Fefferman-Stein sharp maximal function 
(cf. {\bf\cite{FS}}). 
\begin{lemma}\label{lem8}
If $\nabla_{x'}\varphi\in{\rm BMO}(\mathbb{R}^{n-1})$ then 
$\nabla(T\varphi)\in{\rm BMO}(\mathbb{R}^n_+)$ and 
\begin{equation}\label{1.8}
[\nabla(T\varphi)]_{{\rm BMO}(\mathbb{R}^n_+)}
\leq c\,[\nabla_{x'}\varphi]_{{\rm BMO}(\mathbb{R}^{n-1})}.
\end{equation}
\end{lemma}
\noindent{\bf Proof.} Since $(T\varphi)(x',x_n)$ is even with respect to $x_n$,
it suffices to estimate $[\nabla_x(T\varphi)]_{{\rm BMO}(\mathbb{R}^n)}$. 
Let $Q_r$ denote a cube with side-length $r$ centered at the point 
$\eta =(\eta',\eta_n)\in\RR^{n-1}\times\RR$. Also let $Q'_r$ be the 
projection of $Q_r$ on $\mathbb{R}^{n-1}$. Clearly,
\begin{equation}\label{xxx}
\nabla_{x'}(T\varphi)(x',x_n)-\nabla _{x'}\varphi(x')
=x_n^{1-n}\int\limits_{\mathbb{R}^{n-1}}\zeta\Bigl(\frac{\xi -x'}{x_n}\Bigr)
(\nabla\varphi(\xi)-\nabla\varphi(x'))\,d\xi.
\end{equation}
Suppose that $|\eta_n|<2r$ and write 
\begin{eqnarray}\label{1.10}
\int_{Q_r}\Bigl|\nabla_{x'}(T\varphi)(x',x_n)-\nabla_{x'}\varphi(x')\Bigr|\,dx 
& \leq & c\,r^{2-n}\int_{Q'_{4r}}
\int_{Q'_{4r}}|\nabla\varphi(\xi)-\nabla\varphi(z)|\,dz\,d\xi.
\nonumber\\[6pt]
& \leq & c\,r^n[\nabla\varphi]_{{\rm BMO}(\mathbb{R}^{n-1})}.
\end{eqnarray}
Therefore, for $|\eta_n|<2r$
\begin{eqnarray}\label{1.18}
&& \meanint_{\!\!\!Q_r}\meanint_{\!\!\!Q_r}|\nabla_{x'}T\varphi(x)-\nabla_{y'}
T\varphi(y)|\,dxdy
\\[6pt]
&&\qquad\qquad
\leq 2\meanint_{\!\!\!Q_r}|\nabla_{x'}T\varphi(x)-\nabla\varphi(x')|\,dx
+\meanint_{\!\!\!Q'_r}\meanint_{\!\!\!Q'_r}
|\nabla\varphi(x')-\nabla\varphi(y')|\,dx'dy' 
\leq c\,[\nabla\varphi]_{{\rm BMO}(\mathbb{R}^{n-1})}.
\nonumber
\end{eqnarray}
Next, consider the case when $|\eta_n|\geq 2r$ and let $x$ and $y$ 
be arbitrary points in $Q_r(\eta)$. Then, using the generic abbreviation 
$\bar{f}_E:={\displaystyle{\meanint_{\!\!\!E}f}}$, we may write 
\begin{eqnarray}\label{arTT}
|\nabla_{x'}T\varphi(x)-\nabla_{y'}T\varphi(y)|
& \leq & 
\int\limits_{\mathbb{R}^{n-1}}\Bigl|x_n^{1-n}\zeta\Bigl(\frac{\xi-x'}{x_n}
\Bigr)-y_n^{1-n}\zeta\Bigl(\frac{\xi-y'}{y_n}\Bigr)
\Bigr|\Bigl|\nabla\varphi(\xi)- {\overline{\nabla\varphi}}_{Q'_{2|\eta_n|}}
\Bigr|\,d\xi 
\nonumber\\[6pt]
& \leq & \frac{c\,r}{|\eta_n|^n}\int_{Q'_{2|\eta_n|}}
\Bigl|\nabla\varphi(\xi)
-{\overline{\nabla\varphi}}_{Q'_{2|\eta_n|}}\Bigr|\,d\xi 
\nonumber\\[6pt]
& \leq & c\,[\nabla\varphi]_{{\rm BMO}(\mathbb{R}^{n-1})}. 
\end{eqnarray}
\noindent Consequently, for $|\eta_n|\geq 2r$, 
\begin{equation}\label{QrT3}
\meanint_{\!\!\!Q_r}\meanint_{\!\!\!Q_r}
|\nabla_{x'}T\varphi(x)-\nabla_{y'}T\varphi(y)|\,dxdy
\leq c\,[\nabla\varphi]_{{\rm BMO}(\mathbb{R}^{n-1})}
\end{equation}
\noindent which, together with (\ref{1.18}), gives
\begin{equation}\label{QrT4}
[\nabla_{x'}T\varphi]_{{\rm BMO}(\mathbb{R}^{n})}
\leq c\,[\nabla\varphi]_{\rm BMO(\mathbb{R}^{n-1})}.
\end{equation}
\noindent This inequality and (\ref{1.12}) with $|\gamma|=0$ imply (\ref{1.8}).
\hfill$\Box$

\section{The Dirichlet problem in $\mathbb{R}^n_+$ for variable 
coefficient systems}
\setcounter{equation}{0}

\subsection{The setup}
For 
\begin{equation}\label{indices}
1<p<\infty,\quad -\frac{1}{p}<a<1-\frac{1}{p}\quad 
\mbox{and}\quad m\in\NN, 
\end{equation}
\noindent we let $V^{m,a}_p(\mathbb{R}^n_+)$ denote the weighted Sobolev space 
associated with the norm
\begin{equation}\label{defVVV}
\|u\|_{V_p^{m,a}(\mathbb{R}^n_+)}:=\Bigl(\sum_{|\beta|\leq m} 
\int_{\mathbb{R}^n_+}|x_n^{|\beta|-m} D^\beta u(x)|^p\,x_n^{pa}\,dx
\Bigr)^{1/p}.
\end{equation}
\noindent It is easily proved that $C_0^\infty(\mathbb{R}^n_+)$ is dense in 
$V^{m,a}_p(\mathbb{R}^n_+)$. Moreover, by the one-dimensional Hardy's 
inequality (see, for instance, {\bf\cite{Maz1}}, formula (1.3/1)), we have 
(with $s:=1-a-1/p$, as usual), 
\begin{equation}\label{4.60}
\|u\|_{V_p^{m,a}(\mathbb{R}^n_+)}\leq C\,s^{-1}
\Bigl(\sum_{|\beta|=m}\int_{\mathbb{R}^n_+}
|D^\beta u(x)|^p\,x_n^{pa}\,dx\Bigr)^{1/p}
\,\,\,\,\mbox{ for }\,\,u\in C_0^\infty(\mathbb{R}^n_+). 
\end{equation}
\noindent The dual of $V^{m,-a}_{p'}(\mathbb{R}^n_+)$ will be denoted by 
$V^{m,a}_p(\mathbb{R}^n_+)$, where $1/p+1/p'=1$. 

Consider now the operator
\begin{equation}\label{LxDu}
{L}(x,D_x)u:=\sum_{|\alpha|,|\beta|\leq m}
D^\alpha_x({A}_{\alpha\beta}(x)\,x_n^{|\alpha|+|\beta|-2m} D_x^\beta\,u)
\end{equation}
\noindent where ${A}_{\alpha\beta}$ are $\CC^{l\times l}$-valued functions 
in $L_\infty(\mathbb{R}^n_+)$. We shall use the notation 
$\ring{L}(x,D_x)$ for the principal part of ${L}(x,D_x)$, i.e.,
\begin{equation}\label{Lcirc}
\ring{L}(x,D_x)u:=\sum_{|\alpha|=|\beta|=m} 
D^\alpha_x({A}_{\alpha\beta}(x)\,D_x^\beta\,u).
\end{equation}

\subsection{Solvability and regularity result}

\begin{lemma}\label{lem5}
Assume that there exists $\kappa=const>0$ such that the coercivity condition
\begin{equation}\label{B5}
\Re\int_{\mathbb{R}^n_+}\sum_{|\alpha|=|\beta|=m}
\langle A_{\alpha\beta}(x)\,D^\beta u(x),\,D^\alpha u(x)
\rangle_{\mathbb{C}^l}\,dx 
\geq \kappa\sum_{|\gamma|=m}\|D^\gamma\,u\|^2_{L_2(\mathbb{R}^n_+)},
\end{equation}
\noindent holds for all $u\in C^\infty_0(\mathbb{R}^n_+)$, and that 
\begin{equation}\label{B6}
\sum_{|\alpha|=|\beta|= m}\|A_{\alpha\beta}\|_{L_\infty(\mathbb{R}^n_+)} 
\leq \kappa^{-1}.
\end{equation}

{\rm (i)} Let $p\in(1,\infty)$, $a\in(-1/p,1-1/p)$, and suppose that 
\begin{equation}\label{E7}
\frac{1}{s(1-s)}\sum_{{|\alpha|+|\beta|<2m}\atop{0\leq |\alpha|,|\beta|\leq m}}
\|A_{\alpha\beta}\|_{{L_\infty}(\mathbb{R}^n_+)} 
+\sum_{|\alpha|=|\beta|=m}[A_{\alpha\beta}]_{{\rm BMO}(\mathbb{R}^n_+)}
\leq \delta,
\end{equation}
\noindent where $s:=1-a-1/p$, $1/p+1/p'=1$, and $\delta$ satisfies
\begin{equation}\label{E8}
\Bigl(p\,p'+\frac{1}{s(1-s)}\Bigr)\,{\delta}<c(n,m,\kappa)
\end{equation}
\noindent with a sufficiently small constant $c(n,m,\kappa)>0$. Then 
\begin{equation}\label{Liso}
L=L(x,D_x):V_p^{m,a}(\mathbb{R}^n_+)\longrightarrow V_p^{-m,a}(\mathbb{R}^n_+)
\quad\mbox{ isomorphically}. 
\end{equation}

{\rm (ii)} Let $p_i\in (1,\infty)$ and $-1/p_i<a_i<1-1/p_i$, where $i=1,2$. 
Suppose that (\ref{E8}) holds with $p_i$ and $s_i:=1-a_i-1/p_i$ in place of 
$p$ and $s$. Then, if the function $u\in V_{p_1}^{m,a_1}(\mathbb{R}^n_+)$ is 
such that ${L}u\in V_{p_1}^{-m,a_1}(\mathbb{R}^n_+)\cap 
V_{p_2}^{-m,a_2}(\mathbb{R}^n_+)$, it follows that 
$u\in V_{p_2}^{m,a_2}(\mathbb{R}^n_+)$.
\end{lemma}
\noindent{\bf Proof.} The fact that the operator in (\ref{Liso}) is continuous 
is obvious. Also, the existence of a bounded inverse ${L}^{-1}$ for $p=2$ 
and $a=0$ follows from (\ref{B5}) and (\ref{E7})-(\ref{E8}) with $p=2$, $a=0$,
which allow us to implement the Lax-Milgram lemma. 
We shall use the notation $\ring{L}_y$ for the operator $\ring{L}(y,D_x)$, 
corresponding to (\ref{Lcirc}) in which the coefficients have been frozen at 
$y\in\mathbb{R}^n_+$, and the notation $G_y$ for the solution operator
for the Dirichlet problem for $\ring{L}_y$ in $\RR^n_+$ with homogeneous 
boundary conditions. Next, given $u\in V_p^{m,a}(\mathbb{R}^n_+)$, set 
$f:=Lu\in V_p^{-m,a}(\mathbb{R}^n_+)$ so that 
\begin{equation}\label{E10}
\left\{
\begin{array}{l}
{L}(x,D)u=f\qquad\qquad\mbox{in}\,\,\,\mathbb{R}^n_+,\\[6pt]
\displaystyle{\frac{\partial^j\,u}{\partial x_n^j}(x',0)=0}
\qquad\qquad {\rm on}\,\,\mathbb{R}^{n-1},\,\,0\leq j\leq m-1. 
\end{array}
\right.
\end{equation}
\noindent We may now write 
\begin{equation}\label{E12}
u(x)=(G_y f)(x)-(G_{y}(\ring{L}-\ring{L}_y)u)(x)-(G_{y}({L}-\ring{L})u)(x),
\qquad x\in\RR^n_+, 
\end{equation}
\noindent and aim to use (\ref{E12}) in order to express $u$ in terms 
of $f$ (cf. (\ref{IntRepFor})-(\ref{SSS}) below) via integral operators
whose norms we can control. 
First, we claim that whenever $|\gamma|=m$, the norm of the operator 
\begin{equation}\label{ItalianTrick}
V_p^{m,a}(\mathbb{R}^n_+)\ni u
\mapsto D_x^\gamma(G_y(\ring{L}-\ring{L}_y)u)(x)\Bigl|_{x=y}\,
\in L_p(\mathbb{R}^n_+,y_n^{ap}\,dy)
\end{equation}
\noindent does not exceed 
\begin{equation}\label{smallCT}
C\,\Bigl(p\,p'+\frac{1}{s(1-s)}\Bigr)\sum_{|\alpha|=|\beta|=m}
[A_{\alpha\beta}]_{{\rm BMO}(\mathbb{R}^n_+)}.
\end{equation}
\noindent Given the hypotheses under which we operate, 
the expression (\ref{smallCT}) is therefore small if $\delta$ is small. 

In what follows, we denote by $G_y(x,z)$ the integral kernel of $G_y$ and 
integrate by parts in order to move derivatives of the form $D_z^\alpha$ with 
$|\alpha|=m$ from $(\ring{L}-\ring{L}_y)u$ onto $G_y(x,z)$ (the absence 
of boundary terms is due to the fact that $G_y(x,\cdot)$ satisfies homogeneous 
Dirichlet boundary conditions). That (\ref{smallCT}) bounds the norm of 
(\ref{ItalianTrick}) can now be seen by combining Theorem~\ref{th1} with 
(\ref{CC-71}) and Proposition~\ref{Prop4}. 
Let $\gamma$ and $\alpha$ be multi-indices with $|\gamma|=m$, 
$|\alpha|\leq m$ and consider the assignment
\begin{equation}\label{oppsi}
\begin{array}{l}
\displaystyle
C_0^\infty(\RR^n_+)\ni\Psi\mapsto\Bigl(D^\gamma_x\int_{\RR^n_+}
G_y(x,z) D^\alpha_z\frac{\Psi(z)}{z_n^{m-|\alpha|}}\,dz\Bigr)\Bigl|_{x=y}.
\end{array}
\end{equation}
\noindent After integrating by parts, the action of this operator can 
be rewritten in the form
\begin{equation}\label{DgammaInt} 
\Bigl(D^\gamma_x\int_{\RR^n_+}\Bigl[  
\Bigl(\frac{-1}{i}\frac{\partial}{\partial z_n}\Bigr)^{m-|\alpha|}
(-D_z)^\alpha G_y(x,z)\Bigr]\Gamma_\alpha(z)\,dz\Bigr)\Bigl|_{x=y},
\end{equation}
\noindent where 
\begin{equation}\label{Gamma-def}
\Gamma_\alpha(z):=
\left\{
\begin{array}{l}
\Psi(z),\qquad\mbox{if }\,\,\,|\alpha|=m,
\\[10pt]
{\displaystyle{\frac{(-1)^{m-|\alpha|}}{(m-|\alpha|-1)!}
\int_{z_n}^\infty (t-z_n)^{m-|\alpha|-1}\frac{\Psi(z',t)}{t^{m-|\alpha|}}\,dt,
\qquad\mbox{if }\,\,\,|\alpha|<m}}. 
\end{array}
\right.
\end{equation}
\noindent Using  Theorem~\ref{th1} along with (\ref{CC-71}) and 
Proposition~\ref{Prop4}, we may therefore conclude that 
\begin{eqnarray}\label{DZGamma}
&& \Bigl\|\Bigl(D^\gamma_x\int_{\RR^n_+} 
\Bigl[\Bigl(\frac{-1}{i}\frac{\partial}{\partial z_n}\Bigr)^{m-|\alpha|} 
(-D_z)^\alpha G_y(x,z)\Bigr]\Gamma_\alpha(z)\,dz\Bigr)\Bigl|_{x=y}
\Bigr\|_{L_p(\mathbb{R}^n_+,\,y_n^{ap}\,dy)}
\nonumber\\[6pt]
&&\qquad\qquad\qquad\qquad\qquad\qquad\qquad\quad
\leq C\,\Bigl(p\,p'+\frac{1}{s(1-s)}\Bigr)
\|\Gamma_\alpha\|_{L_p(\mathbb{R}^n_+,\,x_n^{ap}\,dx)}.
\end{eqnarray}
\noindent On the other hand, Hardy's inequality gives 
\begin{equation}\label{Cs-1}
\|\Gamma_\alpha\|_{L_p(\mathbb{R}^n_+,\,x_n^{ap}\,dx)} 
\leq\frac{C}{1-s}\,\|\Psi\|_{L_p(\mathbb{R}^n_+,\,x_n^{ap}\,dx)}
\end{equation}
\noindent and, hence, the operator (\ref{oppsi}) can be extended to 
a linear mapping from $C_0^\infty(\mathbb{R}^n_+)$ to 
$L_p(\mathbb{R}^n_+,\,x_n^{ap}\,dx)$ with norm 
$\leq \frac{C}{1-s}\,\Bigl(p\,p'+\frac{1}{s(1-s)}\Bigr)$. 
Next, given $u\in V_p^{m,a}(\mathbb{R}^n_+)$, we let 
$\Psi=\Psi_{\alpha\beta}$ in (\ref{oppsi}) with 
\begin{equation}\label{Pizz}
\Psi_{\alpha\beta}(z):=z_n^{|\beta|-m}\,A_{\alpha\beta}\,D^\beta u(z),
\qquad |\alpha|+|\beta|<2m,
\end{equation}
\noindent and conclude that the norm of the operator 
\begin{equation}\label{altOp}
V_p^{m,a}(\mathbb{R}^n_+)\ni u\mapsto 
D_x^\gamma\,(G_y\,({L}-\ring{L})u)(x)\Bigl|_{x=y}
\in L_p(\mathbb{R}^n_+,y_n^{ap}\,dy)
\end{equation}
\noindent does not exceed
\begin{equation}\label{notExceed}
\frac{C}{1-s}\,\Bigl(p\,p'+\frac{1}{s(1-s)}\Bigr)\sum_{{|\alpha|+|\beta|<2m}
\atop{|\alpha|,|\beta|\leq m}}\|A_{\alpha\beta}\|_{{L_\infty}(\mathbb{R}^n_+)}.
\end{equation}

It is well-known (cf. (1.1.10/6) on p.\,22 of {\bf\cite{Maz1}}) 
that any $u\in V^{m,a}_p(\RR^n_+)$ can be represented as  
\begin{equation}\label{IntUU}
u=K\{D^\sigma u\}_{|\sigma|=m}
\end{equation}
\noindent where $K$ is a linear operator with the property that 
\begin{equation}\label{K-map}
D^\alpha K:L_p(\mathbb{R}^n_+,x_n^{ap}\,dx)\longrightarrow 
L_p(\mathbb{R}^n_+,x_n^{ap}\,dx)
\end{equation}
\noindent is bounded for every multi-index $\alpha$ with $|\alpha|=m$. 
In particular, by (\ref{4.60}), 
\begin{equation}\label{KDu}
\|K\{D^\sigma u\}_{|\sigma|=m}\|_{V_p^{m,a}(\mathbb{R}^n_+)}\leq C\,s^{-1}
\|\{D^\sigma u\}_{|\sigma|=m}\|_{L_p(\mathbb{R}^n_+,x_n^{ap}\,dx)}.
\end{equation}
At this stage, we transform the identity (\ref{E12}) as follows.
First, we express the two $u$'s occurring inside the Green 
operator $G_y$ in the left-hand side of (\ref{E12}) as in (\ref{IntUU}).  
Second, for each $\gamma\in{\mathbb{N}}_0$ with $|\gamma|=m$, 
we apply $D^\gamma$ to both sides of (\ref{E12}) and, finally, set $x=y$. 
The resulting identity reads 
\begin{equation}\label{IntRepFor}
\{D^\gamma u\}_{|\gamma|=m}+S\{D^\sigma u\}_{|\sigma|=m}=Q\,f,
\end{equation}
\noindent where $Q$ is a bounded operator from $V_p^{-m,a}(\mathbb{R}^n_+)$ 
into $L_p(\mathbb{R}^n_+,\,x_n^{ap}\,dx)$ and $S$ is a linear operator
mapping $L_p(\mathbb{R}^n_+,\,x_n^{ap}\,dx)$ into itself. Furthermore, 
on account of (\ref{ItalianTrick})-(\ref{smallCT}), 
(\ref{altOp})-(\ref{notExceed}) and (\ref{KDu}), we can bound 
$\|S\|_{{\mathfrak L}(L_p(\mathbb{R}^n_+,\,x_n^{ap}\,dx))}$ by 
\begin{equation}\label{SSS}
C\,\Bigl(p\,p'+\frac{1}{s(1-s)}\Bigr) \Bigl(\sum_{|\alpha|=|\beta|=m}
[A_{\alpha\beta}]_{{\rm BMO}(\mathbb{R}^n_+)} 
+\frac{1}{s(1-s)}
\sum_{{|\alpha|+|\beta|<2m}\atop{0\leq |\alpha|,|\beta|\leq m}} 
\|A_{\alpha\beta}\|_{{L^\infty}(\mathbb{R}^n_+)}\Bigr).
\end{equation}
Owing to (\ref{E7})-(\ref{E8}) and with the integral representation formula 
(\ref{IntRepFor}) and the bound (\ref{SSS}) in hand, 
a Neumann series argument and standard functional analysis allow us to 
simultaneously settle the claims (i) and (ii) in the statement of the lemma.  
\hfill$\Box$

\section{The Dirichlet problem in a special Lipschitz domain}
\setcounter{equation}{0}

In this section as well as in subsequent ones, we shall work with 
an unbounded domain of the form 
\begin{equation}\label{10.1.26}
G=\{X=(X',X_n)\in\RR^n:\,X'\in\mathbb{R}^{n-1},\,\,X_n>\varphi(X')\},
\end{equation}
\noindent where $\varphi:\RR^{n-1}\to\RR$ is a Lipschitz function.

\subsection{The space ${\rm BMO}(G)$}

The space of functions of bounded mean oscillations in $G$ can be
introduced in a similar fashion to the case $G=\RR^n_+$. Specifically, 
a locally integrable function on $G$ belongs to the space ${\rm BMO}(G)$ if 
\begin{equation}\label{f-BMO}
[f]_{{\rm BMO}(G)}:=\sup\limits_{(B)}\meanint_{\!\!\!B\cap G}
\Bigl|f(X)-\meanint_{\!\!\!B\cap G}f(Y)\,dY\Bigr|\,dX<\infty,
\end{equation}
\noindent where the supremum is taken over all balls $B$ centered at 
points in ${\bar G}$. Much as before, 
\begin{equation}\label{10.26}
[f]_{{\rm BMO}(G)}\sim\sup\limits_{(B)}\meanint_{\!\!\!B\cap G}
\meanint_{\!\!\!B\cap G}\Bigl|f(X)-f(Y)\Bigr|\,dXdY.
\end{equation}
\noindent This implies the equivalence relation
\begin{equation}\label{1.32}
[f]_{{\rm BMO}(G)}\sim [f\circ\lambda]_{{\rm BMO}(\mathbb{R}^n_+)}
\end{equation}
\noindent for each bi-Lipschitz diffeomorphism $\lambda$ of $\mathbb{R}^n_+$ 
onto $G$. As consequences of definitions, we also have 
\begin{eqnarray}\label{1.33}
[\prod_{1\leq j\leq N} f_j]_{{\rm BMO}(G)}
& \leq & c\,\|f\|^{N-1}_{L_\infty(G)}
[f]_{{\rm BMO}(G)},\quad\mbox{where}\,\,f=(f_1,\ldots, f_N),
\\[6pt]
[f^{-1}]_{{\rm BMO}(G)} & \leq & 
c\,\|f^{-1}\|^2_{L_\infty(G)}[f]_{{\rm BMO}(G)}.
\label{1.34}
\end{eqnarray}

\subsection{A bi-Lipschitz map $\lambda:\mathbb{R}^n_+\to G$ and its inverse}

Let $\varphi:\RR^{n-1}\to\RR$ be the Lipschitz function whose graph is 
$\partial G$ and set $M:=\|\nabla\varphi\|_{L_\infty(\RR^{n-1})}$. 
Next, let $T$ be the extension operator defined as in (\ref{10.1.20}) and, 
for a fixed, sufficiently large constant $C>0$, consider the Lipschitz mapping
\begin{equation}\label{lambda}
\lambda:\,\mathbb{R}^n_+\ni(x',x_n)\mapsto (X',X_n)\in\,G
\end{equation}
\noindent defined by the equalities
\begin{equation}\label{10.1.27}
X':=x',\qquad X_n:=C\,M\,x_n+(T\varphi)(x',x_n)
\end{equation}
\noindent (see {\bf\cite{MS}}, \S{6.5.1} and an earlier, less accessible, 
reference {\bf\cite{MS1}}). The Jacobi matrix of $\lambda$ is given by
\begin{equation}\label{matrix}
\lambda'=
\left( 
\begin{array}{cc}
I & 0 \\
\nabla_{x'}(T\varphi) &  CM+\partial (T\varphi)/\partial x_n
\end{array} 
\right),
\end{equation}
\noindent where $I$ is the identity $(n-1)\times(n-1)$-matrix. Since 
$\vert\partial(T\varphi)/\partial x_n\vert\leq c\,M$ by (\ref{10.1.30}), 
it follows that ${\rm det}\,\lambda'>(C-c)\,M>0$.
Thanks to (\ref{Tfi}) and (\ref{lambda})-(\ref{10.1.27}) we have 
\begin{equation}\label{eq-phi}
X_n-\varphi(X')\sim x_n. 
\end{equation}
\noindent Also, based on (\ref{1.8}) we may write
\begin{equation}\label{1.35}
[\lambda']_{{\rm BMO}(\mathbb{R}^n_+)}
\leq c\,[\nabla\varphi]_{{\rm BMO}(\mathbb{R}^{n-1})}
\end{equation}
\noindent and further, by (\ref{10.1.21}) and (\ref{1.2}), 
\begin{equation}\label{1.4}
\|D^\alpha \lambda'(x)\|_{\mathbb{R}^{n\times n}} 
\leq c(M)\,x_n^{-|\alpha|}[\nabla\varphi]_{{\rm BMO}(\mathbb{R}^{n-1})},
\qquad\forall\,\alpha\,:\,|\alpha|\geq 1. 
\end{equation}
Next, by closely mimicking the proof of Proposition~2.6 from {\bf\cite{MS2}} 
it is possible to show the existence of the inverse Lipschitz  
mapping $\varkappa:=\lambda^{-1}:G\to\RR^n_+$. Owing to (\ref{1.32}), the 
inequality (\ref{1.35}) implies
\begin{equation}\label{1.40}
[\lambda'\circ\varkappa]_{{\rm BMO}(G)}
\leq c\,[\nabla\varphi]_{{\rm BMO}(\mathbb{R}^{n-1})}.
\end{equation}
\noindent Furthermore, (\ref{1.4}) is equivalent to
\begin{equation}\label{1.41}
\|(D^\alpha\lambda')(\varkappa(X))\|_{\mathbb{R}^{n\times n}} 
\leq c(M,\alpha)(X_n-\varphi(X'))^{-|\alpha|} 
[\nabla\varphi]_{{\rm BMO}(\mathbb{R}^{n-1})},
\end{equation}
\noindent whenever $|\alpha|>0$. Since 
$\varkappa'=(\lambda'\circ\varkappa)^{-1}$ we obtain from 
(\ref{1.34}) and (\ref{1.40})
\begin{equation}\label{1.36}
[\varkappa']_{{\rm BMO}(G)}
\leq c\,[\nabla\varphi]_{{\rm BMO}(\mathbb{R}^{n-1})}.
\end{equation}

\noindent On the other hand, using $\varkappa'=(\lambda'\circ\varkappa)^{-1}$ 
and (\ref{1.41}) one can prove by induction on the order of differentiation 
that, for all $X\in G$ and $\alpha\in{\mathbb{N}}_0$ with $|\alpha|>0$, 
\begin{equation}\label{1.5}
\|D^\alpha\varkappa'(X)\|_{\mathbb{R}^{n\times n}}
\leq c(M,\alpha)\,(X_n-\varphi(X'))^{-|\alpha|}
[\nabla\varphi]_{{\rm BMO}(\mathbb{R}^{n-1})}.
\end{equation}

\subsection{The space $V_p^{m,a}(G)$}

Analogously to $V_p^{m,a}(\mathbb{R}^n_+)$, we define the weighted Sobolev 
space $V_p^{m,a}(G)$ naturally associated with the norm
\begin{equation}\label{VVV-space}
\|{\mathcal U}\|_{V_p^{m,a}(G)}:=\Bigl(\sum_{|\gamma|\leq m}\int_{G}
|(X_n-\varphi(X'))^{|\gamma|-m}D^\gamma{\mathcal U}(X)|^p
\,(X_n-\varphi(X'))^{pa}\,dX\Bigr)^{1/p}.
\end{equation}
\noindent Replacing the function $X_n-\varphi(X')$ by either 
$\rho(X):={\rm dist}\,(X,\partial G)$, or by the so-called regularized 
distance function $\rho_{\rm reg}(X)$ 
(defined as on pp.\,170-171 of {\bf\cite{St}}), 
yields equivalent norms on $V_p^{m,a}(G)$. A standard localization 
argument involving a cut-off function vanishing near $\partial G$
(for example, take $\eta(\rho_{\rm reg}/\varepsilon)$ where 
$\eta\in C^\infty_0(\RR)$ satisfies $\eta(t)=0$ for $|t|<1$ and $\eta(t)=1$
for $|t|>2$) shows that 
\begin{equation}\label{Dense-VV}
C^\infty_0(G)\hookrightarrow V_p^{m,a}(G)\quad\mbox{densely}.
\end{equation}
\noindent Next, we observe that for each ${\mathcal U}\in C^\infty_0(G)$, 
\begin{equation}\label{equivNr}
C\,s\,\|{\mathcal U}\|_{V_p^{m,a}(G)}\leq \Bigl(\sum_{|\gamma|=m}\int_G 
|D^\gamma {\mathcal U}(X)|^p\,(X_n-\varphi(X'))^{pa}\,dX\Bigr)^{1/p}
\leq\,\|{\mathcal U}\|_{V_p^{m,a}(G)},
\end{equation}
\noindent where, as before, $s=1-a-1/p$. 
Indeed, for each multi-index $\gamma$ with $|\gamma|\leq m$, 
the one-dimensional Hardy's inequality gives 
\begin{eqnarray}\label{Hardy}
&& \int_{G}|(X_n-\varphi(X'))^{|\gamma|-m}D^\gamma{\mathcal U}(X)|^p\,
(X_n-\varphi(X'))^{pa}\,dX 
\nonumber\\[6pt]
&&\qquad\qquad\qquad\qquad
\leq\bigl({C}/{s}\bigr)^p 
\sum_{|\alpha|=m}\int_G|D^\alpha{\mathcal U}(X)|^p\,(X_n-\varphi(X'))^{pa}\,dX,
\end{eqnarray}
\noindent and the first inequality in (\ref{equivNr}) follows readily from it. 
Also, the second inequality in (\ref{equivNr}) is a trivial consequence 
of (\ref{VVV-space}). Going further, we aim to establish that 
\begin{equation}\label{equivNr-bis}
c_1\,\|{u}\|_{V_p^{m,a}(\mathbb{R}^n_+)}
\leq\|{u}\circ\varkappa\|_{V_p^{m,a}(G)}
\leq c_2\,\|{u}\|_{V_p^{m,a}(\mathbb{R}^n_+)},
\end{equation}
\noindent where $c_1$ and $c_2$ do not depend on $p$ and $s$, and 
$\varkappa:G\to\mathbb{R}^n_+$ is the map introduced in \S{5.2}. 
Clearly, it suffices to prove the upper estimate for 
$\|{u}\circ\varkappa\|_{V_p^{m,a}(G)}$ in (\ref{equivNr-bis}). 
As a preliminary matter, we note that 
\begin{eqnarray}\label{1.49}
D^\gamma\bigl({u}(\varkappa(X))\bigr) 
& = & \bigl((\varkappa'^*(X)\xi)_{\xi=D}^\gamma\,{u}\bigr)(\varkappa(X))
\nonumber\\[6pt]
&& +\sum_{1\leq|\tau|<|\gamma|}
(D^\tau{u})(\varkappa(X))\sum_\sigma\,c_\sigma\prod_{i=1}^n\prod_j 
D^{\sigma_{ij}}\varkappa_i(X),
\end{eqnarray}
\noindent where
\begin{equation}\label{1.50}
\sigma=(\sigma_{ij}),\quad
\sum_{i,j}\sigma_{ij}=\gamma,\quad |\sigma_{ij}|\geq 1,\quad 
\sum_{i,j}(|\sigma_{ij}|-1)=|\gamma|-|\tau|.
\end{equation}
\noindent In turn, (\ref{1.49})-(\ref{1.50}) and (\ref{1.5}) 
allow us to conclude that 
\begin{equation}\label{DUU}
|D^\gamma\bigl({u}(\varkappa(X))\bigr)|\leq c\sum_{1\leq |\tau|\leq |\gamma|}
x_n^{|\tau|-|\gamma|}\,|D^\tau{u}(x)|, 
\end{equation}
\noindent which, in view of (\ref{eq-phi}), yields the desired conclusion. 
Finally, if, as usual, $p'=p/(p-1)$, we set 
\begin{equation}\label{dual-VG}
V^{-m,a}_p(G):=\Bigl(V^{m,-a}_{p'}(G)\Bigr)^*.
\end{equation}

\subsection{Solvability and regularity result for the Dirichlet 
problem in the domain $G$}

Let us consider the differential operator 
\begin{equation}\label{E4}
{\mathcal L}\,{\mathcal U}
={\mathcal L}(X,D_X)\,{\mathcal U}=\sum_{|\alpha|=|\beta|=m} 
D^\alpha(\mathfrak{A}_{\alpha\beta}(X)\,D^\beta{\mathcal U}),\qquad X\in G,
\end{equation}
\noindent whose matrix-valued coefficients satisfy
\begin{equation}\label{E4a}
\sum_{|\alpha|=|\beta|=m}\|\mathfrak{A}_{\alpha\beta}\|_{L_\infty(G)} 
\leq \kappa^{-1}.
\end{equation}
\noindent This operator generates the sesquilinear form  
${\mathcal L}(\cdot,\cdot):V_p^{m,a}(G)\times V_{p'}^{m,-a}(G)\to\CC$
where $p'$ is the conjugate exponent of $p$, defined by 
\begin{equation}\label{LUV}
{\mathcal L}({\mathcal U},{\mathcal V}):=\sum_{|\alpha|=|\beta|=m}
\int_G\langle\mathfrak{A}_{\alpha\beta}(X)\,
D^\beta{\mathcal U}(X),\,D^\alpha{\mathcal V}(X)\rangle\,dX. 
\end{equation}
\noindent We assume that 
\begin{equation}\label{B25}
\Re\,{\mathcal L}({\mathcal U},{\mathcal U})
\geq \kappa\sum_{|\gamma|=m}\|D^\gamma\,{\mathcal U}\|^2_{L_2(G)}
\quad\mbox{ holds for all }\,\,{\mathcal U}\in V_2^{m,0}(G).
\end{equation}
\begin{lemma}\label{lem5a}
{\rm (i)} Let $p\in (1,\infty)$, $1/p+1/p'=1$, $-1/p<a<1-1/p$, and 
set $s:=1-a-1/p$. Also, suppose that
\begin{equation}\label{E7a}
[\nabla\varphi]_{{\rm BMO}(\mathbb{R}^{n-1})}
+\sum_{|\alpha|=|\beta|=m}[{\mathfrak A}_{\alpha\beta}]_{{\rm BMO}(G)}
\leq \delta,
\end{equation}
\noindent where $\delta$ satisfies
\begin{equation}\label{E8b}
\Bigl(p\,p'+\frac{1}{s(1-s)}\Bigr)\,\frac{\delta}{s(1-s)}
<C(n,m,\kappa,\|\nabla\varphi\|_{L_\infty(\mathbb{R}^{n-1})})
\end{equation}
\noindent with a sufficiently small constant $C$, independent of $p$ and $s$. 
Then 
\begin{equation}\label{Liso-2}
{\mathcal L}(X,D_X):V_p^{m,a}(G)\longrightarrow V_p^{-m,a}(G)
\quad\mbox{ isomorphically}. 
\end{equation}
{\rm (ii)} Suppose that $p_i\in (1,\infty)$, $-1/p_i<a_i<1-1/p_i$, $i=1,2$. 
and that {\rm (\ref{E8})} holds with $p_i$ and $s_i:=1-a_i-1/p_i$ in place 
of $p$ and $s$. If ${\mathcal U}\in V_{p_1}^{m,a_1}(G)$ and 
${\mathcal L}\,{\mathcal U}\in V_{p_1}^{-m,a_1}(G)\cap V_{p_2}^{-m,a_2}(G)$, 
then ${\mathcal U}\in V_{p_2}^{m,a_2}(G)$.
\end{lemma}
\noindent{\bf Proof.} We shall extensively use the flattening mapping 
$\lambda$ and its inverse studied in \S{5.2}. The assertions (i) and (ii)  
will follow directly from Lemma~\ref{lem5} as soon as we show that 
the operator 
\begin{equation}\label{E20}
L({\mathcal U}\circ\lambda):=({\mathcal L}\,{\mathcal U})\circ\lambda
\end{equation}
\noindent satisfies all the hypotheses in that lemma. The sesquilinear form 
corresponding to the operator $L$ will be denoted by $L(u,v)$.
Set $u(x):={\mathcal U}(\lambda(x))$, $v(x):={\mathcal V}(\lambda(x))$ and
note that (\ref{1.49}) implies
\begin{equation}\label{E40}
D^\beta {\mathcal U}(X)=\bigl((\varkappa'^*(\lambda(x))\xi)_{\xi=D}^\beta\,{u}
\bigr)(x)+\sum_{1\leq|\tau|<|\beta|}K_{\beta\tau}(x)\,x_n^{|\tau|-|\beta|} 
D^\tau u(x),
\end{equation}
\begin{equation}\label{E41}
D^\alpha {\mathcal V}(X)
=\bigl((\varkappa'^*(\lambda(x))\xi)_{\xi=D}^\alpha\,{v}
\bigr)(x)+\sum_{1\leq|\tau|<|\alpha|}K_{\alpha\tau}(x)\,x_n^{|\tau|-|\alpha|} 
D^\tau v(x),
\end{equation}
\noindent where, thanks to (\ref{1.5}), the coefficients $K_{\gamma\tau}$ satisfy
\begin{equation}\label{E42}
\|K_{\gamma\tau}\|_{L_\infty(\mathbb{R}^n_+)} 
\leq c\,[\nabla\varphi]_{{\rm BMO}(\mathbb{R}^{n-1})}.
\end{equation}
\noindent Plugging (\ref{E40}) and (\ref{E41}) into the definition of 
${\mathcal L}({\mathcal U},{\mathcal V})$, we arrive at
\begin{equation}\label{LUV2}
{\mathcal L}({\mathcal U},{\mathcal V})=L_0(u,v)
+\sum_{{1\leq|\alpha|,|\beta|\leq m}\atop{|\alpha|+|\beta|<2m}}
\int_{\mathbb{R}^n_+}\langle {A}_{\alpha\beta}(x)\,x_n^{|\alpha|+|\beta|-2m} 
D^\beta\,u(x),\,D^\alpha v(x)\rangle\,dx,
\end{equation}
\noindent where
\begin{equation}\label{Lzero}
L_0(u,v)=\sum_{|\alpha|=|\beta|=m}\int_{\mathbb{R}^n_+} 
\langle (\mathfrak{A}_{\alpha\beta}\circ\lambda)
((\varkappa'^*\circ\lambda)\xi)_{\xi=D}^\beta\,{u},\, 
((\varkappa'^*\circ\lambda)\xi)_{\xi=D}^\alpha\,{v}\rangle
\,{\rm det}\,\lambda'\,dx.
\end{equation}
\noindent It follows from (\ref{E40})-(\ref{E42}) that the coefficient 
matrices $A_{\alpha\beta}$ obey 
\begin{equation}\label{E43}
\sum_{{1\leq|\alpha|,|\beta|\leq m}\atop{|\alpha|+|\beta|<2m}}
\|A_{\alpha\beta}\|_{L_\infty(\mathbb{R}^n_+)} 
\leq {c}{\kappa}^{-1}\,[\nabla\varphi]_{{\rm BMO}(\mathbb{R}^{n-1})},
\end{equation}
\noindent where $c$ depends on $m$, $n$, and 
$\|\nabla\varphi\|_{L_\infty(\mathbb{R}^{n-1})}$. 
We can write the form $L_0(u,v)$ as
\begin{equation}\label{Lzero-uv}
\sum_{|\alpha| =|\beta| =m}\int_{\mathbb{R}^n_+}
\langle {A}_{\alpha\beta}(x)\,D^\beta u(x),\,D^\alpha v(x)\rangle\,dx,
\end{equation}
\noindent where the coefficient matrices $A_{\alpha\beta}$ are given by
\begin{equation}\label{Aalbet}
A_{\alpha\beta}={\rm det}\,\lambda'\,\sum_{|\gamma|=|\tau|=m}
P_{\alpha\beta}^{\gamma\tau}(\varkappa'\circ\lambda)
({\mathfrak A}_{\gamma\tau}\circ\lambda),
\end{equation}
\noindent for some scalar homogeneous polynomials 
$P_{\alpha\beta}^{\gamma\tau}$ of the elements of the matrix 
$\varkappa'(\lambda(x))$ of degree $2m$. 
In view of (\ref{1.33})-(\ref{1.36}),  
\begin{equation}\label{E44}
\sum_{|\alpha|=|\beta|=m}[A_{\alpha\beta}]_{{\rm BMO}(\mathbb{R}^n_+)} 
\leq c\Bigl(\kappa^{-1}[\nabla\varphi]_{{\rm BMO}(\mathbb{R}^{n-1})} 
+\sum_{|\alpha|=|\beta|=m}[{\mathfrak A}_{\alpha\beta}]_{{\rm BMO}(G)}\Bigr),
\end{equation}
\noindent where $c$ depends on $n$, $m$, and 
$\|\nabla\varphi\|_{L_\infty(\mathbb{R}^{n-1})}$. By (\ref{E43}) 
\begin{equation}\label{L-L}
|L(u,u)-L_0(u,u)|\leq c\,\delta\|u\|_{V_2^{m,0}(\mathbb{R}^n_+)}^2
\end{equation}
\noindent and, therefore,
\begin{equation}\label{ReLzero}
\Re\,L_0(u,u)\geq\Re\,{\mathcal L}({\mathcal U},{\mathcal U})
-c\,\delta\|u\|^2_{V_2^{m,0}(\mathbb{R}^n_+)}.
\end{equation}
\noindent Using (\ref{B25}) and the equivalence 
(cf. the discussion in \S{5.3})
\begin{equation}\label{norm-U}
\|{\mathcal U}\|_{V_2^{m,0}(G)}\sim\|u\|_{V_2^{m,0}(\mathbb{R}^n_+)},
\end{equation}
\noindent we arrive at (\ref{B5}). Thus, all conditions of Lemma~\ref{lem5} 
hold and the result follows. 
\hfill$\Box$

\section{The Dirichlet problem in a bounded Lipschitz domain}
\setcounter{equation}{0}

\subsection{Background}
Let $\Omega$ be a {\it bounded Lipschitz domain} in $\RR^n$ which means
(cf. {\bf\cite{St}}, p.\,189)
that there exists a finite open covering $\{{\mathcal O}_j\}_{1\leq j\leq N}$
of $\partial\Omega$ with the property that, for every $j\in\{1,...,N\}$,
${\mathcal O}_j\cap\Omega$ coincides with the portion of ${\mathcal O}_j$ 
lying in the over-graph of a Lipschitz function $\varphi_j:\RR^{n-1}\to\RR$
(where $\RR^{n-1}\times\RR$ is a new system of coordinates obtained from 
the original one via a rigid motion). We then define the 
{\it Lipschitz constant} of a bounded Lipschitz domain $\Omega\subset\RR^n$ as
\begin{equation}\label{Lip-ct}
\inf\,\Bigl(\max\{\|\nabla\varphi_j\|_{L_\infty(\RR^{n-1})}:\,1\leq j\leq N\}
\Bigr),
\end{equation}
\noindent where the infimum is taken over all possible families 
$\{\varphi_j\}_{1\leq j\leq N}$ as above. 
Such domains are called {\it minimally smooth} in E.\,Stein's book 
{\bf{\cite{St}}}.
It is a classical result that, for a Lipschitz domain $\Omega$, the surface 
measure $d\sigma$ is well-defined on $\partial\Omega$ and an outward pointing 
normal vector $\nu$ exists a.e. on $\partial\Omega$. 
 
We denote by $\rho(X)$ the distance from $X\in\RR^n$ to $\partial\Omega$
and, for $p$, $a$ and $m$ as in (\ref{indices}), introduce the weighted 
Sobolev space $V^{m,a}_p(\Omega)$ naturally associated with the norm
\begin{equation}\label{normU2}
\|{\mathcal U}\|_{V_p^{m,a}(\Omega)}
:=\Bigl(\sum_{0\leq |\beta|\leq m}\int_{\Omega}
|\rho(X)^{|\beta|-m} D^\beta{\mathcal U}(X)|^p\,\rho(X)^{pa}\,dX\Bigr)^{1/p}. 
\end{equation}
\noindent One can check the equivalence of the norms
\begin{equation}\label{equiv-Nr2}
\|{\mathcal U}\|_{V_p^{m,a}(\Omega)}\sim 
\|\rho_{\rm reg}^a\,{\mathcal U}\|_{V_p^{m,0}(\Omega)}, 
\end{equation}
\noindent where $\rho_{\rm reg}(X)$ stands for the regularized distance 
from $X$ to $\partial\Omega$ (in the sense of Theorem~2, p.\,171 in 
{\bf\cite{St}}). Much as with (\ref{Dense-VV}), it is also easily proved that 
$C_0^\infty(\Omega)$ is dense in $V^{m,a}_p(\Omega)$ and that 
\begin{equation}\label{sumDU}
\|{\mathcal U}\|_{V_p^{m,a}(\Omega)}\sim 
\Bigl(\sum_{|\beta|=m}\int_\Omega|D^\beta{\mathcal U}(X)|^p\,\rho(X)^{pa}\,dX
\Bigr)^{1/p}
\end{equation}
\noindent uniformly for ${\mathcal U}\in C_0^\infty(\Omega)$. 
As in (\ref{dual-VG}), we set 
\begin{equation}\label{dual-V}
V^{-m,a}_p(\Omega):=\Bigl(V^{m,-a}_{p'}(\Omega)\Bigr)^*.
\end{equation}

Let us fix a Cartesian coordinates system and consider 
the differential operator 
\begin{equation}\label{E444}
{\mathcal A}\,{\mathcal U}={\mathcal A}(X,D_X)\,{\mathcal U}
:=\sum_{|\alpha|=|\beta|=m}D^\alpha({\mathcal A}_{\alpha\beta}(X)
\,D^\beta{\mathcal U}),\qquad X\in\Omega,
\end{equation}
\noindent with measurable $l\times l$ matrix-valued coefficients. 
The corresponding sesquilinear form will be denoted by 
${\mathcal A}({\mathcal U},{\mathcal V})$. 
Similarly to (\ref{E4a}) and (\ref{B25}) we impose the conditions
\begin{equation}\label{E4b}
\sum_{|\alpha|=|\beta|=m}\|{\mathcal A}_{\alpha\beta}\|_{L_\infty(\Omega)} 
\leq \kappa^{-1}
\end{equation}
\noindent and
\begin{equation}\label{B25b}
\Re\,{\mathcal A}({\mathcal U},{\mathcal U})\geq\kappa\sum_{|\gamma|=m}
\|D^\gamma\,{\mathcal U}\|^2_{L_2(\Omega)}\quad\mbox{for all}\,\,\,
{\mathcal U}\in V_2^{m,0}(\Omega).
\end{equation}

\subsection{Interior regularity of solutions}

\begin{lemma}\label{lem2}
Let $\Omega\subset\RR^n$ be a bounded Lipschitz domain. 
Pick two functions ${\mathcal H},{\mathcal Z}\in C^\infty_0(\Omega)$ 
such that ${\mathcal H}\,{\mathcal Z}={\mathcal H}$, and assume that, 
for a sufficiently small constant $c(m,n,\kappa)>0$, 
\begin{equation}\label{E5}
\sum_{|\alpha|=|\beta|=m}[{\mathcal A}_{\alpha\beta}]_{{\rm BMO}(\Omega)}
\leq\delta,
\end{equation}
\noindent where
 \begin{equation}\label{E6}
\delta\leq \frac{c(m,n,\kappa)}{p\,p'}
\end{equation}
If ${\mathcal U}\in W_q^m(\Omega,loc)$ for a certain $q<p$ and 
${\mathcal A}\,{\mathcal U}\in W_p^{-m}(\Omega,loc)$, then 
${\mathcal U}\in W_p^m(\Omega,loc)$ and
\begin{equation}\label{E3}
\|{\mathcal H}\,{\mathcal U}\|_{W_p^m(\Omega)}
\leq C\,(\|{\mathcal H}\,{\mathcal A}(\cdot,D)\,
{\mathcal U}\|_{W_p^{-m}(\Omega)}
+\|{\mathcal Z}\,{\mathcal U}\|_{W_q^m(\Omega)}).
\end{equation}
\end{lemma}
\noindent{\bf Proof.} We shall use the notation 
${\mathcal A}_Y$ for the operator ${\mathcal A}(Y,D_X)$, where $Y\in\Omega$ 
and the notation $\Phi_Y$ for a fundamental solution of ${\mathcal A}_Y$ 
in $\mathbb{R}^n$. Then, with star denoting the convolution product, 
\begin{equation}\label{HUPhi}
{\mathcal H}\,{\mathcal U}
+\Phi_Y\ast({\mathcal A}-{\mathcal A}_Y)({\mathcal H}\,{\mathcal U}) 
=\Phi_Y\ast({\mathcal H}\,{\mathcal A}{\mathcal U})
+\Phi_Y\ast ([{\mathcal A},{\mathcal H}]({\mathcal Z}\,{\mathcal U}))
\end{equation}
\noindent and, consequently, for each multi-index $\gamma$, $|\gamma|=m$, 
\begin{eqnarray}\label{DHU}
&& D^\gamma({\mathcal H}\,{\mathcal U})+\sum_{|\alpha|=|\beta|=m} 
D^{\alpha+\gamma}\Phi_Y\ast\bigl(({\mathcal A}_{\alpha\beta}
-{\mathcal A}_{\alpha\beta}(Y)) D^\beta({\mathcal H}\,{\mathcal U})\bigr)
\nonumber\\[6pt]
&& \qquad\quad
=D^\gamma\Phi_Y\ast({\mathcal H}\,{\mathcal A}{\mathcal U})
+D^\gamma\Phi_Y\ast([{\mathcal A},{\mathcal H}]({\mathcal Z}\,{\mathcal U})).
\end{eqnarray}
\noindent Writing this equation at the point $Y$ and using (\ref{b1}), 
we obtain
\begin{eqnarray}\label{E3a}
&& (1-C\,p\,p'\delta)
\sum_{|\gamma|=m}\|D^\gamma({\mathcal H}\,{\mathcal U})\|_{L_p(\Omega)}
\nonumber\\[6pt]
&& \qquad\quad
\leq C(p,\kappa)
(\|{\mathcal H}\,{\mathcal A}{\mathcal U}\|_{W_p^{-m}(\Omega)} 
+\|[{\mathcal A},{\mathcal H}]({\mathcal Z}\,
{\mathcal U})\|_{W_p^{-m}(\Omega)}).
\end{eqnarray}
\noindent Let $p'<n$. We have for every ${\mathcal V}\in\ring{W}^m_p(\Omega)$,
the closure of $C^\infty_0(\Omega)$ in $W^m_p(\Omega)$,
\begin{eqnarray}\label{AHZUV}
&& \Bigl|\int_\Omega\langle[{\mathcal A},{\mathcal H}]
({\mathcal Z}\,{\mathcal U}),{\mathcal V}\rangle\,dX\Bigr|
=|{\mathcal A}({\mathcal H}{\mathcal Z}\,{\mathcal U},{\mathcal V})
-{\mathcal A}({\mathcal Z}\,{\mathcal U},{\mathcal H}{\mathcal V})| 
\nonumber\\[6pt]
&& \qquad\qquad\quad
\leq c(\|{\mathcal Z}\,{\mathcal U}\|_{W_p^{m-1}(\Omega)} 
\|{\mathcal V}\|_{W_{p'}^m(\Omega)}
+\|{\mathcal Z}\,{\mathcal U}\|_{W_{\frac{pn}{n+p}}^{m} 
(\Omega)}\|{\mathcal V}\|_{W_{\frac{p'n}{n-p'}}^{m-1}(\Omega)}).
\end{eqnarray}
\noindent By Sobolev's theorem
\begin{equation}\label{ZZU}
\|{\mathcal Z}\,{\mathcal U}\|_{W_p^{m-1}(\Omega)}
\leq c\,\|{\mathcal Z}\,{\mathcal U}\|_{W_{\frac{pn}{n+p}}^{m}(\Omega)} 
\end{equation}
\noindent and
\begin{equation}\label{VWpm}
\|{\mathcal V}\|_{W_{\frac{p'n}{n-p'}}^{m-1}(\Omega)} 
\leq c\,\|{\mathcal V}\|_{W_{p'}^m(\Omega)}.
\end{equation}
\noindent Therefore,
\begin{equation}\label{AHZ}
\Bigl|\int_\Omega\langle
[{\mathcal A},{\mathcal H}]({\mathcal Z}\,{\mathcal U}),{\mathcal V}\rangle 
\,dX\Bigr|
\leq c\,\|{\mathcal Z}\,{\mathcal U}\|_{W_{\frac{pn}{n+p}}^{m}(\Omega)} 
\|{\mathcal V}\|_{W_{p'}^m(\Omega)}
\end{equation}
\noindent which is equivalent to the inequality
\begin{equation}\label{AHZ3}
\|[{\mathcal A},{\mathcal H}]({\mathcal Z}\,{\mathcal U})\|_{W_p^{-m}(\Omega)}
\leq c\,\|{\mathcal Z}\,{\mathcal U}\|_{W_{\frac{pn}{n+p}}^{m} (\Omega)}. 
\end{equation}
\noindent In the case $p'\geq n$, the same argument leads to a 
similar inequality, where $pn/(n+p)$ is replaced by $1+\varepsilon$ with 
an arbitrary $\varepsilon>0$ for $p'>n$ and $\varepsilon=0$ for $p'=n$.
Now, (\ref{E3}) follows from (\ref{E3a}) if $p'\geq n$ and $p'<n$, 
$q\geq pn/(n+p)$. In the remaining case the goal is achieved by iterating 
this argument finitely many times. 
\hfill$\Box$
\begin{corollary}\label{cor3}
Let $p\geq 2$ and suppose that {\rm (\ref{E5})} and {\rm (\ref{E6})} hold. 
If ${\mathcal U}\in W_2^m(\Omega,loc)$ and 
${\mathcal A}\,{\mathcal U}\in W_p^{-m}(\Omega,loc)$, then 
${\mathcal U}\in W_p^m(\Omega,loc)$ and
\begin{equation}\label{E3-aaa}
\|{\mathcal H}\,{\mathcal U}\|_{W_p^m(\Omega)}
\leq C\,(\|{\mathcal Z}\,{\mathcal A}(\cdot,D)\, 
{\mathcal U}\|_{W_p^{-m}(\Omega)}
+\|{\mathcal Z}\,{\mathcal U}\|_{W_2^{m-1}(\Omega)}).
\end{equation}
\end{corollary}
\noindent{\bf Proof.} Let ${\mathcal Z}_0$ denote a real-valued function 
in $C^\infty_0(\Omega)$ such that ${\mathcal H}{\mathcal Z}_0={\mathcal H}$ 
and ${\mathcal Z}_0{\mathcal Z}={\mathcal Z}_0$. By (\ref{E3})
\begin{equation}\label{E3b}
\|{\mathcal H}\,{\mathcal U}\|_{W_p^m(\Omega)}
\leq C\,(\|{\mathcal H}\,{\mathcal A}(\cdot,D)\, 
{\mathcal U}\|_{W_p^{-m}(\Omega)}
+\|{\mathcal Z}_0\,{\mathcal U}\|_{W_2^{m}(\Omega)})
\end{equation}
\noindent and it follows from (\ref{B25b}) that
\begin{equation}\label{ZUW}
\|{\mathcal Z}_0\,{\mathcal U}\|^2_{W_2^m(\Omega)}\leq c\,\kappa^{-1}
\Re\,{\mathcal A}({\mathcal Z}_0\,{\mathcal U},{\mathcal Z}_0\,{\mathcal U}).
\end{equation}
\noindent Furthermore,
\begin{equation}\label{AZ5}
|{\mathcal A}({\mathcal Z}_0\,{\mathcal U},{\mathcal Z}_0\,{\mathcal U})
-{\mathcal A}({\mathcal U},{\mathcal Z}_0^2\,{\mathcal U})|\leq c\,\kappa^{-1}
\|{\mathcal Z}\,{\mathcal U}\|_{W_2^{m-1}(\Omega)}\,
\|{\mathcal Z}_0\,{\mathcal U}\|_{W_2^{m}(\Omega)}.
\end{equation}
\noindent Hence
\begin{equation}\label{ZU6}
\|{\mathcal Z}_0\,{\mathcal U}\|^2_{W_2^m(\Omega)}\leq c\,\kappa^{-1}
(\|{\mathcal Z}\,{\mathcal A}{\mathcal U}\|_{W_2^{-m}(\Omega)}\, 
\|{\mathcal Z}_0^2\,{\mathcal U}\|_{W_2^m(\Omega)}
+\kappa^{-1}\|{\mathcal Z}\,{\mathcal U}\|_{W_2^{m-1}(\Omega)}\,
\|{\mathcal Z}_0\,{\mathcal U}\|_{W_2^{m}(\Omega)})
\end{equation}
\noindent and, therefore, 
\begin{equation}\label{Zu7}
\|{\mathcal Z}_0\,{\mathcal U}\|_{W_2^m(\Omega)}
\leq c\,\kappa^{-1}(\|{\mathcal Z}\,{\mathcal A} 
{\mathcal U}\|_{W_2^{-m}(\Omega)}\,
+\kappa^{-1}\|{\mathcal Z}\,{\mathcal U}\|_{W_2^{m-1}(\Omega)}). 
\end{equation}
\noindent Combining this inequality with (\ref{E3b}) we arrive 
at (\ref{E3-aaa}). 
\hfill$\Box$

\subsection{Invertibility of 
${\mathcal A}:V_p^{m,a}(\Omega)\longrightarrow V_p^{-m,a}(\Omega)$}

Recall the quantities defined in (\ref{e60})-(\ref{e61}).
\begin{theorem}\label{th1a}
Let  $1<p<\infty$, $0<s<1$, and set $a=1-s-1/p$. Furthermore, let $\Omega$ be 
a bounded Lipschitz domain in $\mathbb{R}^n$. Suppose that the differential
operator ${\mathcal A}$ is as in {\rm \S{6.1}} and that, in addition, 
\begin{equation}\label{E16}
\sum_{|\alpha|=|\beta|=m}\{{\mathcal A}_{\alpha\beta}\}_{\ast,\Omega}
+\{\nu\}_{\ast,\partial\Omega}\leq\delta,
\end{equation}
\noindent where
\begin{equation}\label{E17}
\Bigl(p\,p'+\frac{1}{s(1-s)}\Bigr)\frac{\delta}{s(1-s)}\leq c
\end{equation}
\noindent for a sufficiently small constant $c>0$ independent of $p$ and $s$. 
Then the operator
\begin{equation}\label{cal-L}
{\mathcal A}:V_p^{m,a}(\Omega)\longrightarrow V_p^{-m,a}(\Omega)
\quad\mbox{ isomorphically}. 
\end{equation}
\end{theorem}
\noindent{\bf Proof.} We shall proceed in a series of steps starting with
\vskip 0.08in
(i) {\it The construction of the auxiliary domain $G$ and 
operator ${\mathcal L}$}. Let $\varepsilon$ be small enough so that 
\begin{equation}\label{m1}
\sum_{|\alpha|=|\beta|=m}\meanint_{\!\!\!B_r\cap \Omega}
\meanint_{\!\!\!B_r\cap\Omega}
|{\mathcal A}_{\alpha\beta}(X)-{\mathcal A}_{\alpha\beta}(Y)|\,dXdY
\leq 2\delta
\end{equation}
\noindent for all balls in $\{B_r\}_\Omega$ with radii $r<\varepsilon$ and 
\begin{equation}\label{m2}
\meanint_{\!\!\!B_r\cap\partial\Omega}\meanint_{\!\!\!B_r\cap\partial\Omega}\,
\Bigl|\nu(X)-\nu(Y)\,\Bigr|\,d\sigma_Xd\sigma_Y\leq 2\delta
\end{equation}
\noindent for all balls in $\{B_r\}_{\partial\Omega}$ with radii 
$r<\varepsilon$. 
We fix a ball $B_\varepsilon$ in $(B_\varepsilon)_{\partial\Omega}$ 
and assume without loss of generality that, in a suitable system of 
Cartesian coordinates,  
\begin{equation}\label{newGGG}
\Omega\cap B_\varepsilon=\{X=(X',X_n)\in B_\varepsilon:\,X_n>\varphi(X')\}
\end{equation}
\noindent for some Lipschitz function $\varphi:\RR^{n-1}\to\RR$. 
Consider now the unique cube $Q(\varepsilon)$ (relative to this system of 
coordinates) which is inscribed in $B_\varepsilon$ and denote its projection 
onto $\mathbb{R}^{n-1}$ by $Q'(\varepsilon)$. 
Since $\nabla\varphi=-\nu'/\nu_n$, it follows from (\ref{m2}) that 
\begin{equation}\label{m3}
\meanint_{\!\!\!B'_r}\meanint_{\!\!\!B'_r}\,
\Bigl|\nabla\varphi(X')-\nabla\varphi(Y')\,\Bigr|\,dX'dY'\leq c(n)\,\delta,
\end{equation}
\noindent where $B'_r=B_r\cap \mathbb{R}^{n-1}$, $r<\varepsilon$. Let us 
retain the notation $\varphi$ for the mirror extension of the function 
$\varphi$ from $Q'(\varepsilon)$ onto $\mathbb{R}^{n-1}$.  
We extend ${\mathcal A}_{\alpha\beta}$ from $Q(\varepsilon)\cap\Omega$ onto 
$Q(\varepsilon)\backslash\Omega$ by setting
\begin{equation}\label{A=A}
{\mathcal A}_{\alpha\beta}(X)
:={\mathcal A}_{\alpha\beta}(X',-X_n+2\varphi(X')),
\qquad X\in Q(\varepsilon)\backslash\Omega, 
\end{equation}
\noindent and we shall use the notation ${\mathfrak A}_{\alpha\beta}$ for 
the periodic extension of ${\mathcal A}_{\alpha\beta}$ from $Q(\varepsilon)$ 
onto $\mathbb{R}^n$. 

Consistent with the earlier discussion in Section~5, we shall denote the 
special 
Lipschitz domain $\{X=(X',X_n):\,X'\in\mathbb{R}^{n-1},\,X_n>\varphi(X')\}$ by
$G$. One can easily see that, owing to the $2\varepsilon n^{-1/2}$-periodicity
of $\varphi$ and ${\mathcal A}_{\alpha\beta}$, 
\begin{equation}\label{SumA}
\sum_{|\alpha|=|\beta|=m}[{\mathcal A}_{\alpha\beta}]_{{\rm BMO}(G)}
+[\nabla\varphi]_{{\rm BMO}(\mathbb{R}^{n-1})}\leq c(n)\,\delta.
\end{equation}
\noindent With the operator ${\mathcal A}(X,D_X)$ in $\Omega$, we 
associate the auxiliary operator ${\mathcal L}(X,D_X)$ in $G$ 
given by (\ref{E4}). 

\vskip 0.08in
(ii) {\it Uniqueness.} Assuming that ${\mathcal U}\in V_p^{m,a}(\Omega)$ 
satisfies ${\mathcal L}\,{\mathcal U}=0$ in $\Omega$, we shall show that 
${\mathcal U}\in V_2^{m,0}(\Omega)$. This will imply that ${\mathcal U}=0$ 
which proves the injectivity of the operator in (\ref{cal-L}). 
To this end, pick a function ${\mathcal H}\in C_0^\infty(Q(\varepsilon))$ 
and write ${\mathcal L}({\mathcal H}\,{\mathcal U})
=[{\mathcal L},\,{\mathcal H}]\,{\mathcal U}$. 
Also, fix a small $\theta>0$ and select a smooth function $\Lambda$ on 
$\mathbb{R}^1_+$, which is identically $1$ on $[0,1]$ and which 
vanishes identically on $(2,\infty)$. Then by (ii) in Lemma~\ref{lem5a},
\begin{equation}\label{E19}
{\mathcal L}({\mathcal H}\,{\mathcal U})-[{\mathcal L},\,{\mathcal H}]\,
(\Lambda(\rho_{\rm reg}/\theta)\,{\mathcal U})
\in V_2^{-m,0}(G)\cap V_p^{-m,a}(G).
\end{equation}
\noindent Note that the operator 
\begin{equation}\label{LH-1}
[{\mathcal  L},\,{\mathcal H}]\rho_{\rm reg}^{-1}: V_p^{m,a}(G) 
\longrightarrow V_p^{-m,a}(G)
\end{equation}
\noindent is bounded and that the norm of the multiplier 
$\rho_{\rm reg}\,\Lambda(\rho_{\rm reg}/\theta)$ in $V_p^{m,a}(G)$ 
is $O(\theta)$. Moreover, the same is true for $p=2$ and $a=0$. 
The inclusion (\ref{E19}) can be written in the form
\begin{equation}\label{E21}
{\mathcal L}({\mathcal H}\,{\mathcal U})
+{\mathcal M}({\mathcal Z}\,{\mathcal U})\in V_p^{-m,a}(G)\cap V_2^{-m,0}(G),
\end{equation}
\noindent where ${\mathcal Z}\in C^\infty_0(\mathbb{R}^n)$, 
${\mathcal Z}\,{\mathcal H}={\mathcal H}$ and ${\mathcal M}$ is a 
linear operator mapping 
\begin{equation}\label{V2Vp}
V_p^{m,a}(G)\to V_p^{-m,a}(G)\quad {\rm and}\quad  
V_2^{m,0}(G)\to V_2^{-m,0}(G)
\end{equation}
\noindent with both norms of order $O(\theta)$. 

Select a finite covering of $\overline{\Omega}$ by cubes $Q_j(\varepsilon)$ 
and let $\{{\mathcal H}_j\}$ be a smooth partition of unity subordinate to 
$\{Q_j(\varepsilon)\}$. 
Also, let ${\mathcal Z}_j\in C_0^\infty(Q_j(\varepsilon))$ 
be such that ${\mathcal H}_j{\mathcal Z}_j={\mathcal H}_j$. 
By $G_j$ we denote the special Lipschitz domain generated by the 
cube $Q_j(\varepsilon)$ as in part (i) of the present proof. 
The corresponding operators ${\mathcal L}$ and ${\mathcal M}$ will 
be denoted by ${\mathcal L}_j$ and ${\mathcal M}_j$, respectively. 
It follows from (\ref{E21}) that 
\begin{equation}\label{Hu8}
{\mathcal H}_j\,{\mathcal U}
+\sum_k({\mathcal L}_j^{-1}\,{\mathcal M}_j\,{\mathcal Z}_j
\,{\mathcal Z}_k)({\mathcal H}_k\,{\mathcal U})
\in V_p^{m,a}(\Omega)\cap V_2^{m,0}(\Omega).
\end{equation}
\noindent Taking into account that the norms of the matrix operator 
${\mathcal L}_j\,{\mathcal M}_j\,{\mathcal Z}_j\,{\mathcal Z}_k$ 
in the spaces $V_p^{m,a}(\Omega)$ and $V_2^{m,0}(\Omega)$ are $O(\theta)$, 
we may take $\theta>0$ small enough and obtain 
${\mathcal H}_j\,{\mathcal U}\in V_2^{m,0}(\Omega)$, i.e., 
${\mathcal U}\in V_2^{m,0}(\Omega)$. Therefore, 
${\mathcal L}:V_p^{m,a}(\Omega)\to V_p^{-m,a}(\Omega)$ is injective.

\vskip 0.08in
(iii) {\it A priori estimate}. Let $p\geq 2$ and assume that 
${\mathcal U}\in V_p^{m,a}(\Omega)$. Referring to Corollary~\ref{cor3} and 
arguing as in part (ii) of the present proof, we arrive at the equation
\begin{equation}\label{m32}
{\mathcal H}_j\,{\mathcal U}
+\sum_k({\mathcal L}_j^{-1}\,{\mathcal M}_j\,{\mathcal Z}_j\,
{\mathcal Z}_k)({\mathcal H}_k\,{\mathcal U})={\mathcal F},
\end{equation}
\noindent whose right-hand side satisfies
\begin{equation}\label{FVO}
\|{\mathcal F}\|_{V_p^{m,a}(\Omega)}\leq c
(\|{\mathcal A}\,{\mathcal U}\|_{V_p^{-m,a}(\Omega)}
+\|{\mathcal U}\|_{W_2^{m-1}(\omega)}),
\end{equation}
\noindent for some domain $\omega$ with ${\overline\omega}\subset\Omega$. 
Since the $V_p^{m,a}(\Omega)$-norm of the sum in (\ref{m32}) does not 
exceed $ C\theta\|{\mathcal U}\|_{V_p^{m,a}(\Omega)}$, we obtain the estimate 
\begin{equation}\label{m32-bis}
\|{\mathcal U}\|_{V_p^{m,a}(\Omega)}
\leq c\,(\|{\mathcal A}\,{\mathcal U}\|_{V_p^{-m,a}(\Omega)}
+\|{\mathcal U} \|_{W_2^{m-1}(\omega)}).
\end{equation}

\vskip 0.08in
(iv) {\it End of proof.} Let $p\geq 2$. The range of the operator 
${\mathcal A}:V_p^{m,a}(\Omega)\to V_p^{-m,a}(\Omega)$ is 
closed by (\ref{m32}) and the compactness of the restriction operator: 
$V_p^{m,a}(\Omega) \to W_2^{m-1}(\omega)$. Since the coefficients of the 
adjoint operator ${\mathcal L}^*$ satisfy the same conditions as those of 
${\mathcal L}$, we may conclude that
${\mathcal L}^*: V_{p'}^{m,a}(\Omega)\to V_{p'}^{-m,-a}(\Omega)$ 
is injective. Therefore, 
${\mathcal L}:V_{p}^{m,a}(\Omega)\to V_{p}^{-m,-a}(\Omega)$ is surjective. 
Being also injective, ${\mathcal L}$ is isomorphic if $p\geq 2$. 
Thus ${\mathcal L}^*$ is isomorphic for $p'\leq 2$ so 
${\mathcal L}$ is isomorphic for $p\geq 2$. 
\hfill$\Box$

\section{Traces and extensions}
\setcounter{equation}{0}

\subsection{Higher order Besov spaces on Lipschitz surfaces} 

Let $\Omega\subset\RR^n$ be a bounded Lipschitz domain and, for $m\in\NN$, 
$1<p<\infty$ and $-1/p<a<1-1/p$, consider a new space, 
$W_p^{m,a}(\Omega)$, consisting of functions ${\mathcal U}\in L_p(\Omega,loc)$ 
with the property that $\rho^{a}D^\alpha{\mathcal U}\in L_p(\Omega)$ for all 
multi-indices $\alpha$ with $|\alpha|=m$. We equip $W_p^{m,a}(\Omega)$ 
with the norm 
\begin{equation}\label{newW}
\|{\mathcal U}\|_{W_p^{m,a}(\Omega)}
:=\sum_{|\alpha|=m}\|D^\alpha{\mathcal U}\|_{L_p(\Omega,\,\rho(X)^{ap}\,dX)} 
+\|{\mathcal U}\|_{L_p(\omega)},
\end{equation}
\noindent where $\omega$ is an open non-empty domain, 
$\overline{\omega}\subset\Omega$. An equivalent norm is given by
the expression in (\ref{W-Nr}). We omit the standard proof of the fact that 
\begin{equation}\label{dense}
C^\infty({\overline\Omega})\hookrightarrow W_p^{m,a}(\Omega)
\quad\mbox{densely}. 
\end{equation}

Recall that for $p\in(1,\infty)$ and $s\in(0,1)$ the Besov space 
$B_p^s(\partial\Omega)$ is then defined via the requirement (\ref{Bes-xxx}). 
The nature of our problem requires that we work with Besov spaces
(defined on Lipschitz boundaries) which exhibit a higher order
of smoothness. In accordance with {\bf\cite{JW}}, {\bf{\cite{St}}}, 
{\bf{\cite{Wh}}}, we now make the following definition. 
\begin{definition}\label{def1} 
For $p\in(1,\infty)$, $m\in\NN$ and $s\in(0,1)$, define the (higher order) 
Besov space $\dot{B}^{m-1+s}_p(\partial\Omega)$ as the collection of all 
families $\dot{f}=\{f_\alpha\}_{|\alpha|\leq m-1}$ of 
measurable functions defined on $\partial\Omega$, such that if 
\begin{equation}\label{reminder}
R_\alpha(X,Y):=f_\alpha(X)-\sum_{|\beta|\leq m-1-|\alpha|}\frac{1}{\beta!}\,
f_{\alpha+\beta}(Y)\,(X-Y)^\beta,\qquad X,Y\in\partial\Omega,
\end{equation}
\noindent for each multi-index $\alpha$ of length $\leq m-1$, then 
\begin{eqnarray}\label{Bes-Nr}
\|\dot{f}\|_{\dot{B}^{m-1+s}_p(\partial\Omega)}
& := & \sum_{|\alpha|\leq m-1}\|f_\alpha\|_{L_p(\partial\Omega)}
\\[6pt]
&& +\sum_{|\alpha|\leq m-1}\Bigl(\int_{\partial\Omega}\int_{\partial\Omega}
\frac{|R_\alpha(X,Y)|^p}{|X-Y|^{p(m-1+s-|\alpha|)+n-1}}\,d\sigma_Xd\sigma_Y
\Bigr)^{1/p}<\infty.
\nonumber
\end{eqnarray}
\end{definition}
\noindent It is standard to prove that $\dot{B}^{m-1+s}_p(\partial\Omega)$ is a 
Banach space. Also, trivially, for any constant $\kappa>0$, 
\begin{equation}\label{Bes-Nr-eq}
\sum_{|\alpha|\leq m-1}\|f_\alpha\|_{L_p(\partial\Omega)}
+\sum_{|\alpha|\leq m-1}\,\,\,\Bigl(\,\,\,\,\,\,\,\,\,
\int\!\!\!\!\!\!\!\!\!\!\!\!\!\!\!\!\!\!\!
\int\limits_{{X,Y\in\partial\Omega}\atop{|X-Y|<\kappa}}
\frac{|R_\alpha(X,Y)|^p}{|X-Y|^{p(m-1+s-|\alpha|)+n-1}}\,d\sigma_Xd\sigma_Y
\Bigr)^{1/p}
\end{equation}
\noindent is an equivalent norm on $\dot{B}^{m-1+s}_p(\partial\Omega)$. 

A few notational conventions which will occasionally simplify the presentation
are as follows. Given a family of functions 
$\dot{f}=\{f_\alpha\}_{|\alpha|\leq m-1}$
on $\partial\Omega$ and $X\in\RR^n$, $Y,Z\in\partial\Omega$, set 
\begin{eqnarray}\label{PPP}
P_\alpha(X,Y) & := & \sum_{|\beta|\leq m-1-|\alpha|}\frac{1}{\beta!}\,
f_{\alpha+\beta}(Y)\,(X-Y)^\beta,\qquad \forall\,\alpha\,:\,|\alpha|\leq m-1,
\\[6pt]
P(X,Y) & := & P_{(0,...,0)}(X,Y).
\label{PPP-xk}
\end{eqnarray}
\noindent Then
\begin{equation}\label{PR-0}
R_\alpha(Y,Z)=f_\alpha(Y)-P_\alpha(Y,Z),
\qquad\forall\,\alpha\,:\,|\alpha|\leq m-1,
\end{equation}
\noindent and the following elementary identities hold for each 
multi-index $\alpha$ of length $\leq m-1$:
\begin{eqnarray}\label{PR-1}
i^{|\beta|}D^{\beta}_XP_{\alpha}(X,Y) & = & P_{\alpha+\beta}(X,Y),
\qquad|\beta|\leq m-1-|\alpha|,
\\[6pt]
P_{\alpha}(X,Y) -P_{\alpha}(X,Z) & =& 
\sum_{|\beta|\leq m-1-|\alpha|}\frac{1}{\beta!}R_{\alpha+\beta}(Y,Z)
(X-Y)^\beta.
\label{PR}
\end{eqnarray}
\noindent See, e.g., p.\,177 in {\bf\cite{St}} for the last formula. 

We now discuss how $\dot{B}^{m-1+s}_p(\partial\Omega)$ behaves 
under multiplication by a smooth function with compact support.
\begin{lemma}\label{mod-C}
For each $p\in(1,\infty)$, $m\in\NN$ and $s\in(0,1)$, the Besov space 
$\dot{B}^{m-1+s}_p(\partial\Omega)$ is a module over $C^\infty_0(\RR^n)$, 
granted that for each 
$\dot{f}=\{f_\alpha\}_{|\alpha|\leq m-1}\in\dot{B}^{m-1+s}_p(\partial\Omega)$
and $\psi\in C^\infty_0(\RR^n)$ we set 
\begin{equation}\label{cutoff}
\psi\,\dot{f}:=\Bigl\{\sum_{\beta+\gamma=\alpha}
\frac{\alpha !}{\beta !\gamma !}i^{|\beta|}\,{\rm Tr}\,[D^\beta\psi]
f_\gamma\Bigr\}_{|\alpha|\leq m-1}.
\end{equation}
\end{lemma}
\noindent{\bf Proof.} Let 
$\dot{f}=\{f_\alpha\}_{|\alpha|\leq m-1}\in\dot{B}^{m-1+s}_p(\partial\Omega)$
and $\psi\in C^\infty_0(\RR^n)$ be arbitrary and set 
\begin{equation}\label{dotpsi}
\dot{\psi}:=\Bigl\{i^{|\alpha|}\,{\rm Tr}\,[D^\alpha\psi]
\Bigr\}_{|\alpha|\leq m-1}.
\end{equation}
\noindent Denote by $\tilde{R}_\alpha(X,Y)$ the remainder (\ref{reminder}) 
written for $\psi\,\dot{f}$ in place of 
$\dot{f}=\{f_\alpha\}_{|\alpha|\leq m-1}$. Also, let 
$\tilde{P}_\alpha(X,Y)$ and $\tilde{P}(X,Y)$ be the polynomials defined 
in (\ref{PPP}), (\ref{PPP-xk}) with the components of $\dot{f}$ replaced 
by those of $\dot{\psi}$. 

Next, fix some $\alpha\in\NN_0^n$ with 
$|\alpha|\leq m-1$ and, for each $X,Y\in\partial\Omega$ write  
\begin{eqnarray}\label{Rpsif}
&& \tilde{R}_\alpha(X,Y)=\sum_{\beta+\gamma
=\alpha}\frac{\alpha !}{\beta !\gamma !}i^{|\beta|}\,D^\beta\psi(X)f_\gamma(X)
\nonumber\\[6pt]
&&\qquad\qquad
-\sum_{|\delta|\leq m-1-|\alpha|}\frac{1}{\delta !}
\Bigl[\sum_{\sigma+\tau=\alpha+\delta}\frac{(\alpha+\delta)!}{\sigma !\tau !}
i^{|\sigma|}\,D^\sigma\psi(Y)f_\tau(Y)\Bigr](X-Y)^\delta.
\end{eqnarray}
\noindent The crux of the matter is establishing that for a fixed 
$\kappa>0$ there exists $C>0$ independent of $\dot{f}$ such that 
\begin{equation}\label{Ral33}
\int\!\!\!\!\!\!\!\!\!\!\!\!\!\!\!\!\!\!\!
\int\limits_{{X,Y\in\partial\Omega}\atop{|X-Y|<\kappa}}
\frac{|\tilde{R}_\alpha(X,Y)|^p}
{|X-Y|^{p(m-1+s-|\alpha|)+n-1}}\,d\sigma_Xd\sigma_Y
\leq C\|\dot{f}\|^p_{\dot{B}^{m-1+s}_p(\partial\Omega)}.
\end{equation}

To get started, we note that the first sum in (\ref{Rpsif}) can be 
further expanded as 
\begin{eqnarray}\label{Rpsif-2}
\sum_{\beta+\gamma=\alpha}\frac{\alpha !}{\beta !\gamma !}
i^{|\beta|}\,D^\beta\psi(X)f_\gamma(X) &=& 
\sum_{\beta+\gamma=\alpha}\frac{\alpha !}{\beta !\gamma !}
\Bigl[i^{|\beta|}\,D^\beta\psi(X)-\tilde{P}_\beta(X,Y)\Bigr]f_\gamma(X)
\nonumber\\[6pt]
&&+\sum_{\beta+\gamma=\alpha}\frac{\alpha !}{\beta !\gamma !}
\tilde{P}_\beta(X,Y)\Bigl[f_\gamma(X)-P_\gamma(X,Y)\Bigr]
\nonumber\\[6pt]
&&+\sum_{\beta+\gamma=\alpha}\frac{\alpha !}{\beta !\gamma !}
\tilde{P}_\beta(X,Y)P_\gamma(X,Y).
\end{eqnarray}
\noindent Now, if $|X-Y|\leq\kappa$, the first sum in the right-hand 
side of (\ref{Rpsif-2}) is pointwise dominated by 
$C|X-Y|^{m-|\alpha|}\sum_{|\gamma|\leq m-1}|f_\gamma(X)|$ hence,
when raised to the $p$-th power and then 
multiplied by $|X-Y|^{-p(m-1+s-|\alpha|)+n-1}$, it is dominated by 
\begin{equation}\label{R-123}
C\sum_{|\gamma|\leq m-1}\qquad
\int\!\!\!\!\!\!\!\!\!\!\!\!\!\!\!\!\!\!\!
\int\limits_{{X,Y\in\partial\Omega}\atop{|X-Y|<\kappa}}\frac{|f_\gamma(X)|^p}
{|X-Y|^{n-1-p(1-s)}}\,d\sigma_Xd\sigma_Y
\leq C\sum_{|\gamma|\leq m-1}\|f_\gamma\|^p_{L_p(\partial\Omega)},
\end{equation}
\noindent which suits our goal. Similarly, for $|X-Y|\leq\kappa$, 
the second sum in the right-hand side of (\ref{Rpsif-2}) is 
$\leq C\sum_{|\gamma|\leq |\alpha|}|R_{\gamma}(X,Y)|$ thus, as before,
its contribution in the context of estimating the left-hand side of 
(\ref{Ral33}) does not exceed 
\begin{equation}\label{R-1234}
C\sum_{|\gamma|\leq |\alpha|}\qquad
\int\!\!\!\!\!\!\!\!\!\!\!\!\!\!\!\!\!\!\!
\int\limits_{{X,Y\in\partial\Omega}\atop{|X-Y|<\kappa}}
\frac{|R_\gamma(X,Y)|^p}
{|X-Y|^{p(m-1+s-|\gamma|)+n-1}}\,d\sigma_Xd\sigma_Y
\leq C\|\dot{f}\|^p_{\dot{B}^{m-1+s}_p(\partial\Omega)}. 
\end{equation}
\noindent As for the last sum in the right-hand side of (\ref{Rpsif-2}),
we employ (\ref{PR-1}), (\ref{PPP-xk}), Leibniz's rule and the definitions 
of $P(X,Y)$, $\tilde{P}(X,Y)$, in order to successively transform this sum into
\begin{eqnarray}\label{lastSum}
&& \sum_{\beta+\gamma=\alpha}\frac{\alpha !}{\beta !\gamma !}
i^{|\beta|}D^{\beta}_X\tilde{P}(X,Y)i^{|\gamma|}D^{\gamma}P(X,Y)
=i^{|\alpha|}D^{\alpha}_X[\tilde{P}(X,Y)P(X,Y)]
\nonumber\\[6pt]
&&\qquad =i^{|\alpha|}D^{\alpha}_X\left(\sum_{|\sigma|,|\tau|\leq m-1}
\frac{1}{\sigma !\tau !}i^{|\sigma|}D^{\sigma}\psi(Y)f_{\tau}(Y)
(X-Y)^{\sigma+\tau}\right).
\end{eqnarray}
\noindent Note that $D^{\alpha}_X[(X-Y)^{\sigma+\tau}]=0$ 
unless it is possible to select $\delta\in\NN_0^n$ such that 
$\sigma+\tau=\alpha+\delta$. In the latter situation, we use
$i^{|\alpha|}D^{\alpha}_X[(X-Y)^{\alpha+\delta}]
=\frac{(\alpha+\delta)!}{\delta !}(X-Y)^\delta$ and re-write the last 
expression in (\ref{lastSum}) as 
\begin{equation}\label{S-Last}
\sum_{|\delta|\leq 2(m-1)-|\alpha|}\frac{1}{\delta !}\Bigl[
\sum_{\sigma+\tau=\alpha+\delta}\frac{(\alpha+\delta)!}{\sigma !\tau !}
i^{|\sigma|}D^{\sigma}\psi(Y)f_{\tau}(Y)\Bigr](X-Y)^\delta. 
\end{equation}

Now, the second sum in (\ref{Rpsif}) cancels the portion from (\ref{S-Last}) 
corresponding to the case when $|\delta|\leq m-1-|\alpha|$, and
the remaining terms in this sum are 
$\leq C|X-Y|^{m-|\alpha|}\sum_{|\gamma|\leq m-1}|f_\gamma(X)|$. 
Consequently, in the context of (\ref{Ral33}), their contribution 
is estimated as we did in (\ref{R-123}). 

This analysis establishes (\ref{Ral33}). Since, trivially,  
$\|(\psi\dot{f})_{\alpha}\|_{L_p(\partial\Omega)}\leq C
\sum_{|\gamma|\leq m-1}\|f_\gamma\|_{L_p(\partial\Omega)}$, the
proof of the lemma is finished. 
\hfill$\Box$
\vskip 0.08in

Typically, the previous lemma is used to localize functions 
$\dot{f}\in\dot{B}^{m-1+s}_p(\partial\Omega)$ in such a way that the
supports of their components are contained in suitably small 
open subsets of $\partial\Omega$, where the boundary can be described
as a graph of a real-valued Lipschitz function defined in $\RR^{n-1}$. 
Such an argument, involving a smooth partition of unity, is 
standard and will often be used tacitly hereafter. 

We next discuss a special case of the general trace result we have in mind. 
\begin{lemma}\label{trace-1}
For each $1<p<\infty$, $-1/p<a<1-1/p$ and $s:=1-a-1/p$, the trace operator
\begin{equation}\label{TR-1}
{\rm Tr}:W^{1,a}_p(\Omega)\longrightarrow B^s_p(\partial\Omega)
\end{equation}
\noindent is well-defined, linear, bounded, onto and has $V^{1,a}_p(\Omega)$ 
as its null-space. Furthermore, there exists a linear, continuous mapping
\begin{equation}\label{Extension}
{\mathcal E}:B^s_p(\partial\Omega)\longrightarrow W^{1,a}_p(\Omega),
\end{equation}
\noindent called extension operator, such that ${\rm Tr}\circ{\mathcal E}=I$
(i.e., a bounded, linear right-inverse of trace). 
\end{lemma}
\noindent{\bf Proof.} By a standard argument involving a smooth partition 
of unity it suffices to deal with the case when $\Omega$ is the domain 
lying above the graph of a Lipschitz function $\varphi:\RR^{n-1}\to\RR$. 
Composing with the bi-Lipschitz homeomorphism 
$\RR^n_+\ni(X',X_n)\mapsto(X',\varphi(X')+X_n)\in\Omega$ further reduces
matters to the case when $\Omega=\RR^n_+$, in which situation the claims 
in the lemma have been proved in (greater generality) by S.V.\,Uspenski\u{\i} 
in {\bf\cite{Usp}} (a paper preceded by the significant work of E.\,Gagliardo 
in {\bf\cite{Ga}} in the unweighted case). 
\hfill$\Box$
\vskip 0.08in

We need to establish an analogue of Lemma~\ref{trace-1} for 
higher smoothness spaces. While for $\Omega=\RR^n_+$ this has been done by 
S.V.\,Uspenski\u{\i} in {\bf\cite{Usp}}, the flattening argument used in 
Lemma~\ref{trace-1} is no longer effective in this context. 
Let us also mention here that a result similar in spirit, valid for any 
Lipschitz domain $\Omega$ but with $B^{m-1+s+1/p}(\Omega)$ in place of
$W^{m,a}_p(\Omega)$ has been proved by A.\,Jonsson and H.\,Wallin in {\bf\cite{JW}} 
(in fact, in this latter context, these authors have dealt with much more
general sets than Lipschitz domains). The result which serves 
our purposes is as follows. 
\begin{proposition}\label{trace-2}
For $1<p<\infty$, $-1/p<a<1-1/p$, $s:=1-a-1/p\in(0,1)$ and $m\in\NN$, 
define the {\rm higher} {\rm order} trace operator
\begin{equation}\label{TR-11}
{\rm tr}_{m-1}:W^{m,a}_p(\Omega)\longrightarrow
\dot{B}^{m-1+s}_p(\partial\Omega)
\end{equation}
\noindent by setting 
\begin{equation}\label{Tr-DDD}
{\rm tr}_{m-1}\,\,{\mathcal U}
:=\Bigl\{i^{|\alpha|}\,{\rm Tr}\,[D^\alpha\,{\mathcal U}]\Bigr\}
_{|\alpha|\leq m-1},
\end{equation}
\noindent where the traces in the right-hand side are taken in the sense
of Lemma~\ref{trace-1}. Then {\rm (\ref{TR-11})}-{\rm (\ref{Tr-DDD})} is 
a well-defined, linear, bounded operator, which is onto and has 
$V^{m,a}_p(\Omega)$ as its null-space. Moreover, it has a bounded, linear 
right-inverse, i.e., there exists a linear, continuous operator 
\begin{equation}\label{Ext-222}
{\mathcal E}:\dot{B}^{m-1+s}_p(\partial\Omega)
\longrightarrow W^{m,a}_p(\Omega)
\end{equation}
\noindent such that 
\begin{equation}\label{Ext-333}
\dot{f}=\{f_\alpha\}_{|\alpha|\leq m-1}\in\dot{B}^{m-1+s}_p(\partial\Omega)
\Longrightarrow i^{|\alpha|}\,{\rm Tr}\,[D^\alpha({\mathcal E}\,\dot{f})]
=f_\alpha,\quad\forall\,\alpha\,:\,|\alpha|\leq m-1. 
\end{equation}
\end{proposition}

In order to facilitate the exposition, we isolate a couple of preliminary
results prior to the proof of Proposition~\ref{trace-2}. The first is
analogous to a Taylor remainder formula proved by H.\,Whitney in 
{\bf{\cite{Wh2}}} using a different set of compatibility conditions than 
(\ref{B-CC}) below. 
\begin{lemma}\label{Lemma-R}
Assume that $\varphi:\RR^{n-1}\to\RR$ is a Lipschitz function 
and define $\Phi:\RR^{n-1}\to\partial\Omega\hookrightarrow\RR^n$ by setting 
$\Phi(X'):=(X',\varphi(X'))$ at each $X'\in\RR^{n-1}$. Next, consider the 
special Lipschitz domain $\Omega:=\{X=(X',X_n)\in\RR^n:\,X_n>\varphi(X')\}$
and, for some fixed $m\in\NN$, a system of sufficiently nice functions 
$\{f_\alpha\}_{|\alpha|\leq m-1}$ with the property that 
\begin{equation}\label{B-CC}
\frac{\partial}{\partial X_k}[f_\alpha(\Phi(X'))]
=\sum_{j=1}^{n}f_{\alpha+e_j}(\Phi(X'))\partial_k\Phi_j(X'),
\qquad 1\leq k\leq n-1,\quad |\alpha|\leq m-2,
\end{equation}
\noindent where $\{e_j\}_j$ is the canonical orthonormal basis in $\RR^n$. 
Then, with $R_\alpha(X,Y)$ defined as in {\rm (\ref{reminder})}, 
the following identity holds for each multi-index $\alpha$ 
of length $\leq m-2$ and a.e. $X',Y'\in\RR^{n-1}$: 
\begin{eqnarray}\label{RRR=id}
&& R_\alpha(\Phi(X'),\Phi(Y'))
=\sum_{j=1}^n\sum_{|\gamma|=m-2-|\alpha|}\frac{1}{\gamma !}\int_0^1\Bigl[
f_{\alpha+\gamma+e_j}(\Phi(Y'+t(X'-Y')))-f_{\alpha+\gamma+e_j}(\Phi(Y'))\Bigr]
\nonumber\\[6pt]
&&\qquad\qquad\qquad
\times(\Phi(X')-\Phi(Y'+t(X'-Y')))^{\gamma}
\nabla\Phi_{j}(Y'+t(X'-Y'))\cdot(X'-Y')\,dt.
\end{eqnarray}
\end{lemma}
\noindent{\bf Proof.} We shall prove that for any system of functions
$\{f_\alpha\}_{|\alpha|\leq m-1}$ which satisfies (\ref{B-CC}), any 
multi-index $\alpha\in\NN_0^n$ with $|\alpha|\leq m-2$ and any $l\in\NN$
with $l\leq m-1-|\alpha|$, there holds
\begin{eqnarray}\label{FF=id}
&& f_\alpha(\Phi(X'))-\sum_{|\beta|\leq l}\frac{1}{\beta!}
f_{\alpha+\beta}(\Phi(Y'))(\Phi(X')-\Phi(Y'))^\beta 
\nonumber\\[6pt]
&&\quad=
\sum_{j=1}^n\sum_{|\gamma|=l-1}\frac{1}{\gamma !}\int_0^1\Bigl[
f_{\alpha+\gamma+e_j}(\Phi(Y'+t(X'-Y')))-f_{\alpha+\gamma+e_j}(\Phi(Y'))\Bigr]
\nonumber\\[6pt]
&&
\qquad\times(\Phi(X')-\Phi(Y'+t(X'-Y')))^{\gamma}
\nabla\Phi_{j}(Y'+t(X'-Y'))\cdot(X'-Y')\,dt.
\end{eqnarray}
\noindent Clearly, (\ref{RRR=id}) follows from (\ref{reminder}) 
and (\ref{FF=id}) by taking $l:=m-1-|\alpha|$. 

In order to justify (\ref{FF=id}) we proceed by induction on $l$. 
Concretely, when $l=1$ we may write, based on (\ref{B-CC}) and 
the Fundamental Theorem of Calculus,  
\begin{eqnarray}\label{RRR=0}
&& f_{\alpha}(\Phi(X'))-f_{\alpha}(\Phi(Y'))
-\sum_{j=1}^n f_{\alpha+e_j}(\Phi(Y'))(\Phi_j(X')-\Phi_j(Y'))
\\[6pt]
&&\quad =\int_0^1\frac{d}{dt}\Bigl[f_\alpha(\Phi(Y'+t(X'-Y')))\Bigr]\,dt
-\sum_{j=1}^n f_{\alpha+e_j}(\Phi(Y'))
\int_0^1\frac{d}{dt}\Bigl[\Phi_j(Y'+t(X'-Y'))\Bigr]\,dt
\nonumber\\[6pt]
&&\quad =\sum_{j=1}^n\Bigl\{\int_0^1\Bigl[f_{\alpha+e_j}(\Phi(Y'+t(X'-Y')))
-f_{\alpha+e_j}(\Phi(Y'))\Bigr]
\nabla\Phi_j(Y'+t(X'-Y'))\cdot(X'-Y')\,dt\Bigr\},
\nonumber
\end{eqnarray}
\noindent as wanted. 

To prove the version of (\ref{FF=id}) when $l$ 
is replaced by $l+1$ we split the sum in the left-hand side of (\ref{FF=id}),
written for $l+1$ in place of $l$, according to whether $|\beta|\leq l$
or $|\beta|=l+1$ and denote the expressions created in this fashion 
by $S_1$ and $S_2$, respectively. Next, based on 
(\ref{B-CC}) and the Fundamental Theorem of Calculus, we write  
\begin{eqnarray}\label{FTC}
&& f_{\alpha+\gamma+e_j}(\Phi(Y'+t(X'-Y')))
-f_{\alpha+\gamma+e_j}(\Phi(Y'))
\\[6pt]
&&
=\sum_{k=1}^n\int_0^{t}f_{\alpha+\gamma+e_j+e_k}
(\Phi(Y'+\tau(X'-Y')))
\nabla\Phi_{k}(Y'+\tau(X'-Y'))\cdot(X'-Y')\,d\tau
\nonumber
\end{eqnarray}
\noindent and use the induction hypothesis to conclude that 
\begin{eqnarray}\label{FF=id-2}
&& S_1=\sum_{j,k=1}^n\sum_{|\gamma|=l-1}\frac{1}{\gamma !}\int_0^1
\int_0^{t}f_{\alpha+\gamma+e_j+e_k}(\Phi(Y'+\tau(X'-Y')))
(\Phi(X')-\Phi(Y'+t(X'-Y')))^{\gamma}
\nonumber\\[6pt]
&&
\qquad
\times\nabla\Phi_{j}(Y'+t(X'-Y'))\cdot(X'-Y')
\nabla\Phi_{k}(Y'+\tau(X'-Y'))\cdot(X'-Y')\,d\tau\,dt.
\end{eqnarray}
\noindent Thus, if we set 
\begin{equation}\label{F=Psi}
F_j(t):=\Phi_j(X')-\Phi_j(Y'+t(X'-Y')),\qquad 1\leq j\leq n,
\end{equation}
\noindent and use an elementary identity to the effect that for any 
$\RR^n$-valued function $F=(F_1,...,F_n)$, 
\begin{equation}\label{F=Der}
\frac{d}{dt}\Bigl[\frac{1}{\beta!}F(t)^\beta\Bigr]
=\sum_{\gamma+e_j=\beta}\frac{1}{\gamma!}F(t)^\gamma F_j'(t),\qquad
\forall\,\beta\in\NN_0^n,
\end{equation}
\noindent we may express $S_1$ in the form 
\begin{eqnarray}\label{another}
&& S_1=-\sum_{k=1}^n\sum_{|\beta|=l}\int_0^1\int_0^{t}
[f_{\alpha+\beta+e_k}(\Phi(Y'+\tau(X'-Y')))
-f_{\alpha+\beta+e_k}(\Phi(Y'))]
\nonumber\\[6pt]
&&
\qquad\qquad\qquad\qquad\quad
\times\frac{d}{dt}\Bigl[\frac{1}{\beta!}F(t)^\beta\Bigr]
\nabla\Phi_{k}(Y'+\tau(X'-Y'))\cdot(X'-Y')\,d\tau\,dt
\nonumber\\[6pt]
&&
\qquad
+\sum_{k=1}^n\sum_{|\beta|=l}
f_{\alpha+\beta+e_k}(\Phi(Y'))\int_0^1\int_0^{t}
\frac{d}{dt}\Bigl[\frac{1}{\beta!}F(t)^\beta\Bigr]F_k'(\tau)\,d\tau\,dt.
\end{eqnarray}
Note that after changing the order of integration and using the Fundamental
Theorem of Calculus, the first double sum above corresponds precisely to 
the expression in the right-hand side of (\ref{FF=id}) written for $l+1$
in place of $l$. By once again changing the order of integration 
and relying on (\ref{F=Der}), it becomes apparent that the second double sum 
in (\ref{another}) is $-S_2$. Thus, $S_1+S_2$ matches
the right-hand side of (\ref{FF=id}) with $l$ replaced by $l+1$, 
proving (\ref{FF=id}).  
\hfill$\Box$
\begin{corollary}\label{Cor-R}
Under the assumptions of Lemma~\ref{Lemma-R}, 
for each multi-index $\alpha$ of length $\leq m-2$ the following estimate holds
\begin{equation}\label{RRR=est}
\Bigl(\int_{\partial\Omega}\int_{\partial\Omega}
\frac{|R_\alpha(X,Y)|^p}{|X-Y|^{p(m-1+s-|\alpha|)+n-1}}\,d\sigma_Xd\sigma_Y
\Bigr)^{1/p}\leq C\sum_{|\gamma|=m-1}\|f_\gamma\|_{B^s_p(\partial\Omega)},
\end{equation}
\noindent where the constant $C$ depends only on $n$, $p$, 
$s$ and $\|\nabla\varphi\|_{L_\infty(\RR^{n-1})}$. 
\end{corollary}
\noindent{\bf Proof.} The identity (\ref{RRR=id}) gives 
\begin{eqnarray}\label{pw-est-R}
&& |R_\alpha(\Phi(X'),\Phi(Y'))|
\\[6pt]
&&\qquad\qquad
\leq C|X'-Y'|^{m-1-|\alpha|}
\sum_{|\gamma|=m-1}\int_0^1
|f_{\gamma}(\Phi(Y'+t(X'-Y')))-f_{\gamma}(\Phi(Y'))|\,dt
\nonumber
\end{eqnarray}
\noindent for each $X',Y'\in\RR^{n-1}$, where the constant $C$ 
depends only on $n$ and $\|\nabla\Phi\|_{L_\infty(\RR^{n-1})}$ which, in turn,
is controlled in terms of $\|\nabla\varphi\|_{L_\infty(\RR^{n-1})}$. 
If we now integrate the $p$-th power of both sides in (\ref{pw-est-R}) for 
$X',Y'\in\RR^{n-1}$, use Fubini's Theorem and make the change of 
variables $Z':=Y'+t(X'-Y')$, we obtain
\begin{eqnarray}\label{r-vs-om}
&& \int_{\partial\Omega}\int_{\partial\Omega}
\frac{|R_\alpha(X,Y)|^p}{|X-Y|^{p(m-1-|\alpha|+s)+n-1}}\,d\sigma_Xd\sigma_Y
\leq C\int_{\RR^{n-1}}\int_{\RR^{n-1}}\frac{|R_\alpha(\Phi(X'),\Phi(Y'))|^p}
{|X'-Y'|^{p(m-1-|\alpha|+s)+n-1}}\,dX'dY'
\nonumber\\[6pt]
&&\qquad 
\leq C\sum_{|\gamma|=m-1}\int_{\RR^{n-1}}\int_{\RR^{n-1}}\int_0^1
\frac{|f_{\gamma}(\Phi(Y'+t(X'-Y')))-f_{\gamma}(\Phi(Y'))|^p}
{|X'-Y'|^{ps+n-1}}\,dt\,dX'dY'
\nonumber\\[6pt]
&& \qquad
\leq C\sum_{|\gamma|=m-1}\int_0^1\int_{\RR^{n-1}}\int_{\RR^{n-1}}
t^{ps}\,\frac{|f_{\gamma}(\Phi(Z'))-f_{\gamma}(\Phi(Y'))|^p}
{|Z'-Y'|^{ps+n-1}}\,dZ'dY'dt
\nonumber\\[6pt]
&& \qquad 
\leq C\sum_{|\gamma|=m-1}\|f_\gamma\|^p_{B^s_p(\partial\Omega)},
\end{eqnarray}
\noindent since $|\nabla\Phi(X')|\sim 1$ and 
$|\Phi(X')-\Phi(Y')|\sim |X'-Y'|$, uniformly for $X',Y'\in\RR^{n-1}$. 
\hfill$\Box$
\vskip 0.08in
After this preamble, we are in a position to present the 

\vskip 0.08in
\noindent{\bf Proof of Proposition~\ref{trace-2}.} We divide the proof 
into a series of steps, starting with 
\vskip 0.08in
\noindent{\tt Step I}: {\it The well-definiteness of trace.}  
Thanks to (\ref{dense}), it suffices to study the action of the trace
operator of a function ${\mathcal U}\in C^\infty(\overline{\Omega})
\hookrightarrow W^{m,a}_p(\Omega)$. If we now set 
\begin{equation}\label{ff-aa}
f_\alpha:=i^{|\alpha|}\,{\rm Tr}\,[D^\alpha\,{\mathcal U}],\qquad
\forall\,\alpha\,:\,|\alpha|\leq m-1,
\end{equation}
\noindent it follows from Lemma~\ref{trace-1} that these trace functions 
are well-defined and, in fact, 
\begin{equation}\label{falpha}
\sum_{|\alpha|\leq m-1}\|f_\alpha\|_{B^s_p(\partial\Omega)}
\leq C\|\,{\mathcal U}\|_{W^{m,a}_p(\Omega)}.
\end{equation}
In order to prove that $\dot{f}:=\{f_\alpha\}_{|\alpha|\leq m-1}$ 
belongs to $\dot{B}^{m-1+s}_p(\partial\Omega)$, let $R_\alpha(X,Y)$ 
be as in (\ref{reminder}). Our goal is to show that for 
every multi-index $\alpha$ with $|\alpha|\leq m-1$, 
\begin{equation}\label{B-alpha}
\Bigl(\int_{\partial\Omega}\int_{\partial\Omega}
\frac{|R_\alpha(X,Y)|^p}{|X-Y|^{p(m-1-|\alpha|+s)+n-1}}\,d\sigma_Xd\sigma_Y
\Bigr)^{1/p}\leq C\|\,{\mathcal U}\|_{W^{m,a}_p(\Omega)}. 
\end{equation}
\noindent To this end, we first observe that if $|\alpha|=m-1$ then the 
expression in the left-hand side of (\ref{B-alpha}) is majorized by 
$C\Bigl(\int_{\partial\Omega}\int_{\partial\Omega}
|f_\alpha(X)-f_\alpha(Y)|^p|X-Y|^{-(ps+n-1)}\,d\sigma_Xd\sigma_Y
\Bigr)^{1/p}$ which, by (\ref{falpha}), is indeed 
$\leq C\|\,{\mathcal U}\|_{W^{m,a}_p(\Omega)}$. To treat the case 
when $|\alpha|\leq m-2$ we assume that $\Omega$ is locally
represented as $\{X:\,X_n>\varphi(X')\}$ for some Lipschitz function 
$\varphi:\RR^{n-1}\to\RR$ and, as before, set $\Phi(X'):=(X',\varphi(X'))$, 
$X'\in\RR^{n-1}$. Then (\ref{B-CC}) holds, thanks to (\ref{ff-aa}),
for every multi-index $\alpha$ of length $\leq m-2$. 
Consequently, Corollary~\ref{Cor-R} applies and, in concert with 
(\ref{falpha}), yields (\ref{B-alpha}). This proves that the operator 
(\ref{TR-11})-(\ref{Tr-DDD}) is well-defined and bounded. 

\vskip 0.08in
\noindent{\tt Step II}: {\it The extension operator.}  
We introduce a co-boundary operator 
${\mathcal E}$ which acts on $\dot{f}=\{f_\alpha\}_{|\alpha|\leq m-1}\in
\dot{B}^{m-1+s}_p(\partial\Omega)$ according to 
\begin{equation}\label{def-Ee}
({\mathcal E}\dot{f})(X)
=\int_{\partial\Omega}{\mathcal K}(X,Y)\,P(X,Y)\,d\sigma_Y,
\qquad X\in\Omega,
\end{equation}
\noindent where $P(X,Y)$ is the polynomial associated with $\dot{f}$ as 
in (\ref{PPP-xk}). The integral kernel ${\mathcal K}$ is assumed to satisfy
\begin{eqnarray}\label{ker-prp}
&& \int_{\partial\Omega}{\mathcal K}(X,Y)\,d\sigma_Y=1
\qquad\mbox{for all }\,X\in\Omega,
\\[6pt]
&&
|D^\alpha_X{\mathcal K}(X,Y)|\leq c_\alpha\,\rho(X)^{1-n-|\alpha|},
\quad\forall\,X\in\Omega,\,\,\forall\,Y\in\partial\Omega,
\label{more-Kp}
\end{eqnarray}
\noindent where $\alpha$ is an arbitrary multi-index, and 
\begin{equation}\label{last-Kp}
{\mathcal K}(X,Y)=0\quad\mbox{if }\,\,|X-Y|\geq 2\rho(X).
\end{equation}
\noindent One can take, for instance, the kernel
\begin{equation}\label{K-def}
{\mathcal K}(X,Y):=\eta\left(\frac{X-Y}{\varkappa\rho_{\rm reg}(X)}\right)
\left(\int_{\partial\Omega}\eta\left(\frac{X-Z}{\varkappa\rho_{\rm reg}(X)} 
\right)d\sigma_Z\right)^{-1},
\end{equation}
\noindent where $\eta\in C^\infty_0(B_2)$, $\eta=1$ on $B_1$, $\eta\geq 0$ and 
$\varkappa$ is a positive constant depending on the Lipschitz constant 
of $\partial\Omega$. Here, as before, $\rho_{\rm reg}(X)$ stands for 
the regularized distance from $X$ to $\partial\Omega$.  

For each $X\in\Omega$ and $Z\in\partial\Omega$ and for every multi-index 
$\gamma$ with $|\gamma|=m$ we then obtain
\begin{equation}\label{N0}
D^\gamma{\mathcal E}\dot{f}(X)
=\sum_{{\alpha+\beta=\gamma}\atop{|\alpha|\geq 1}} 
\frac{\gamma!}{\alpha!\beta!}\int_{\partial\Omega} 
D^\alpha_X{\mathcal K}(X,Y)\,(P_\beta(X,Y)-P_\beta (X,Z))\,d\sigma_Y.
\end{equation}
\noindent Fix $\mu>1$ and denote by $B(X,R)$ the ball of radius $R$ 
centered at $X$. We may then estimate 
\begin{eqnarray}\label{N1}
|D^\gamma{\mathcal E}\dot{f}(X)|^p &\leq & C
\sum_{{\alpha+\beta=\gamma}\atop{|\alpha|\geq 1}} 
\rho(X)^{-p|\alpha|}\meanint_{B(X,\mu\rho(X))\cap\partial\Omega} 
|P_\beta(X,Y)-P_\beta(X,Z)|^p\,d\sigma_Y
\nonumber\\[6pt]
& \leq & C\sum_{{\alpha+\beta=\gamma}\atop{|\alpha|\geq 1}}  
\sum_{|\beta|+|\delta|\leq m-1}\rho(X)^{-p|\alpha|}
\meanint_{B(X,\mu\rho(X))\cap\partial\Omega} 
|R_{\delta+\beta}(Y,Z)|^p\,|X-Y|^{p|\delta|}\,d\sigma_Y,
\nonumber\\[6pt]
& \leq & C\sum_{|\tau|\leq m-1}\rho(X)^{p(|\tau|-m)}
\meanint_{B(X,\mu\rho(X))\cap\partial\Omega}|R_{\tau}(Y,Z)|^p\,d\sigma_Y,
\end{eqnarray} 
\noindent where we have used H\"older's inequality and (\ref{PR}).
Averaging the extreme terms in (\ref{N1}) for $Z$ in 
$B(X,\mu\rho(X))\cap\partial\Omega$, we arrive at
\begin{equation}\label{N2}
|D^\gamma{\mathcal E}\dot{f}(X)|^p \leq C 
\sum_{|\tau|\leq m-1}\rho(X)^{p(|\tau|-m)-2(n-1)}\,\,
\int\!\!\!\!\!\!\!\!\!\!\!\!\!\!\!\!\!\!\!
\int\limits_{{Y,Z\in\partial\Omega}\atop{|X-Y|,|X-Z|<\mu\rho(X)}}
\!\!\!\!\!\!|R_\tau(Y,Z)|^p\,d\sigma_Yd\sigma_Z.
\end{equation}

For each multi-index $\gamma$ of length $m$ we may then estimate:
\begin{eqnarray}\label{bigstep}
&& \int_{\Omega}|D^\gamma{\mathcal E}\dot{f}(X)|^p\rho(X)^{p(1-s)-1}\,dX
\nonumber\\[6pt]
&& 
\qquad\leq C\sum_{|\tau|\leq m-1}
\int_{\partial\Omega}\int_{\partial\Omega}|R_\tau(Y,Z)|^p
\Bigl(\int\limits_{{X\in\Omega}\atop{|X-Y|,|X-Z|<\mu\rho(X)}}
\rho(X)^{p(-m+1-s+|\tau|)-2n+1}\,dX\Bigr)\,d\sigma_Yd\sigma_Z
\nonumber\\[6pt]
&& 
\qquad\leq C\sum_{|\tau|\leq m-1}
\int_{\partial\Omega}\int_{\partial\Omega}
\frac{|R_\tau(Y,Z)|^p}{|Y-Z|^{p(m-1-|\tau|+s)+n-1}}\,d\sigma_Yd\sigma_Z
\nonumber\\[6pt]
&& 
\qquad\leq C\|\dot{f}\|^p_{\dot{B}^{m-1+s}_p(\partial\Omega)},
\end{eqnarray}
\noindent by (\ref{Bes-Nr}), where we have used the readily verified 
fact that there exists $C>0$ such that 
\begin{equation}\label{New-Es}
\int\limits_{{X\in\Omega}\atop{|X-Y|,|X-Z|<\mu\rho(X)}}
\rho(X)^{p(-m+1-s+|\tau|)-2n+1}\,dX
\leq C|Y-Z|^{p(-m+1+|\tau|-s)-n+1},
\end{equation}
\noindent for any $Y,Z\in\partial\Omega$ and $\tau\in\NN_0^n$ with 
$|\tau|\leq m-1$. This proves that the operator (\ref{def-Ee}) 
is well-defined and bounded in the context of (\ref{Ext-222}). 

\vskip 0.08in
\noindent{\tt Step III}: {\it The right-invertibility property.} 
We shall now show that the operator (\ref{def-Ee}) is a right-inverse for 
the trace operator (\ref{TR-11}), i.e., whenever 
$\dot{f}=\{f_\gamma\}_{|\gamma|\leq m-1}\in
\dot{B}^{m-1+s}_p(\partial\Omega)$, there holds 
\begin{equation}\label{N3a}
f_\gamma=i^{|\gamma|}\,{\rm Tr}[D^\gamma{\mathcal E}\dot{f}]
\end{equation}
\noindent for every multi-index $\gamma$ of length $\leq m-1$. 
To this end, for $|\gamma|\leq m-1$ we write 
\begin{equation}\label{N4} 
D^\gamma{\mathcal E}\dot{f}(X)-{\mathcal E}_\gamma\dot{f}(X) 
=\sum_{{\alpha+\beta=\gamma}\atop{|\alpha|\geq 1}}
\frac{\gamma!}{\alpha!\beta!}\int_{\partial\Omega} 
D^\alpha_X {\mathcal K}(X,Y)(P_\beta(X,Y)-P_\beta(X,Z))\,d\sigma_Y,
\end{equation}
\noindent where
\begin{equation}\label{N3b}
{\mathcal E}_\gamma\dot{f}(X)
:=\int_{\partial\Omega}{\mathcal K}(X,Y)\, P_\gamma(X,Y)\,d\sigma_Y,
\qquad X\in\Omega. 
\end{equation}
\noindent Estimating the right-hand side in (\ref{N4}) in the same way 
as we did with the right-hand side of (\ref{N0}), we obtain using 
the boundedness of $\partial\Omega$ 
\begin{eqnarray}\label{Dgamma-E}
\int_{\partial\Omega}
|D^\gamma{\mathcal E}\dot{f}(X)-{\mathcal E}_\gamma\dot{f}(X)|^p 
\rho(X)^{-ps-1}\,dX
&\leq & C\sum_{|\tau|\leq m-1}
\int_{\partial\Omega}\int_{\partial\Omega}\frac{|R_\tau(Y,Z)|^p}
{|Y-Z|^{p(|\gamma|+s-|\tau|)+n-1}}\,d\sigma_Yd\sigma_Z
\nonumber\\[6pt]
& \leq & C\,\|\dot{f}\|^p_{\dot{B}^{m-1+s}_p(\partial\Omega)}.
\end{eqnarray}
\noindent In a similar fashion, we check that
\begin{eqnarray}\label{sim-fash}
&& \int_{\partial\Omega}|\nabla (D^\gamma{\mathcal E}\dot{f}(X)
-{\mathcal E}_\gamma\dot{f}(X))|^p\rho(X)^{p-ps-1}\,dX
\nonumber\\[6pt]
&&\qquad
\leq C\sum_{|\tau|\leq m-1}
\int_{\partial\Omega}\int_{\partial\Omega}\frac{|R_\tau(Y,Z)|^p}
{|Y-Z|^{p(|\gamma|+s-|\tau|)+n-1}}\,d\sigma_Yd\sigma_Z
\leq C\,\|\dot{f}\|^p_{\dot{B}^{m-1+s}_p(\partial\Omega)}.
\end{eqnarray}
\noindent The two last inequalities imply 
$D^\gamma{\mathcal E}\dot{f}-{\mathcal E}\dot{f}\in V^{1,a}_p(\Omega)$ 
and, therefore,
\begin{equation}\label{N4a}
{\rm Tr}\,(D^\gamma{\mathcal E}\dot{f}-{\mathcal E}_\gamma\dot{f})=0.
\end{equation}
Going further, let us set
\begin{equation}\label{EgX}
Eg(X):=\int_{\partial\Omega}{\mathcal K}(X,Y)\,g(Y)\,d\sigma_Y, 
\qquad X\in\Omega.
\end{equation}
\noindent A simpler version of the reasoning in Step~II yields 
that $E$ maps $B^s_p(\partial\Omega)$ boundedly into $W^{1,a}_p(\Omega)$.
Also, a standard argument based on the Poisson kernel-like behavior of 
${\mathcal K}(X,Y)$ shows that ${\rm Tr}\,Eg=g$ for each 
$g\in B^s_p(\partial\Omega)$. 
Based on (\ref{PPP})-(\ref{PPP-xk}) and (\ref{N3b}) we have
\begin{eqnarray}\label{E-E}
&& |{\mathcal E}_\gamma\dot{f}(X)-Ef_\gamma(X)|^p
+\rho(X)^p|\nabla({\mathcal E}_\gamma\dot{f}(X)-Ef_\gamma(X))|^p 
\nonumber\\[6pt]
&& \qquad\qquad
\leq C\sum_{{|\beta|\leq m-1-|\gamma|}\atop{|\beta|\geq 1}} 
\rho(X)^{p|\beta|}\meanint_{B(X,\mu\rho(X))\cap\partial\Omega} 
|f_{\gamma+\beta}(Y)|^p\,d\sigma_Y.
\end{eqnarray}
\noindent Consequently, for an arbitrary Whitney cube $Q\subset\Omega$
of side-length $l$, we have
\begin{eqnarray}\label{n1}
&& \int_{Q}|{\mathcal E}_\gamma\dot{f}(X)-Ef_\gamma(X)|^p\rho(X)^{-ps-1}\,dX
+\int_{Q}
|\nabla({\mathcal E}_\gamma\dot{f}(X)-Ef_\gamma(X))|^p\rho(X)^{p-ps-1}\,dX
\nonumber\\[6pt]
&& \qquad\qquad\qquad\qquad\qquad\qquad
\leq C\sum_{{|\beta|\leq m-1-|\gamma|}\atop{|\beta|\geq 1}}
l^{p(|\beta|-s)}\int_{\partial\Omega\cap \varkappa Q}
|f_{\gamma+\beta}(Y)|^p\,d\sigma_Y,
\end{eqnarray}
\noindent where $\varkappa\,Q$ is the concentric dilate of 
$Q$ by some fixed factor $\varkappa>1$. Summing over all cubes $Q$ of 
a Whitney decomposition of $\Omega$ we find
\begin{equation}\label{Sm-wb}
\|{\mathcal E}_\gamma\dot{f}-Ef_\gamma\|_{V_p^{1,a}(\Omega)} 
\leq C\sum_{|\alpha|\leq m-1}\|f_\alpha\|_{L_p(\partial\Omega)}
\end{equation}
\noindent which implies
\begin{equation}\label{N5}
{\rm Tr}\,({\mathcal E}_\gamma\dot{f}-Ef_\gamma) =0.
\end{equation}
\noindent Finally, combining (\ref{N5}), (\ref{N4a}), and 
${\rm Tr}\,Ef_\gamma=f_\gamma$, we arrive at (\ref{N3a}). 

\vskip 0.08in
\noindent{\tt Step IV}: {\it The kernel of the trace.} 
We now turn to the task of identifying the null-space of the trace operator 
(\ref{TR-11})-(\ref{Tr-DDD}). For each $k\in\NN_0$ we 
denote by ${\mathcal P}_k$ the collection of all vector-valued, complex 
coefficient polynomials of degree $\leq k$ (and agree that 
${\mathcal P}_k= 0$ whenever $k$ is a negative integer). 
The claim we make at this stage is that the null-space of the operator
\begin{equation}\label{TRW}
W^{m,a}_p(\Omega)\ni {\mathcal W}\mapsto 
\Bigl\{i^{|\gamma|}\,{\rm Tr}\,[D^\gamma{\mathcal W}]\Bigr\}_{|\gamma|=m-1}
\in B^s_p(\partial\Omega)
\end{equation}
\noindent is given by 
\begin{equation}\label{Null-tr}
{\mathcal P}_{m-2}+V^{m,a}_p(\Omega).
\end{equation}
\noindent The fact that the null-space of the trace operator
(\ref{TR-11})-(\ref{Tr-DDD}) is $V^{m,a}_p(\Omega)$ follows readily from this.

That (\ref{Null-tr}) is included in the null-space of the operator (\ref{TRW})
is obvious. The opposite inclusion amounts to showing that if 
${\mathcal W}\in W_p^{m,a}(\Omega)$ is such that 
${\rm Tr}\,[D^\gamma{\mathcal W}]=0$ for all $\gamma\in{\mathbb{N}}_0$ with 
$|\gamma|=m-1$, then there exists $P_{m-2}\in {\mathcal P}_{m-2}$ with the 
property that ${\mathcal W}-P_{m-2}\in V_p^{m,a}(\Omega)$.  
To this end, we note that the case $m=1$ is a consequence of (\ref{Hardy}) 
and consider next the case $m=2$, i.e., when 
\begin{equation}\label{WWT}
{\mathcal W}\in W_p^{2,a}(\Omega),\qquad{\rm Tr}\,[\nabla{\mathcal W}]=0
\quad\mbox{on}\,\,\partial\Omega.
\end{equation}
\noindent Assume that $\{{\mathcal W}_j\}_{j\geq 1}$ is a sequence of smooth 
(even polynomial) vector-valued functions in $\overline{\Omega}$, 
approximating ${\mathcal W}$ in $W_p^{2,a}(\Omega)$. In particular, 
\begin{equation}\label{w1}
{\rm Tr}\,[\nabla{\mathcal W}_j]\to 0\,\,\,\mbox{ in }\,\,\,L_p(\partial\Omega)
\,\,\,\mbox{ as }\,\,j\to\infty. 
\end{equation}
\noindent If in a neighborhood of a point on $\partial\Omega$ the domain
$\Omega$ is given by $\{X:\,X_n>\varphi(X')\}$ for some Lipschitz function 
$\varphi$, the following chain rule holds for the gradient of the function 
$w_j:B'\ni X'\mapsto{\mathcal W}_j(X',\varphi(X'))$, where $B'$ is a 
$(n-1)$-dimensional ball: 
\begin{equation}\label{w2}
\nabla w_j(X')=
\Bigl(\nabla_{Y'}{\mathcal W}_j(Y',\varphi(X'))\Bigr)\Bigl|_{Y'= X'} 
+\Bigl(\frac{\partial}{\partial Y_n}{\mathcal W}_j(X',Y_n)\Bigr)
\Bigl|_{Y_n=\varphi(X')}\,\nabla\varphi(X').
\end{equation}
\noindent Since the sequence $\{w_j\}_{j\geq 1}$ is bounded in $L_p(B')$ and 
$\nabla w_j\to 0$ in $L_p(B')$, it follows that there exists a subsequence 
$\{j_i\}_i$ such that $w_{j_i}\to const$ in $L_p(B')$ 
(see Theorem~1.1.12/2 in {\bf\cite{Maz1}}). 
Hence, ${\rm Tr}\,{\mathcal W}=P_0=const$ on $\partial\Omega$. 
In view of ${\rm Tr}\,[{\mathcal W}-P_0]=0$ and 
${\rm Tr}\,[\nabla{\mathcal W}]=0$, we may conclude that 
${\mathcal W}-P_0\in V_p^{2,a}(\Omega)$ by Hardy's inequality. 
The general case follows in an inductive fashion, by reasoning as before 
with $D^\alpha{\mathcal W}$ with $|\alpha|=m-2$ in place of ${\mathcal W}$.
\hfill$\Box$
\vskip 0.08in

We now present a short proof of (\ref{Bes-X}), based on 
Proposition~\ref{trace-2}. 
\begin{proposition}\label{B-EQ}
Assume that $1<p<\infty$, $s\in(0,1)$ and $m\in\NN$. Then 
\begin{equation}\label{Bes-XYZ}
\Bigl\{\{i^{|\alpha|}D^\alpha\,{\mathcal V}|_{\partial\Omega}\}
_{|\alpha|\leq m-1}:\,{\mathcal V}\in C^\infty_0(\RR^n)\Bigr\}
\,\,\mbox{ is dense in }\,\,\dot{B}^{m-1+s}_p(\partial\Omega)
\end{equation}
\noindent and 
\begin{equation}\label{eQ-11}
\|\dot{f}\|_{\dot{B}^{m-1+s}_p(\partial\Omega)}\sim\sum_{|\alpha|\leq m-1}
\|f_\alpha\|_{B^{s}_p(\partial\Omega)},
\end{equation}
\noindent uniformly for 
$\dot{f}=\{f_\alpha\}_{|\alpha|\leq m-1}\in\dot{B}^{m-1+s}_p(\partial\Omega)$.
As a consequence, {\rm (\ref{Bes-X})} holds. 
\end{proposition}
\noindent{\bf Proof.} That (\ref{Bes-XYZ}) holds is a direct consequence 
of (\ref{dense}) and Proposition~\ref{trace-2}. Next, this density result 
and (\ref{RRR=est}) yield the left-pointing inequality in (\ref{eQ-11}). 
As for the opposite inequality, let  
$\dot{f}=\{f_\alpha\}_{|\alpha|\leq m-1}\in\dot{B}^{m-1+s}_p(\partial\Omega)$
and, with $a:=1-s-1/p$, consider 
${\mathcal U}:={\mathcal E}(\dot{f})\in W^{m,a}_p(\Omega)$. 
Then Lemma~\ref{trace-1} implies that, for each multi-index 
$\alpha$ of length $\leq m-1$, the function 
$f_\alpha=i^{|\alpha|}\,{\rm Tr}\,[D^\alpha{\mathcal U}]$ belongs 
to $B^s_p(\partial\Omega)$, plus a naturally accompanying norm estimate. 
\hfill$\Box$
\vskip 0.08in
We include one more equivalent characterization of the space 
$\dot{B}^{m-1+s}_p(\partial\Omega)$, in the spirit of work in 
{\bf\cite{AP}}, {\bf\cite{PV}}, {\bf\cite{Ve2}}. 
To state it, recall that $\{e_j\}_j$ is the canonical 
orthonormal basis in $\RR^n$. 
\begin{proposition}\label{CC-Aray}
Assume that $1<p<\infty$, $s\in(0,1)$ and $m\in\NN$. Then 
\begin{equation}\label{eQ0}
\{f_\alpha\}_{|\alpha|\leq m-1}\in\dot{B}^{m-1+s}_p(\partial\Omega)
\Longleftrightarrow
\left\{
\begin{array}{l}
f_\alpha\in B^{s}_p(\partial\Omega),\quad \forall\,\alpha\,:\,|\alpha|\leq m-1
\\[10pt]
\qquad\qquad\qquad
\mbox{ and }
\\[10pt]
(\nu_j\partial_k-\nu_k\partial_j)f_{\alpha}=
\nu_jf_{\alpha+e_k}-\nu_kf_{\alpha+e_j}
\\[6pt]
\forall\,\alpha\,:\,|\alpha|\leq m-2,\quad\forall\,j,k\in\{1,...,n\}.
\end{array}
\right.
\end{equation}
\end{proposition}
\noindent{\bf Proof.} The left-to-right implication is a consequence 
of (\ref{eQ-11}) and of the fact that (\ref{ff-aa}) holds for 
some ${\mathcal U}\in W^{m,a}_p(\Omega)$ (cf. Proposition~\ref{trace-2}). 
As for the opposite implication, we multiply $\{f_{\alpha}\}_{\alpha}$ (as 
in (\ref{cutoff})), with a function $\psi\in C^\infty_0(\RR^n)$ which  
can be assumed to satisfy ${\rm supp}\,\psi\cap\partial\Omega
={\rm supp}\,\psi\cap\{(X',X_n):\,X_n>\varphi(X')\}$ for some Lipschitz
function $\varphi:\RR^{n-1}\to\RR$. Furthermore, in this latter case, 
the compatibility conditions in (\ref{eQ0}) become equivalent to (\ref{B-CC}). 
Thus, given $\dot{f}:=\{f_\alpha\}_{|\alpha|\leq m-1}$ whose components 
are as in the right-hand side of (\ref{eQ0}), we may proceed as in the proof 
of Corollary~\ref{Cor-R} and use the fact that 
$f_\alpha\in B^{s}_p(\partial\Omega)$ for each $\alpha$ to conclude that 
$\dot{f}\in\dot{B}^{m-1+s}_p(\partial\Omega)$.  
\hfill$\Box$

\subsection{The space of Dirichlet data and the main trace theorem}

In this subsection, we study the mapping properties of the assignment 
${\mathcal U}\mapsto\Bigl\{\partial^k{\mathcal U}/\partial\nu^k\Bigr\}
_{0\leq k\leq m-1}$, for ${\mathcal U}\in W^{m,a}_p(\Omega)$.
In order to facilitate the subsequent discussion, for each 
$\ell\in{\mathbb N}$ we consider polynomial functions
$P^{\alpha\beta}_{\gamma jk}$ such that 
\begin{equation}\label{nuD-tang}
\nu^\beta D^\alpha -\nu^\alpha D^\beta 
=\sum_{|\gamma|=\ell-1}\sum_{j,k=1}^n P^{\alpha\beta}_{\gamma jk}(\nu)
\frac{\partial}{\partial\tau_{jk}}D^\gamma,\qquad
\forall\,\alpha,\beta\in{\mathbb N}_0^n\,:\,|\alpha|=|\beta|=\ell, 
\end{equation}
\noindent where $\partial/\partial\tau_{jk}$ is the tangential derivative 
given by 
\begin{equation}\label{tang-tau}
\frac{\partial}{\partial\tau_{jk}}:=\nu_j\frac{\partial}{\partial x_k}
-\nu_k\frac{\partial}{\partial x_j},\qquad 1\leq j,k\leq n.
\end{equation}
\noindent We shall also need Sobolev spaces of order one on 
$\partial\Omega$. Concretely, let 
$\nabla_{\rm tan}:=(\sum_j\nu_j\partial/\partial\tau_{jk})_{1\leq k\leq n}$
stand for the tangential gradient on the surface $\partial\Omega$ and, 
for $1<p<\infty$, introduce the space 
\begin{equation}\label{Lp1}
L_p^1(\partial\Omega):=\{f:\,\|f\|_{L^1_p(\partial\Omega)}:=
\|f\|_{L_p(\partial\Omega)}+\|\nabla_{\rm tan}f\|_{L_p(\partial\Omega)}
<\infty\}. 
\end{equation}
\noindent Our main trace/extension result then reads as follows. 
\begin{theorem}\label{P1-az} 
Let $\Omega$ be a bounded Lipschitz domain in $\RR^n$ and 
for $1<p<\infty$, $-1/p<a<1-1/p$, $s:=1-a-1/p\in(0,1)$ and $m\in{\mathbb N}$, 
consider the mapping 
\begin{equation}\label{newtrace}
{\rm Tr}_{m-1}:W^{m,a}_p(\Omega)\ni{\mathcal U}\mapsto
\Bigl\{\frac{\partial^k{\mathcal U}}{\partial\nu^k}\Bigr\}
_{0\leq k\leq m-1}\in L_p(\partial\Omega).
\end{equation}
\noindent Then its null-space is precisely $V^{m,a}_p(\Omega)$, and
its image can be characterized as follows. 

Given $(g_0,g_1,...,g_{m-1})\in L_p(\partial\Omega)$, set 
\begin{equation}\label{newt-CCCC}
f_{(0,...,0)}:=g_0
\end{equation}
\noindent and, inductively, having defined $\{f_\gamma\}_{|\gamma|\leq\ell-1}$
for some $\ell\in\{1,...,m-1\}$, consider 
\begin{equation}\label{newtrace-CC}
f_\alpha:=\nu^\alpha g_\ell+i\sum_{|\beta|=\ell}\sum_{|\gamma|=\ell-1}
\sum_{j,k=1}^n\,\frac{\ell !}{\beta !}\, 
\nu^\beta\,P^{\alpha\beta}_{\gamma jk}(\nu)\,
\frac{\partial f_\gamma}{\partial\tau_{jk}},\qquad
\forall\,\alpha\in{\mathbb N}_0^n\,:\,|\alpha|=\ell,
\end{equation}
\noindent where $P^{\alpha\beta}_{\gamma jk}$ is any fixed family 
of polynomials satisfying {\rm (\ref{nuD-tang})}.
Then there exists ${\mathcal U}\in W^{m,a}_p(\Omega)$ 
such that ${\rm Tr}_{m-1}\,{\mathcal U}=\{g_k\}_{0\leq k\leq m-1}$
if and only if 
\begin{equation}\label{CC-new}
f_\alpha\in L^1_p(\partial\Omega)\,\,\,\mbox{ if }\,\,\,|\alpha|\leq m-2,
\,\,\,\,\mbox{ and }\,\,\,\,f_\alpha\in B^s_p(\partial\Omega)\,\,\,\mbox{ if }
\,\,\,|\alpha|=m-1. 
\end{equation}

Furthermore, for ${\rm Tr}_{m-1}$ as in {\rm (\ref{newtrace})}, set 
\begin{eqnarray}\label{Im-newTr}
&& \dot{W}^{m-1+s}_p(\partial\Omega):=\mbox{ the image of the operator 
${\rm Tr}_{m-1}$ in {\rm (\ref{newtrace})}},
\\[4pt] 
&& \|g\|_{\dot{W}^{m-1+s}_p(\partial\Omega)}:=
\sum_{|\alpha|\leq m-1}\|f_\alpha\|_{B^s_p(\partial\Omega)}
\label{Im-newTr-norm}
\end{eqnarray}
\noindent if the families $g:=\{g_k\}_{0\leq k\leq m-1}$ and 
$\{f_\alpha\}_{|\alpha|\leq m-1}$ are related to oneanother as 
in {\rm (\ref{newt-CCCC})}-{\rm (\ref{newtrace-CC})}. Then this space
is independent of the particular choice of polynomials 
$P^{\alpha\beta}_{\gamma jk}$ satisfying {\rm (\ref{nuD-tang})} and 
the operator 
\begin{equation}\label{newtrace-22}
{\rm Tr}_{m-1}:W^{m,a}_p(\Omega)\longrightarrow
\dot{W}^{m-1+s}_p(\partial\Omega)
\end{equation}
\noindent is well-defined, bounded, and has a right-inverse. 
That is, there exists a bounded, linear operator 
\begin{equation}\label{ImTr-vc}
{\rm Ext}:\dot{W}^{m-1+s}_p(\partial\Omega)\longrightarrow W^{m,a}_p(\Omega)
\end{equation}
\noindent such that ${\rm Tr}_{m-1}\circ{\rm Ext}=I$, the identity. 
\end{theorem}

Prior to the proof of Theorem~\ref{P1-az} we will establish a useful  
approximation result, extending work done in {\bf{\cite{AP}}} for the 
case $m=2$. 
\begin{lemma}\label{DENSE}
Let $\Omega$ be a bounded Lipschitz domain in $\RR^n$ with outward unit 
normal $\nu$, and fix $1<p<\infty$, $m\in\NN$. Also, assume that 
$\{f_\alpha\}_{|\alpha|\leq m-1}$ is a family of functions satisfying 
\begin{equation}\label{f-CC-X-op}
f_{\alpha}\in L_p^1(\partial\Omega),
\qquad\forall\,\alpha\in\NN_0^n\,:\,|\alpha|\leq m-1,
\end{equation}
\noindent and 
\begin{equation}\label{f-CC-op}
\nu_jf_{\alpha+e_k}-\nu_kf_{\alpha+e_j}
=\frac{\partial f_{\alpha}}{\partial\tau_{jk}},\qquad
\forall\,\alpha\,:\,|\alpha|\leq m-2,\quad\forall\,j,k\in\{1,...,n\}. 
\end{equation}
\noindent Then there exists a sequence of functions 
$F^{\varepsilon}\in C^\infty_0(\RR^n)$, $\varepsilon>0$, such that
\begin{equation}\label{approx-f-op}
i^{|\alpha|}\,{\rm Tr}\,[D^\alpha F^{\varepsilon}]\longrightarrow f_\alpha
\,\,\mbox{ in }\,\,L^1_p(\partial\Omega)\,\,\mbox{ as }\,\,\varepsilon\to 0,
\qquad\forall\,\alpha\in\NN_0^n\,:\,|\alpha|\leq m-1. 
\end{equation}
\end{lemma}
\noindent{\bf Proof.} Since both the hypotheses 
(\ref{f-CC-X-op})-(\ref{f-CC-op}) and the conclusion (\ref{approx-f-op}) are
stable under multiplication by a smooth function with
compact support as in (\ref{cutoff}), there is no 
loss of generality in assuming that $\Omega=\{X:\,X_n>\varphi(X')\}$ 
for some Lipschitz function $\varphi:\RR^{n-1}\to\RR$, and that 
the functions $f_\alpha$ have compact support. In this setting, 
\begin{equation}\label{f-YYY}
\|f_{\alpha}\|_{L_p^1(\partial\Omega)}\sim
\|f_{\alpha}(\cdot\,,\varphi(\cdot))\|_{L_p^1(\RR^{n-1})},\qquad
\forall\,\alpha\in\NN_0^n\,:\,|\alpha|\leq m-1,
\end{equation}
\noindent and the compatibility conditions (\ref{f-CC-op}) can be 
written in the form
\begin{equation}\label{f-XXX-op}
\begin{array}{l}
\frac{\partial}{\partial X_j}\Bigl[f_{\alpha}(X',\varphi(X'))\Bigr]
=f_{\alpha+e_j}(X',\varphi(X'))+\partial_j\varphi(X')
f_{\alpha+e_n}(X',\varphi(X'))
\\[16pt]
\quad\mbox{for a.e. }X'\in\RR^{n-1},\quad
\forall\,\alpha\in\NN_0^n\,:\,|\alpha|\leq m-2,\quad 1\leq j\leq n-1.
\end{array}
\end{equation}

Next, fix a nonnegative function $\eta\in C^\infty_0(\RR^{n-1})$
which integrates to one and, for each $\varepsilon>0$, set 
$\eta_\varepsilon(X):=\varepsilon^{-n}\eta(X/\varepsilon)$. Then, for each 
$|\alpha|\leq m-1$, $\varepsilon>0$, and $X=(X',X_n)\in\RR^n$, consider
\begin{equation}\label{FFF-XXX}
F^{\varepsilon}_{\alpha}(X):=\sum_{k=0}^{m-1-|\alpha|}
\frac{1}{k!}\Bigl[(X_n-\varphi(\cdot))^k
f_{\alpha+k\,e_n}(\cdot\,,\varphi(\cdot))
\Bigr]\ast\eta_{\varepsilon}(X').
\end{equation}
\noindent Based on (\ref{f-XXX-op}), a straightforward calculation gives that 
whenever $r:=m-1-|\alpha|\geq 1$ and $1\leq j\leq n-1$,
\begin{eqnarray}\label{FFF-XXX-2}
&& \partial_jF^{\varepsilon}_{\alpha}(X)=F^{\varepsilon}_{\alpha+e_j}(X)
+\frac{1}{(r-1)!}\Bigl((X_n-\varphi(\cdot))^{r-1}\partial_{j}\varphi(\cdot)
f_{\alpha+r\,e_n}(\cdot\,,\varphi(\cdot))\Bigr)
\ast\eta_\varepsilon(X')
\nonumber\\[6pt]
&& \qquad\qquad\quad
+\frac{1}{r!}\Bigl((X_n-\varphi(\cdot))^{r}
f_{\alpha+r\,e_n}(\cdot\,,\varphi(\cdot))\Bigr)
\ast\partial_{j}(\eta_\varepsilon)(X').
\end{eqnarray}
\noindent Thus, if $r:=m-1-|\alpha|\geq 1$, after moving the derivative 
off of $\eta_\varepsilon$ in the last term above we arrive at the recurrence 
formula 
\begin{equation}\label{FFF-XXX-3}
\partial_jF^{\varepsilon}_{\alpha}(X)=
\left\{
\begin{array}{l}
F^{\varepsilon}_{\alpha+e_j}(X)
+\frac{1}{r!}\Bigl((X_n-\varphi(\cdot))^{r}
\partial_{j}(f_{\alpha+r\,e_n}(\cdot\,,\varphi(\cdot)))\Bigr)
\ast\eta_\varepsilon(X')\,\,\mbox{ if }\,\,j<n,
\\[16pt]
F^{\varepsilon}_{\alpha+e_n}(X)\,\,\mbox{ if }\,\,j=n.
\end{array}
\right.
\end{equation}
Hence, if we now define $F^{\varepsilon}(X):=F^{\varepsilon}_{(0,...,0)}(X)$, 
an inductive argument based on (\ref{FFF-XXX-3}) shows that, for any 
multi-index $\alpha=(\alpha',\alpha_n)\in\NN_0^{n-1}\times\NN$ of length 
$|\alpha'|+\alpha_n\leq m-1$, the difference between 
$i^{|\alpha|}D^{\alpha}F^{\varepsilon}(X)$ and 
$F^{\varepsilon}_{\alpha}(X)$ can be expressed as a finite, constant 
coefficient linear combination of terms of the type 
\begin{equation}\label{FFF-XXX-4}
\varepsilon^{-|\gamma|}\Bigl((X_n-\varphi(\cdot))^{m-1-|\beta|-\alpha_n}
\partial_{j}(f_{\delta}(\cdot\,,\varphi(\cdot)))\Bigr)
\ast(\partial^{\gamma}\eta)_\varepsilon(X'),
\end{equation}
\noindent where
\begin{equation}\label{cFFs}
1\leq j\leq n-1,\quad 
\beta,\,\gamma\in\NN_0^{n-1}\,\,\mbox{ are such that }\,\,
e_j+\beta+\gamma=\alpha',
\,\,\mbox{ and }\,\,\delta\in\NN_0^{n},\quad|\delta|=m-1.  
\end{equation}
\noindent Consequently, (\ref{approx-f-op}) will follow easily once we 
establish that for every $\alpha\in\NN_0^n$ of length $\leq m-1$, 
\begin{equation}\label{FFF-XXX-5}
F^{\varepsilon}_{\alpha}(\cdot,\,\varphi(\cdot))
\to f_{\alpha}(\cdot,\,\varphi(\cdot))\,\,\mbox{ in }\,\,L_p^1(\RR^{n-1})
\,\,\mbox{ as }\,\,\varepsilon\to 0,
\end{equation}
\noindent and that, whenever the indices are as in (\ref{cFFs}) and 
$X_n=\varphi(X')$, the expression in (\ref{FFF-XXX-4}) converges to zero 
in $L_p^1(\RR^{n-1})$ as $\varepsilon\to 0$. 
As regards (\ref{FFF-XXX-5}), we begin by noting that 
\begin{equation}\label{FFF-XXX-6}
F^{\varepsilon}_{\alpha}(X',\varphi(X'))=\sum_{k=0}^{m-1-|\alpha|}
\frac{1}{k!}\int_{\RR^{n-1}}\Bigl(\varphi(X')-\varphi(Y')\Bigr)^k
f_{\alpha+k\,e_n}(Y',\varphi(Y'))\,\eta_{\varepsilon}(X'-Y')\,dY'.
\end{equation}
\noindent Thus, clearly, $F^{\varepsilon}_{\alpha}(\cdot\,,\varphi(\cdot))\to 
f_{\alpha}(\cdot\,,\varphi(\cdot))$ in $L_p(\RR^{n-1})$ as $\varepsilon\to 0$,
and (\ref{FFF-XXX-5}) is proved as soon as we show that 
$\nabla_{X'}[F^{\varepsilon}_{\alpha}(\cdot\,,\varphi(\cdot))]\to 
\nabla_{X'}[f_{\alpha}(\cdot\,,\varphi(\cdot))]$ in $L_p(\RR^{n-1})$ 
as $\varepsilon\to 0$. We remark that this is obviously true if 
$|\alpha|=m-1$ (cf. (\ref{FFF-XXX-6})), so we consider the case 
when $|\alpha|\leq m-2$. In this situation, we use (\ref{FFF-XXX-6}) 
to compute, for each $k\in\{1,...,n-1\}$, 
\begin{eqnarray}\label{FFF-XXX-7}
\frac{\partial}{\partial X_k}\Bigl[F^{\varepsilon}_{\alpha}(X',\varphi(X'))
\Bigr]
&=&\int_{\RR^{n-1}}\frac{\partial}{\partial Y_k}\Bigl[
f_{\alpha+\,e_n}(Y',\varphi(Y'))\Bigr]\,\eta_{\varepsilon}(X'-Y')\,dY'
\nonumber\\[6pt]
&& +\int_{\RR^{n-1}}(\partial_k\varphi(X')-\partial_k\varphi(Y'))
f_{\alpha+\,e_n}(Y',\varphi(Y'))\,\eta_{\varepsilon}(X'-Y')\,dY'
\nonumber\\[6pt]
&& +{\mathcal R}(X')
\end{eqnarray}
\noindent with 
$|{\mathcal R}(X')|\leq C
\varepsilon\sum_{|\gamma|\leq m-1}|f_{\gamma}(\cdot\,,\varphi(\cdot))|
\ast\eta_{\varepsilon}(X')$. In particular, $\|{\mathcal R}\|_{L_p(\RR^{n-1})}$
goes to zero as $\varepsilon\to 0$. Consequently,
$\frac{\partial}{\partial X_k}\Bigl[F^{\varepsilon}_{\alpha}(X',\varphi(X'))
\Bigr]\to 
\frac{\partial}{\partial X_k}\Bigl[f_{\alpha}(X',\varphi(X'))\Bigr]$
in $L_p(\RR^{n-1})$ as $\varepsilon\to 0$, proving (\ref{FFF-XXX-5}). 

Let us now consider the expression (\ref{FFF-XXX-4}) when  
$\alpha\in\NN_0^n$ with $|\alpha|\leq m-1$ is fixed and the conditions 
in (\ref{cFFs}) hold. A direct estimate shows that, when $X_n=\varphi(X')$, 
the $L_p$-norm of this quantity is $\leq C\varepsilon^{m-|\alpha|}\to 0$ 
as $\varepsilon\to 0$. To finish the proof, take $X_n=\varphi(X')$ 
in (\ref{FFF-XXX-4}) and, for an arbitrary $j\in\{1,...,n-1\}$, apply 
$\partial/\partial X_j$. Much as we just did, the $L_p$-norm of this 
quantity is $\leq C\varepsilon^{m-|\alpha|-1}$. Thus, provided that 
$|\alpha|\leq m-2$, this converges to zero as $\varepsilon\to 0$. 
The most delicate case is when $|\alpha|=m-1$. In this situation, we write 
out the terms obtained as a result of making $X_n=\varphi(X')$ in 
(\ref{FFF-XXX-4}) and then applying $\partial/\partial X_k$ for some fixed
$k\in\{1,...,n-1\}$. They are 
\begin{equation}\label{term-X2}
\varepsilon^{-r}\,\int_{\RR^{n-1}}\Bigl(\varphi(X')-\varphi(Y')\Bigr)^{r}
g(Y')(\partial^{\gamma+e_k}\eta)_{\varepsilon}(X'-Y')\,dY',
\end{equation}
\noindent and 
\begin{equation}\label{term-X1}
r\,\varepsilon^{1-r}\,
\int_{\RR^{n-1}}\Bigl(\varphi(X')-\varphi(Y')\Bigr)^{r-1}
\partial_k\varphi(X')g(Y')(\partial^{\gamma}\eta)_{\varepsilon}(X'-Y')\,dY',
\end{equation}
\noindent where we have set $r:=m-1-|\beta|-\alpha_n$ and 
$g:=\partial_{j}(f_{\delta}(\cdot\,,\varphi(\cdot)))$. 
Above, we have used the fact that $|\alpha|=m-1$ forces $|\gamma|=r-1$. 
Our goal is to prove that the $L_p$-norm of 
the sum between (\ref{term-X1}) and (\ref{term-X2}), viewed as functions in 
$X\in\RR^{n-1}$, converges to zero as $\varepsilon\to 0$. To this end, write 
$\varphi(X')-\varphi(Y')=\Delta(X',Y')\,|X'-Y'|+\nabla\varphi(X')\cdot(X'-Y')$,
where 
\begin{equation}\label{term-X3}
\Delta(X',Y'):=\frac{\varphi(X')-\varphi(Y')+\varphi(X')(Y'-X')}{|X'-Y'|},
\end{equation}
\noindent then expand 
\begin{eqnarray}\label{term-X4}
\Bigr(\varphi(X')-\varphi(Y')\Bigr)^r &=&\sum_{a+b=r}\frac{r!}{a!b!}\,
\Delta(X',Y')^a\,|X'-Y'|^a\,(\nabla\varphi(X')\cdot(X'-Y'))^b
\\[6pt]
&=&\sum_{a+b=r}\sum_{{|\sigma|=b}\atop{\sigma\in\NN_0^{n-1}}}
\frac{r!}{a!\sigma !}\,
\Delta(X',Y')^a\,|X'-Y'|^a\,(\nabla\varphi(X'))^\sigma\,(X'-Y')^\sigma.
\nonumber
\end{eqnarray}
\noindent Plugging this back into (\ref{term-X2}) finally yields 
\begin{equation}\label{term-X2b}
\sum_{a+b=r}\sum_{{|\sigma|=b}\atop{\sigma\in\NN_0^{n-1}}}
\frac{r!}{a!\sigma !}\,
(\nabla\varphi(X'))^\sigma\,\int_{\RR^{n-1}}\Delta(X',Y')^a
g(Y')(\Theta_{\gamma+e_k}^{\sigma,a})_{\varepsilon}(X'-Y')\,dY',
\end{equation}
\noindent where we have used the notation
\begin{equation}\label{th-XF}
\Theta_{\tau}^{\sigma,a}(X'):=(X')^\sigma\,|X'|^a\,(\partial^\tau\eta)(X'),
\qquad X'\in\RR^{n-1},\,\,\sigma,\,\tau\in\NN_0^{n-1},\,\,a\in\NN_0.
\end{equation}
\noindent Each integral above is pointwise dominated by 
$C(\|\nabla\varphi\|_{L_\infty}){\mathcal M}g(X')$ uniformly with respect 
to $\varepsilon>0$ (recall that ${\mathcal M}$ is the Hardy-Littlewood 
maximal operator), and converges to zero as $\varepsilon\to 0$ whenever 
$a>0$ and $X'$ is a differentiability point for the function $\varphi$. 
Thus, since $\varphi$ is almost everywhere differentiable, by a well-known theorem 
of H.\,Rademacher, and since ${\mathcal M}$ is bounded on $L_p$ if $1<p<\infty$,
Lebesgue's Dominated Convergence Theorem gives that all integrals in 
(\ref{term-X2b}) corresponding to $a>0$ converge to zero in $L_p(\RR^{n-1})$ 
as $\varepsilon\to 0$. 

On the other hand, in the context of (\ref{term-X2b}), $a=0$ forces 
$|\sigma|=r=|\gamma+e_k|$. Note that in general, if $a=0$ and 
$|\sigma|=|\tau|$, definition (\ref{th-XF}) and repeated integrations 
by parts yield 
\begin{equation}\label{VanMom}
\int_{\RR^{n-1}}\Theta_{\tau}^{\sigma,0}(X')\,dX'
=\int_{\RR^{n-1}}(X')^{\sigma}(\partial^{\tau}\eta)(X')\,dX'
=(-1)^{|\tau|}\sigma!\,\delta_{\sigma\tau},
\end{equation}
\noindent where $\delta_{\sigma\tau}$ is the Kronecker symbol. Consequently, 
as $\varepsilon\to 0$, the portion of (\ref{term-X2b}) corresponding to $a=0$ 
(and, hence, the {\it entire} expression in (\ref{term-X2b})) converges 
in $L_p(\RR^{n-1})$ to 
\begin{equation}\label{Last-11}
(-1)^r\,r!\,(\nabla\varphi)^{\gamma+e_k}g.
\end{equation}
The analysis of (\ref{term-X1}) closely parallels that of (\ref{term-X2}).
In fact, given the close analogy between (\ref{term-X1}) and (\ref{term-X2}),
in order to compute the limit of the former in $L_p$ as 
$\varepsilon\to 0$, we only need to make the following changes 
in (\ref{Last-11}): replace $\gamma+e_k$ by $\gamma$, $r$ by $r-1$ 
and then multiply the result by $r\,\partial_k\varphi$. The resulting 
expression is precisely the opposite of (\ref{Last-11}), and this
finishes the proof. 
\hfill$\Box$
\vskip 0.08in
After this preamble, we are now ready to present the
\vskip 0.08in
\noindent{\bf Proof of Theorem~\ref{P1-az}.} The fact that $V^{m,a}_p(\Omega)$
is the null-space of ${\rm Tr}_{m-1}$ follows from Proposition~\ref{trace-2} 
once we notice that (\ref{newtrace}) is the composition between (\ref{TR-11}) 
and 
\begin{equation}\label{ass}
\dot{B}^{m-1+s}_p(\partial\Omega)\ni\dot{f}=\{f_\alpha\}_{|\alpha|\leq m-1}
\mapsto 
\Bigl\{\sum_{|\alpha|=k}\frac{k!}{\alpha!}\,\nu^\alpha\,f_\alpha\Bigr\}
_{0\leq k\leq m-1}\in L_p(\partial\Omega),
\end{equation}
\noindent and that the assignment (\ref{ass}) is one-to-one. 
The latter claim can be justified with the help of the identity 
\begin{equation}\label{m-xxx}
D^\alpha=i^{-|\alpha|}\,\nu^\alpha\,
\frac{\partial^{|\alpha|}}{\partial\nu^{|\alpha|}}
+\sum_{|\beta|=|\alpha|-1}\sum_{j,k=1}^n p^{\alpha\beta}_{jk}(\nu)
\frac{\partial}{\partial\tau_{jk}}D^\beta,
\end{equation}
\noindent where $p^{\alpha\beta}_{jk}$ are polynomial functions. 
Indeed, let $\dot{f}\in \dot{B}^{m-1+s}_p(\partial\Omega)$ be 
mapped to zero by the assignment (\ref{ass}) and consider 
${\mathcal U}:={\mathcal E}(\dot{f})\in W^{m,a}_p(\Omega)$. Then 
$f_\alpha=i^{|\alpha|}\,{\rm Tr}\,\,[D^\alpha\,{\mathcal U}]$ 
on $\partial\Omega$ for each $\alpha$ with $|\alpha|\leq m-1$ and, 
granted the current hypotheses, $\partial^k{\mathcal U}/\partial\nu^k=0$ 
for $k=0,1,...,m-1$. Consequently, (\ref{m-xxx}) and induction on $|\alpha|$
yield that ${\rm Tr}\,\,[D^\alpha\,{\mathcal U}]=0$ on $\partial\Omega$ 
whenever $|\alpha|\leq m-1$. Thus, ultimately, $f_\alpha=0$ 
for each $\alpha$ with $|\alpha|\leq m-1$, as desired. 
In turn, the identity (\ref{m-xxx}) can be proved by writing 
\begin{eqnarray}\label{m-xxx2}
i^{|\alpha|}\,D^\alpha & = & \prod_{j=1}^n\Bigl(\frac{\partial}{\partial x_j}
\Bigr)^{\alpha_j}
=\prod_{j=1}^n\Bigl[\sum_{k=1}^n\xi_k\Bigl(\xi_k
\frac{\partial}{\partial x_j}-\xi_j\frac{\partial}{\partial x_k}\Bigr)
+\sum_{k=1}^n\xi_j\xi_k\frac{\partial}{\partial x_k}\Bigr]^{\alpha_j}
\Bigl|_{\xi=\nu}
\nonumber\\[6pt]
& = & \prod_{j=1}^n\Bigl[\sum_{l=0}^{\alpha_j}\frac{\alpha_j!}{l!(\alpha_j-l)!}
\Bigl(\sum_{k=1}^n\xi_k\Bigl(\xi_k
\frac{\partial}{\partial x_j}-\xi_j\frac{\partial}{\partial x_k}\Bigr)\Bigr)
^{\alpha_j-l}\nu_j^{l}\frac{\partial^l}{\partial\nu^l}\Bigr]\Bigl|_{\xi=\nu}
\\[6pt]
& = & \prod_{j=1}^n\Bigl[
\nu_j^{\alpha_j}\frac{\partial^{\alpha_j}}{\partial\nu^{\alpha_j}}+
\sum_{l=0}^{\alpha_j-1}\frac{\alpha_j!}{l!(\alpha_j-l)!}
\Bigl(\sum_{k=1}^n\xi_k\Bigl(\xi_k
\frac{\partial}{\partial x_j}-\xi_j\frac{\partial}{\partial x_k}\Bigr)\Bigr)
^{\alpha_j-l}\nu_j^{l}\frac{\partial^l}{\partial\nu^l}\Bigr]\Bigl|_{\xi=\nu}
\nonumber
\end{eqnarray}
\noindent and noticing that 
$\prod_{j=1}^n\nu_j^{\alpha_j}\partial^{\alpha_j}/\partial\nu^{\alpha_j}
=\nu^\alpha\partial^{|\alpha|}/\partial\nu^{|\alpha|}$, 
whereas $(\xi_k\partial/\partial x_j-\xi_j\partial/\partial x_k)|_{\xi=\nu}
=-\partial/\partial\tau_{jk}$. 
Parenthetically, let us point out here that the identity (\ref{m-xxx2}) 
readily proves the existence of some polynomial function 
$P^{\alpha\beta}_{\gamma jk}$ such that (\ref{nuD-tang}) holds. 

Turning to the characterization of the image of the operator (\ref{newtrace}),
assume that $g_k\in L_p(\partial\Omega)$, $0\leq k\leq m-1$, are such that
the functions $f_\alpha$ defined as in (\ref{newt-CCCC})-(\ref{newtrace-CC}) 
belong to $B^s_p(\partial\Omega)$. The claim that we make is that 
$\dot{f}:=\{f_\alpha\}_{|\alpha|\leq m-1}\in\dot{B}^{m-1+s}_p(\partial\Omega)$
and 
\begin{equation}\label{f->g}
g_k=\sum_{|\alpha|=k}\frac{k!}{\alpha!}\,\nu^\alpha\,f_\alpha,\qquad
0\leq k\leq m-1.
\end{equation}
\noindent Regarding the first part of the claim, by (\ref{CC-new}) and 
Proposition~\ref{CC-Aray} it suffices to show that (\ref{f-CC-op}) holds. 
We shall prove by induction on $\ell:=|\alpha|\in\{0,...,m-2\}$. 
Based on (\ref{newt-CCCC})-(\ref{newtrace-CC}), we compute 
\begin{equation}\label{fj}
f_{\alpha}=\nu^\alpha g_1+\sum_{|\beta|=1}\nu^\beta
\frac{\partial g_0}{\partial\tau_{\beta\alpha}},
\qquad \forall\,\alpha\,:\,|\alpha|=1,
\end{equation}
\noindent from which the version of (\ref{f-CC-op}) with $\alpha=(0,...,0)$ is 
immediate. To prove the induction step, assume that (\ref{f-CC-op}) holds 
whenever $|\alpha|\leq \ell-1$. By Lemma~\ref{DENSE}, there exists  
$F_\varepsilon\in C^\infty_0(\RR^n)$, $\varepsilon>0$, such that 
\begin{equation}\label{approx-f}
|\gamma|\leq\ell\Longrightarrow 
i^{|\gamma|}\,{\rm Tr}\,[D^\gamma F_\varepsilon]\to f_\gamma
\,\,\mbox{ in }\,\,L^1_p(\partial\Omega)\,\,\mbox{ as }\,\,\varepsilon\to 0.
\end{equation}
\noindent From (\ref{approx-f}) and (\ref{nuD-tang}) it follows that for 
each $\alpha\in{\mathbb N}_0^n$ with $|\alpha|=\ell+1$,  
\begin{equation}\label{nuD-tang-2}
f_\alpha:=\nu^\alpha g_\ell+i^{\ell+1}\lim_{\varepsilon\to 0}
\sum_{|\beta|=\ell+1}\frac{(\ell+1)!}{\beta !}\,\nu^\beta\,
\Bigl(\nu^\beta\,{\rm Tr}\,[D^\alpha F_\varepsilon]
-\nu^\alpha\,{\rm Tr}\,[D^\beta F_\varepsilon]\Bigr)
\quad\mbox{ in }\,\,L_p(\partial\Omega).
\end{equation}
\noindent Next, fix an arbitrary $\alpha\in{\mathbb N}_0^n$ with 
$|\alpha|=\ell$, choose $j,k\in\{1,...,n\}$, and consider the identity 
(\ref{nuD-tang-2}) written twice, with $\alpha+e_k$ and $\alpha+e_j$, 
respectively, in place of $\alpha$. If we multiply the first such identity 
by $\nu_j$, the second one by $\nu_k$ and then subtract them from one another, 
we arrive at 
\begin{equation}\label{nuD-tang-3}
\nu_j f_{\alpha+e_k}-\nu_k f_{\alpha+e_j}=i^{\ell}\lim_{\varepsilon\to 0}
\sum_{|\beta|=\ell+1}\frac{(\ell+1)!}{\beta !}\,\nu^{2\beta}\,
\frac{\partial}{\partial\tau_{jk}}{\rm Tr}\,[D^\alpha F_\varepsilon].
\end{equation}
\noindent By (\ref{approx-f}), the above limit is 
$i^{-\ell}\partial f_\alpha/\partial\tau_{jk}$ and this finishes the proof of
the induction step. Thus (\ref{f-CC-op}) holds and, as a result, 
$\dot{f}:=\{f_\alpha\}_{|\alpha|\leq m-1}\in\dot{B}^{m-1+s}_p(\partial\Omega)$,
as desired. As for (\ref{f->g}), if we set 
\begin{equation}\label{UUU-f}
{\mathcal U}:={\mathcal E}\dot{f}\in W^{m,a}_p(\Omega),
\end{equation}
\noindent it follows from (\ref{Ext-333}), (\ref{newtrace-CC}) and 
(\ref{nuD-tang}) that 
\begin{equation}\label{solve-ff}
f_{\alpha}=\nu^\alpha g_k+\sum_{|\beta|=k}\frac{k!}{\beta !}
\,\nu^{\beta}\,(\nu^\beta f_\alpha-\nu^\alpha f_\beta),
\qquad\,\forall\,\alpha\,:\,|\alpha|=k,
\end{equation}
\noindent from which we deduce that $\nu^\alpha g_k=
\nu^\alpha\sum_{|\beta|=k}\frac{k!}{\beta !}\nu^{\beta}f_\beta$
for each multi-index $\alpha$ of length $k$. Multiplying both sides of 
this equality by $\frac{k!}{\alpha !}\nu^\alpha$ and summing over all 
$\alpha\in{\mathbb N}_0^n$ with $|\alpha|=k$ finally yields (\ref{f->g}). 

Going further, from (\ref{UUU-f}), (\ref{f->g}) and (\ref{nuk}), we may 
conclude that $\{g_k\}_{0\leq k\leq m-1}={\rm Tr}_{m-1}\,{\mathcal U}$, 
which proves that the family $\{g_k\}_{0\leq k\leq m-1}$ belongs to the 
image of the mapping (\ref{newtrace}). 
Conversely, if $\{g_k\}_{0\leq k\leq m-1}={\rm Tr}_{m-1}\,{\mathcal U}$
for some function ${\mathcal U}\in W^{m,a}_p(\Omega)$, it follows from 
(\ref{nuD-tang}) and (\ref{newtrace-CC}) that 
$f_\alpha=i^{|\alpha|}{\rm Tr}[D^\alpha{\mathcal U}]$ for $|\alpha|\leq m-1$.
Consequently, $f_\alpha\in B^s_p(\partial\Omega)$ if $|\alpha|\leq m-1$ and
$f_\alpha\in L^1_p(\partial\Omega)$ for $|\alpha|\leq m-2$, 
thanks to (\ref{eQ0}). This finishes the proof of the fact that (\ref{CC-new})
characterizes the image of the operator (\ref{newtrace}). 
That the space (\ref{Im-newTr}) is independent of the choice of polynomials 
$P^{\alpha\beta}_{\gamma jk}$ satisfying {\rm (\ref{nuD-tang})} is
implicit in the above reasoning. 
Finally, the results in \S{7.1} imply that the operator (\ref{newtrace-22}) 
is bounded. Since as a byproduct of the above proof, the assignment 
\begin{equation}\label{ass-xx}
\dot{B}^{m-1+s}_p(\partial\Omega)\ni\dot{f}=\{f_\alpha\}_{|\alpha|\leq m-1}
\mapsto 
\Bigl\{\sum_{|\alpha|=k}\frac{k!}{\alpha!}\,\nu^\alpha\,f_\alpha\Bigr\}
_{0\leq k\leq m-1}\in \dot{W}^{m-1+s}_p(\partial\Omega)
\end{equation}
\noindent is an isomorphism, we may take ${\rm Ext}$ in (\ref{ImTr-vc})
to be the composition between the operator (\ref{Ext-222}) and the 
inverse of the mapping (\ref{ass-xx}). This finishes the proof of the 
theorem.
\hfill$\Box$

A specific implementation of the algorithm 
(\ref{newt-CCCC})-(\ref{newtrace-CC}) is discussed below. 
\begin{corollary}\label{C0-az} 
Assume that $\Omega$ be a bounded Lipschitz domain in $\RR^n$ and 
fix $1<p<\infty$ and $-1/p<a<1-1/p$, $s:=1-a-1/p\in(0,1)$, $m\in{\mathbb N}$. 
For a family $(g_0,g_1,...,g_{m-1})\in L_p(\partial\Omega)$ set 
$f_{(0,...,0)}:=g_0$ and, inductively, if $\{f_\gamma\}_{|\gamma|\leq\ell-1}$
have already been defined for some $\ell\in\{1,...,m-1\}$, set 
\begin{equation}\label{nt-VV}
f_\alpha:=\nu^\alpha g_\ell+\frac{\alpha!}{\ell !}
\sum_{{\mu+\delta+e_j=\alpha}\atop{|\theta|=|\delta|}}
\frac{|\delta|!}{\delta !}\, \frac{|\mu|!}{\mu !}\,\frac{|\theta|!}{\theta !}
\,\nu^{\delta+\theta}\,(\nabla_{\rm tan}f_{\mu+\theta})_j,
\qquad\forall\,\alpha\in{\mathbb N}_0^n\,:\,|\alpha|=\ell, 
\end{equation}
\noindent where $(\cdot)_j$ is the $j$-th component. 
Then $\dot{g}=(g_0,g_1,...,g_{m-1})$ belongs to 
$\dot{W}^{m-1+s}_p(\partial\Omega)$ 
if and only if $\dot{f}:=\{f_{\alpha}\}_{|\alpha|\leq m-1}$ belongs to 
$\dot{B}^{m-1+s}_p(\partial\Omega)$, in which case (\ref{f->g}) also holds. 
\end{corollary}
\noindent{\bf Proof.} For any two multi-indices $\alpha,\beta\in\NN_0^n$ of 
length $\ell$ written as $\alpha=e_{j_1}+\cdots e_{j_\ell}$ and 
$\beta=e_{k_1}+\cdots e_{k_\ell}$, a direct calculation yields 
\begin{equation}\label{indix-2}
\nu^\beta D^\alpha-\nu^\alpha D^\beta
=i^{\ell}\sum_{r=0}^{\ell-1}\nu_{k_1}\cdots\nu_{k_{\ell-r-1}}\nu_{j_{\ell-r+1}}
\cdots\nu_{j_{\ell}}\frac{\partial}{\partial\tau_{k_{\ell-r}j_{\ell-r}}}
\partial_{j_1}\cdots\partial_{j_{\ell-r-1}}\partial_{k_{\ell-r+1}}
\cdots\partial_{j_{\ell}}.
\end{equation}
\noindent In order to be able to re-write (\ref{indix-2}) in multi-index 
notation, it is convenient to symmetrize the right-hand side of this identity 
by adding up all its versions obtained by permuting the indices 
$j_1,...,j_\ell$ and $k_1,...,k_\ell$. In this fashion, 
we obtain 
\begin{equation}\label{indix-3}
\nu^\beta D^\alpha-\nu^\alpha D^\beta
=\frac{1}{i}\frac{\alpha !}{\ell !}\frac{\beta !}{\ell !}
\sum_{r=0}^{\ell-1}\sum_{{\mu+\delta+e_j=\alpha,\,|\delta|=r}
\atop{\gamma+\theta+e_k=\beta,\,|\theta|=r}}
\frac{(\ell-r-1)!}{\mu !}\frac{r!}{\delta !}
\frac{(\ell-r-1)!}{\gamma !}\frac{r!}{\theta !}
\nu^{\gamma+\delta}\frac{\partial}{\partial\tau_{kj}}D^{\mu+\theta}.
\end{equation}
\noindent This is a particular version of (\ref{nuD-tang}), where the 
intervening polynomials are identified explicitly. If we now implement the 
algorithm (\ref{newt-CCCC})-(\ref{newtrace-CC}), for each 
$\alpha\in{\mathbb N}_0^n$ with $|\alpha|=\ell$ we arrive at 
\begin{equation}\label{indix-4}
f_\alpha=\nu^\alpha g_\ell+\frac{\alpha!}{\ell !}\sum_{|\beta|=\ell}
\sum_{r=0}^{\ell-1}\sum_{{\mu+\delta+e_j=\alpha,\,|\delta|=r}
\atop{\gamma+\theta+e_k=\beta,\,|\theta|=r}}
\frac{(\ell-r-1)!}{\mu !}\frac{r!}{\delta !}
\frac{(\ell-r-1)!}{\gamma !}\frac{r!}{\theta !}
\nu^{\gamma+\delta}\frac{\partial f_{\mu+\theta}}{\partial\tau_{kj}}.
\end{equation}
\noindent Next, we replace $\beta$ by $\gamma+\theta+e_k$, eliminating the sum 
over $\beta$, and make use of the identities 
\begin{equation}\label{indix-5}
\sum_{|\gamma|=\ell-r-1}\frac{(\ell-r-1)!}{\gamma !}\nu^{2\gamma}=1,\qquad
\sum_{k=1}^n\nu_k\frac{\partial f}{\partial\tau_{kj}}=(\nabla_{\rm tan}\,f)_j,
\end{equation}
\noindent in order to transform (\ref{indix-4}) into 
\begin{equation}\label{indix-6}
f_\alpha=\nu^\alpha g_\ell+\frac{\alpha!}{\ell !}
\sum_{r=0}^{\ell-1}\sum_{{\mu+\delta+e_j=\alpha}
\atop{|\theta|=|\delta|=r}}
\frac{(\ell-r-1)!}{\mu !}\frac{r!}{\delta !}\frac{r!}{\theta !}
\nu^{\delta+\theta}(\nabla_{\rm tan}\,f_{\mu+\theta})_j,
\end{equation}
\noindent which is equivalent to (\ref{nt-VV}). 
\hfill$\Box$
\vskip 0.08in

The space (\ref{Im-newTr}) takes a particularly simple form when $m=2$. 
Indeed, as a direct consequence of (\ref{nt-VV}) in which we take $\ell=1$
we have: 
\begin{corollary}\label{WWW-2}
For each Lipschitz domain $\Omega\subset\RR^n$ and each $1<p<\infty$, 
$s\in(0,1)$, 
\begin{equation}\label{W2-az}
\dot{W}^{1+s}_p(\partial\Omega)
=\{(g_0,g_1)\in L^1_p(\partial\Omega)\oplus L_p(\partial\Omega):\,
\nu g_1+\nabla_{\rm tan}\,g_0\in B^s_p(\partial\Omega)\}.
\end{equation}
\end{corollary}
\noindent This has been conjectured to hold (when $s=1-1/p$) by 
A.\,Buffa and G.\,Geymonat on p.\,703 of {\bf{\cite{BG}}}.

Finally, we comment on how (\ref{Im-newTr}) relates to more classical spaces 
of higher order traces when $\Omega\subset\RR^n$ has a smoother boundary
than mere Lipschitz. Specifically, fix $m\in\NN$, $m\geq 2$, 
$p\in(1,\infty)$, $s\in(0,1)$ and assume that $\partial\Omega$ is locally 
given by graphs of Lipschitz function $\varphi:\RR^{n-1}\to\RR$ 
with the additional property that $\nabla\varphi$ belongs to 
$MB^{m-2+s}_p(\RR^{n-1})$, the space of (pointwise) multipliers for the 
Besov space $B^{m-2+s}_p(\RR^{n-1})$ (cf. {\bf{\cite{MS}}}, {\bf{\cite{MS2}}}).
Then, for each non-integer $1<\mu\leq m-1+s$, one can coherently define the 
space $B^{\mu}_p(\partial\Omega)$ by starting from $B^{\mu}_p(\RR^{n-1})$ and 
then transporting it to $\partial\Omega$ via a smooth partition of unity 
argument and by locally flattening the boundary. In fact, we arrive at 
the same space by taking the image of the trace operator on $\partial\Omega$, 
acting from $B^{\mu+1/p}_p(\RR^n)$. 
\begin{proposition}\label{WWW-1}
Assume that $p\in(1,\infty)$, $s\in (0,1)$
and that $\Omega\subset\RR^n$ is a Lipschitz domain whose boundary 
is locally described by means of graphs of real-valued functions in $\RR^{n-1}$
whose gradients belong to $MB^{m-2+s}_p(\RR^{n-1})$. Then 
\begin{equation}\label{image}
\dot{W}^{m-1+s}_p(\partial\Omega)
=\prod_{k=0}^{m-1}B^{m-1-k+s}_p(\partial\Omega).
\end{equation}
\noindent In particular, this is the case if 
$\partial\Omega\in C^{m-1,\theta}$ for some $\theta>s$. 
\end{proposition}
\noindent{\bf Proof.} In one direction, (\ref{eQ0}) and 
lifting theorems imply that if $\{f_\alpha\}_{|\alpha|\leq m-1}\in 
\dot{B}^{m-1+s}_p(\partial\Omega)$ then 
$f_\alpha\in B^{m-1-|\alpha|+s}_p(\partial\Omega)$ for each $\alpha$ with 
$|\alpha|\leq m-1$. Hence,  
$g_k:=\sum_{|\alpha|=k}\frac{k!}{\alpha!}\,\nu^\alpha\,f_\alpha
\in B^{m-1-k+s}_p(\partial\Omega)$ for each $k\in\{0,...,m-1\}$, 
so the left-to-right inclusion in (\ref{image}) follows from the fact that 
(\ref{ass-xx}) is an isomorphism. 

As for the opposite implication, given 
$\{g_k\}_{0\leq k\leq m-1}\in\oplus_{k=0}^{m-1}B^{m-1-k+s}_p(\partial\Omega)$,
define $\{f_\alpha\}_{|\alpha|\leq m-1}$ as in 
(\ref{newt-CCCC})-(\ref{newtrace-CC}). 
Granted the current assumptions on $\partial\Omega$, an argument based on 
induction and the fact that $\nu\in MB^{m-2+s}_p(\partial\Omega)$ 
shows that $f_\alpha\in B^{m-1-|\alpha|+s}_p(\partial\Omega)$
for each $|\alpha|\leq m-1$. In particular, (\ref{CC-new}) holds 
which proves that $\{g_k\}_{0\leq k\leq m-1}\in
\dot{W}^{m-1+s}_p(\partial\Omega)$. This shows that the right-to-left 
inclusion in (\ref{image}) also holds, thus finishing the proof of 
the proposition. 
\hfill$\Box$
\vskip 0.08in

\section{Proof of the main result}
\setcounter{equation}{0}

\subsection{The inhomogeneous Dirichlet problem}

Theorem~\ref{Theorem} is a particular case of Theorem~\ref{Theorem1}, 
concerning the solvability of the inhomogeneous Dirichlet problem 
\begin{equation}\label{m9}
\left\{
\begin{array}{l}
{\mathcal A}(X,D_X)\,{\mathcal U}={\mathcal F}
\qquad\mbox{in}\,\,\Omega,
\\[15pt] 
{\displaystyle{\frac{\partial^k{\mathcal U}}{\partial\nu^k}}}
=g_k\quad\,\,\mbox{on}\,\,\partial\Omega,\,\,\,\,\,0\leq k\leq m-1,
\end{array}
\right.
\end{equation}
\noindent in the space $W_p^{m,a}(\Omega)$. Note that for any operator 
${\mathcal A}$ as in \S{6.1} we have 
\begin{equation}\label{ADXW}
{\mathcal A}(X,D_X):W_p^{m,a}(\Omega)\longrightarrow V_p^{-m,a}(\Omega)
\end{equation}
\noindent boundedly. Thus, granted the membership of ${\mathcal U}$ solving 
(\ref{m9}) to $W^{m,a}_p(\Omega)$, it follows from Theorem~\ref{P1-az} 
that necessarily ${\mathcal F}\in V_p^{-m,a}(\Omega)$ and 
$g:=\{g_k\}_{0\leq k\leq m-1}\in \dot{W}^{m-1+s}_p(\partial\Omega)$. 
Moreover, 
\begin{equation}\label{estUU-vvv}
\|g\|_{\dot{W}^{m-1+s}_p(\partial\Omega)}
+\|{\mathcal F}\|_{V_p^{-m,a}(\Omega)}
\leq C\|{\mathcal U}\|_{W_p^{m,a}(\Omega)}.
\end{equation}
The converse direction makes the object of the theorem below. 
\begin{theorem}\label{Theorem1}
Let $\Omega$ be a bounded Lipschitz domain in $\RR^n$ and assume that the 
operator ${\mathcal A}$ is as in \S{6.1}. Then there exists $c>0$ such 
that if {\rm (\ref{a0})} is satisfied then the Dirichlet problem 
{\rm (\ref{m9})} has a unique solution ${\mathcal U}\in W_p^{m,a}(\Omega)$ 
for any given 
${\mathcal F}\in V_p^{-m,a}(\Omega)$ and $g:=\{g_k\}_{0\leq k\leq m-1}\in 
\dot{W}^{m-1+s}_p(\partial\Omega)$. Furthermore, there exists 
$C=C(\partial\Omega,{\mathcal A},p,s)>0$ 
such that 
\begin{equation}\label{estUU}
\|{\mathcal U}\|_{W_p^{m,a}(\Omega)}
\leq C\Big(\|g\|_{\dot{W}^{m-1+s}_p(\partial\Omega)}
+\|{\mathcal F}\|_{V_p^{-m,a}(\Omega)}\Bigr). 
\end{equation}
\end{theorem}
\noindent{\bf Proof.} We seek a solution for (\ref{m9}) in the form 
${\mathcal U}={\rm Ext}(g)+{\mathcal W}$, where ${\rm Ext}$ denotes the 
extension operator from Theorem~\ref{P1-az} and  
${\mathcal W}\in V^{m,a}_p(\Omega)$. Note that, by Theorem~\ref{P1-az}, 
this membership automatically entails $\partial^k {\mathcal W}/\partial\nu^k=0$
on $\partial\Omega$ for $k=0,1,...,m-1$, so it suffices to take 
\begin{equation}\label{m10}
{\mathcal W}:={\mathcal A}(X,D_X)^{-1}\Bigl(
{\mathcal F}-{\mathcal A}(X,D_X){\rm Ext}\,(g)\Bigr)\in 
V^{m,a}_p(\Omega)
\end{equation}
\noindent which, by (\ref{ADXW}) and Theorem~\ref{th1a}, is meaningful.  
As for uniqueness, let ${\mathcal U}\in W_p^{m,a}(\Omega)$ solve (\ref{m9}) 
with ${\mathcal F}=0$ and $g_k=0$, $0\leq k\leq m-1$. Then the function 
${\mathcal U}$ belongs to $V_p^{m,a}(\Omega)$, thanks to 
Theorem~\ref{P1-az}, and is a null-solution of ${\mathcal A}(X,D_X)$. 
In turn, Theorem~\ref{th1a} gives that ${\mathcal U}=0$, as desired. 
Finally, (\ref{estUU}) is a consequence of the results in \S{7}.
\hfill$\Box$

We conclude this subsection with a remark pertaining to the presence of lower 
order terms. More specifically, granted Theorem~\ref{Theorem1}, a standard 
perturbation argument (cf., e.g., {\bf{\cite{Ho}}}) proves the following. 
Assume that 
\begin{equation}\label{E444-bis}
{\mathcal A}(X,D_X)\,{\mathcal U}
:=\sum_{0\leq |\alpha|,|\beta|\leq m}D^\alpha({\mathcal A}_{\alpha\beta}(X)
\,D^\beta{\mathcal U}),\qquad X\in\Omega,
\end{equation}
\noindent where the top part of ${\mathcal A}(X,D_X)$ satisfies the 
hypotheses made in Theorem~\ref{Theorem} and the lower order terms are bounded.
Then, assuming that either (\ref{cond-coeffx2}) or (\ref{a0}) holds, the 
Dirichlet problem (\ref{m9}) is Fredholm with index zero, in the sense that 
the operator 
\begin{equation}\label{Fred-AA}
W^{m,a}_p(\Omega)\ni{\mathcal U}\mapsto
\Bigl({\mathcal A}(X,D_X)\,{\mathcal U}\,\,,\,\,
\{\partial^k{\mathcal U}/\partial\nu^k\}_{0\leq k\leq m-1}\Bigr)
\in V_p^{-m,a}(\Omega)\oplus\dot{W}^{m-1+s}_p(\partial\Omega)
\end{equation}
\noindent is so. Furthermore, the estimate 
\begin{equation}\label{estUU-bis}
\|{\mathcal U}\|_{W_p^{m,a}(\Omega)}
\leq C\, \Bigl(\|{\mathcal F}\|_{V_p^{-m,a}(\Omega)} 
+\|g\|_{\dot{W}^{m-1+s}_p(\partial\Omega)}+\|{\mathcal U}\|_{L_p(\Omega)}\Bigr)
\end{equation}
\noindent holds for any solution ${\mathcal U}\in W_p^{m,a}(\Omega)$
of (\ref{m9}).

\subsection{Further comments and the sharpness of Theorem~\ref{Theorem1}}

A byproduct of our proof of Theorem~\ref{Theorem1} is the following. 
Assume that $\Omega\subset\RR^n$, ${\mathcal L}$ are as in the first 
paragraph of the statement of Theorem~\ref{Theorem1}. Then there exists 
$\varepsilon>0$, depending only on the $L_\infty$-norm of the coefficients 
and the ellipticity constant of ${\mathcal L}$, with the property that 
the Dirichlet problem (\ref{m9}) with data from 
$\dot{W}^{m-1+s}_p(\partial\Omega)$ has a unique solution in 
$W^{m,a}_p(\Omega)$ granted that 
\begin{equation}\label{cond-coeffx2}
|2^{-1}-p^{-1}|<\varepsilon\quad\mbox{ and }\quad|a|<\varepsilon. 
\end{equation}
\noindent To justify this claim, we rely on Theorem~\ref{P1-az} and, using 
the same strategy as before, reduce matters to proving that the operator 
(\ref{cal-L}) is an isomorphism. When $a=0$ and $p=2$, our assumptions on 
${\mathcal A}(X,D_X)$ and the classical Lax-Milgram lemma ensure that
this is indeed the case. Then the stability theory from {\bf{\cite{KM}}}, 
{\bf{\cite{Sn}}} allows us to perturb this result, i.e., conclude that 
(\ref{cal-L}) is an isomorphism whenever (\ref{cond-coeffx2}) holds, 
as soon as we show that the scale $V^{m,a}_p(\Omega)$ is stable under 
complex interpolation. That is, if $1<p_i<\infty$, $-1/p_i<a_i<1-1/p_i$, 
$i\in\{0,1\}$, $\theta\in(0,1)$, $1/p=(1-\theta)/p_0+\theta/p_1$ and 
$a=(1-\theta)a_0+\theta a_1$, then 
\begin{equation}\label{interpol1}
[V^{m,a_0}_{p_0}(\Omega),V^{m,a_1}_{p_1}(\Omega)]_\theta=V^{m,a}_p(\Omega),
\end{equation}
\noindent where $[\cdot,\cdot]_\theta$ denotes the usual complex interpolation
bracket. In the proof of (\ref{interpol1}) we may assume that $\Omega$ is 
a special Lipschitz domain and, further, that $\Omega=\RR^n_+$, by making
the change of variables described in \S{5.2}-\S{5.3}. In this latter 
setting, it will be useful to note that 
\begin{equation}\label{lpch}
[L_{p_0}(\RR^n_+,\,x_n^{a_0 p_0}\,dx),
L_{p_1}(\RR^n_+,\,x_n^{a_1 p_1}\,dx)]_\theta
=L_{p}(\RR^n_+,\,x_n^{a p}\,dx),
\end{equation}
\noindent granted that the indices involved are as before, which follows 
from well-known interpolation results for Lebesgue spaces with change 
of measure (cf. Theorem~5.5.3 on p.\,120 in {\bf{\cite{BL}}}). 
Then (\ref{interpol1}) follows easily from (\ref{lpch}), the fact 
that for each $\alpha\in\NN_0^n$ with $|\alpha|=m$ the operator 
$D^\alpha$ maps the scale $V^{m,a}_{p}(\RR^n_+)$ boundedly into 
the scale $L_{p}(\RR^n_+,\,x_n^{a p}\,dx)$,  
and (\ref{IntUU})-(\ref{KDu}). This finishes the proof of the claim made
at the beginning of this subsection. 

In turn, the aforementioned result can be viewed as an extension of a 
well-known theorem of N.\,Meyers, who has treated the case $m=1$, $l=1$ 
in {\bf{\cite{Mey}}}. The example given in \S{5} of {\bf{\cite{Mey}}} 
shows that the membership of $p$ to a small neighborhood of $2$ is a 
necessary condition, even when $\partial\Omega$ is smooth, if 
the coefficients $A_{\alpha\beta}$ are merely bounded. For higher 
order operators we make use of an example originally due to V.G.\,Maz'ya 
{\bf{\cite{Maz-Ctr}}} (cf. also the contemporary article by E.\,De Giorgi 
{\bf{\cite{DG}}}). Specifically, when 
$m\in\NN$ is even, consider the divergence-form equation
\begin{equation}\label{Maz-Op}
\Delta^{\frac{1}{2}m-1}{\mathcal L}_4\,\Delta^{\frac{1}{2}m-1}{\mathcal U}=0
\quad\mbox{in }\,\,\Omega:=\{X\in\RR^n:\,|X|<1\},
\end{equation}
\noindent where ${\mathcal L}_4$ is the fourth order operator  
\begin{eqnarray}\label{Maz-Op2}
{\mathcal L}_4(X,D_X)\,{\mathcal U} & := & 
a\,\Delta^2{\mathcal U}
+b\sum_{i,j=1}^n\Delta\Bigl(\frac{X_iX_j}{|X|^2}
\partial_i\partial_j\,{\mathcal U}\Bigr)
+b\sum_{i,j=1}^n\partial_i\partial_j\Bigl(\frac{X_iX_j}{|X|^2}\,
\Delta\,{\mathcal U}\Bigr)
\nonumber\\[6pt]
&&+c\sum_{i,j,k,l=1}^n\partial_k\partial_l\Bigl(\frac{X_iX_jX_kX_l}{|X|^4}
\partial_i\partial_j
\,{\mathcal U}\Bigr).
\end{eqnarray}
\noindent Obviously, the coefficients of ${\mathcal L}_4(X,D_X)$ are 
bounded, and if the parameters $a,b,c\in\RR$, $a>0$, are chosen such 
that $b^2<ac$ then ${\mathcal L}$ along with 
$\Delta^{\frac{1}{2}m-1}{\mathcal L}_4\,\Delta^{\frac{1}{2}m-1}$ 
are strongly elliptic. Now, if $W^s_p$ denotes the usual $L_p$-based 
Sobolev space of order $s$, it has been observed in {\bf{\cite{Maz-Ctr}}} 
that the function ${\mathcal U}(X):=|X|^{\theta+m-2}\in W^m_2(\Omega)$ has 
${\rm Tr}\,{\mathcal U}\in C^\infty(\partial\Omega)$ and is a
weak solution of (\ref{Maz-Op}) for the choice 
\begin{equation}\label{Maz-Op3}
\theta:=2-\frac{n}{2}+\sqrt{\frac{n^2}{4}-\frac{(n-1)(bn+c)}{a+2b+c}}.
\end{equation}
\noindent Thus, if $a:=(n-2)^2+\varepsilon$, $b:=n(n-2)$, $c:=n^2$, 
$\varepsilon>0$, the strong ellipticity condition is satisfied and 
$\theta=\theta(\varepsilon)$ becomes 
$2-n/2+n\,\varepsilon^{1/2}/2\sqrt{4(n-1)^2+\varepsilon}$. 
However, ${\mathcal U}\in W^m_p(\Omega)$ if and only if 
$p<n/(2-\theta(\varepsilon))$, and the bound $n/(2-\theta(\varepsilon))$ 
approaches $2$ when $\varepsilon\to 0$. An analogous example can be produced 
when $m>1$ is odd, starting with a sixth order operator 
${\mathcal L}_6(X,D_X)$ from {\bf{\cite{Maz-Ctr}}}. 
In the above context, given that 
$W^1_n(\Omega)\hookrightarrow{\rm VMO}(\Omega)$, it is significant to point 
out that both for the example in {\bf{\cite{Mey}}}, when $n=2$, and for 
(\ref{Maz-Op}) when $n\geq 3$, the coefficients have their 
gradients in weak-$L_n$ yet they fail to belong to $W^1_n(\Omega)$.

Of course, condition (\ref{A-bdd}) ensures that the left-hand side of 
(\ref{a0}) is always finite but it is its actual size which determines
whether for a given pair of indices $s$, $p$, the problem 
(\ref{e0}), (\ref{W-Nr}), (\ref{data-B}) is well-posed. Note that 
the maximum value that the right-hand side of (\ref{a0}) takes for
$0<s<1$ and $1<p<\infty$ occurs precisely when $p=2$ and $a:=1-s-1/p=0$. 
As (\ref{a0}) shows, the set of pairs $(s,1/p)\in(0,1)\times(0,1)$ for which 
(\ref{m9}) is well-posed in the context of Theorem~\ref{Theorem1} 
exhausts the entire square $(0,1)\times(0,1)$ as the distance from 
$\nu$ and the $A_{\alpha\beta}$'s to {\rm VMO} tends to zero 
(while the Lipschitz constant of $\Omega$ and the ellipticity constant of
${\mathcal L}$ stay bounded). That the geometry of the Lipschitz domain 
$\Omega$ intervenes in this process through a condition such as (\ref{a0}) 
confirms a conjecture made 
by P.\,Auscher and M.\,Qafsaoui in {\bf{\cite{AQ}}}.

While the main aim of the present work is the consideration of 
higher-order operators with coefficients in $L_\infty$, 
Theorem~\ref{Theorem1} (and, with it, Theorem~\ref{Theorem}) 
is new even in the case when $m=1$ and 
$A_{\alpha\beta}\in\CC^{l\times l}$ (i.e., for second order, constant 
coefficient systems). It provides a complete answer to the issue of 
well-posedness of the problem (\ref{m9}) in the sense that the small mean 
oscillation condition, depending on $p$ and $s$, is in the nature of best 
possible if one insists on allowing arbitrary indices $p$ and $s$. 
This can be seen by considering the following Dirichlet problem for the 
Laplacian in a domain $\Omega\subset\RR^n$:
\begin{equation}\label{LapJK}
\Delta\,{\mathcal U}=0\mbox{ in }\Omega,\quad
{\rm Tr}\,{\mathcal U}=g\in B^s_p(\partial\Omega),\quad
|{\mathcal U}|+|\nabla {\mathcal U}|\in L_p(\Omega,\,\rho(X)^{p(1-s)-1}\,dX).
\end{equation}
\noindent It has long been known that, already in the case when 
$\partial\Omega$ exhibits one cone-like singularity, the well-posedness 
of (\ref{LapJK}) prevents $(s,1/p)$ from being an arbitrary point 
in $(0,1)\times(0,1)$. At a more sophisticated level, the work 
of D.\,Jerison and C.\,Kenig in the 1990's shows that (\ref{LapJK}) is 
well-posed in an arbitrary, given Lipschitz domain $\Omega$ if and only if 
the point $(s,1/p)$ belongs to a certain open subregion of  
$(0,1)\times(0,1)$, determined exclusively 
by the geometry of the domain $\Omega$ (cf. {\bf\cite{JK}}).

\vskip 0.10in
\noindent --------------------------------------
\vskip 0.20in
\begin{minipage}[t]{7.5cm}

\noindent {\tt Vladimir Maz'ya}

\noindent Department of Mathematics

\noindent Ohio State University 

\noindent Columbus, OH 43210, USA

\noindent {\tt e-mail}: {\it vlmaz@math.ohio-state.edu}

and 

\noindent Department of Mathematical Sciences

\noindent University of Liverpool

\noindent Liverpool L69 3BX, UK

and

\noindent Department of Mathematics

\noindent Link\"oping University

\noindent Link\"oping SE-581 83, Sweden

\end{minipage}
\hfill
\begin{minipage}[t]{7.5cm}

\noindent {\tt Marius Mitrea}

\noindent Department of Mathematics

\noindent University of Missouri at Columbia

\noindent Columbia, MO 65211, USA

\noindent {\tt e-mail}: {\it marius@math.missouri.edu}

\vskip 0.08in

\noindent {\tt Tatyana Shaposhnikova}

\noindent Department of Mathematics

\noindent Link\"oping University

\noindent Link\"oping SE-581 83, Sweden

\noindent {\tt e-mail}: {\it tasha@mai.liu.se}

and

\noindent Department of Mathematics

\noindent Ohio State University 

\noindent Columbus, OH 43210, USA

\noindent {\tt e-mail}: {\it tasha@math.ohio-state.edu}

\end{minipage}


\begin{thebibliography}{999}
\small

\bibitem{AP} V.\,Adolfsson and J.\,Pipher, {\it The inhomogeneous Dirichlet 
problem for $\Delta^2$ in Lipschitz domains}, J. Funct. Anal., 159 (1998),  
no.\,1, 137--190. 

\bibitem{Ag1} S.\,Agmon, {\it Multiple layer potentials and the Dirichlet 
problem for higher order elliptic equations in the plane. I},  Comm. Pure 
Appl. Math., 10 (1957), 179--239.

\bibitem{Ag2} S.\,Agmon, {\it The $L^{p}$ approach to the Dirichlet 
problem. I. Regularity theorems}, Ann. Scuola Norm. Sup. Pisa (3), 13 (1959),
405--448.

\bibitem{ADN} S.\,Agmon, A.\,Douglis and L.\,Nirenberg, {\it Estimates near
the boundary for solutions of elliptic partial differential equations
satisfying general boundary conditions. II}, Comm. Pure Appl. Math.,
17 (1964), 35--92.

\bibitem{AQ} P.\,Auscher and M.\,Qafsaoui, {\it Observations on estimates
for divergence elliptic equations with $VMO$ coefficients}, Boll. Unione Mat.
Ital. Sez. B Artic. Ric. Mat., 5 (2002), 487--509.

\bibitem{BL} J.\,Bergh and J.\,L\"{o}fstr\"{o}m, {\it Interpolation Spaces.
An Introduction}, Springer-Verlag, 1976.

\bibitem{Br} F.E.\,Browder, {\it The Dirichlet problem for linear elliptic 
equations of arbitrary even order with variable coefficients}, 
Proc. Nat. Acad. Sci. U.\,S.\,A., 38 (1952), 230--235.

\bibitem{BG} A.\,Buffa and G.\,Geymonat, {\it On traces of functions in 
$W^{2,p}(\Omega)$ for Lipschitz domains in ${\mathbb R}^3$}, 
C. R. Acad. Sci. Paris S\'er. I Math., 332 (2001), no. 8, 699--704.

\bibitem{By1} S.\,Byun, {\it Elliptic equations with ${\rm BMO}$ coefficients
in Lipschitz domains}, Trans. Amer. Math. Soc., 357 (2005), 1025--1046.

\bibitem{CaPe} L.A.\,Caffarelli and I.\,Peral, {\it On estimates for elliptic 
equations in divergence form}, Comm. Pure Appl. Math., 51 (1998), 1--21.

\bibitem{CFL1} F.\,Chiarenza, M.\,Frasca and P.\,Longo, {\it Solvability
of the Dirichlet problem for nondivergence elliptic equations with
$VMO$ coefficients}, Trans. Amer. Math. Soc., 336 (1993), 841--853.

\bibitem{CRW} R.\,Coifman, R.\,Rochberg, and G.\,Weiss, {\it Factorization
theorems for Hardy spaces in several variables}, Ann. Math., 103 (1976),
611--635.

\bibitem{Dah} B.\,E.\,Dahlberg, {\it Poisson semigroups and singular 
integrals}, Proc. Amer. Math. Soc., 97 (1986), 41--48.

\bibitem{DK} B.\,E.\,Dahlberg and C.\,E.\, Kenig, {\it Hardy spaces and the 
Neumann problem in $L\sp p$ for Laplace's equation in Lipschitz domains}, 
Ann. of Math., 125 (1987), no.\,3, 437--465. 

\bibitem{DKPV} B.\,E.\,Dahlberg, C.\,E.\,Kenig, J.\,Pipher and 
G.C.\,Verchota, {\it Area integral estimates for higher order elliptic 
equations and systems}, Ann. Inst. Fourier (Grenoble), 47 (1997), no.\,5, 
1425--1461.

\bibitem{DG} E.\,De Giorgi, {\it Un esempio di estremali discontinue per 
un problema variazionale di tipo ellittico}, Boll. Un. Mat. Ital. (4), 
1 (1968), 135--137.

\bibitem{Faz} G.\,Di Fazio, {\it Estimates for divergence form elliptic
equations with discontinuous coefficients}, Boll. Un. Mat. Ital A(7), 10
(1996), 409--420.

\bibitem{Fe} C.\,Fefferman, {\it A sharp form of Whitney's extension theorem},
Annals of Math., 161 (2005), 509--577.

\bibitem{FS} C.\,Fefferman and E.M.\,Stein, {\it $H^{p}$ spaces of several 
variables}, Acta Math., 129 (1972), no.\,3-4, 137--193. 

\bibitem{Ga} E.\,Gagliardo, {\it Caratterizzazioni delle tracce sulla 
frontiera relative ad alcune classi di funzioni in $n$ variabili}, 
Rend. Sem. Mat. Univ. Padova, 27 (1957), 284--305.

\bibitem{Gar} L.\,G{\aa}rding, {\it Dirichlet's problem for linear elliptic 
partial differential equations}, Math. Scand., 1 (1953), 55--72.

\bibitem{Gr} P.\,Grisvard, {\it Elliptic Problems in Nonsmooth Domains, 
Monographs and Studies in Mathematics}, Vol.\,24, Pitman, Boston, MA, 1985.

\bibitem{GK} G.\,Grubb and N.J.\,Kokholm, {\it A global calculus of 
parameter-dependent pseudodifferential boundary problems in $L_p$ Sobolev 
spaces}, Acta Math., 171 (1993), 165--229. 

\bibitem{Hed-1} L.I.\,Hedberg, {\it Spectral synthesis in Sobolev spaces, 
and uniqueness of solutions of the Dirichlet problem}, Acta Math., 147 (1981),
no.\,3-4, 237--264. 

\bibitem{Hed-2} L.I.\,Hedberg, {\it On the Dirichlet problem for higher-order 
equations}, Conference on harmonic analysis in honor of Antoni Zygmund, 
pp.\,620--633, Wadsworth Math. Ser., Wadsworth, Belmont, CA, 1983.

\bibitem{HL} S.\,Hofmann and J.L.\,Lewis, {\it $L^2$ solvability and 
representation by caloric layer potentials in time-varying domains}, 
Ann. of Math., 144 (1996), no.\,2, 349--420.

\bibitem{Ho} L.\,H\"ormander, {\it Linear Partial Differential Operators}, 
Springer Verlag, Berlin-New York, 1976. 

\bibitem{IS} T.\,Iwaniec and C.\,Sbordone, {\it Riesz transforms and 
elliptic PDEs with VMO coefficients}, J. Anal. Math., 74 (1998), 183--212.

\bibitem{JK} D.\,Jerison and C.\,Kenig, {\it The inhomogeneous Dirichlet 
problem in Lipschitz domains}, J. Funct. Anal., 130 (1995), 161--219.

\bibitem{JW} A.\,Jonsson and H.\,Wallin, {\it Function Spaces on Subsets of
$\RR^n$}, University of Umea, 1984.

\bibitem{Ke} C.E.\,Kenig, {\it  Harmonic Analysis Techniques for Second 
Order Elliptic Boundary Value Problems}, CBMS Regional Conference Series in 
Mathematics, Vol.\,83, AMS, Providence, RI, 1994.

\bibitem{KM} V.\,Kozlov and V.\,Maz'ya, {\it Asymptotic formula for solutions
to elliptic equations near Lipschitz boundary}, Ann. Mat. Pura Appl., 184  
(2005),  185--213.

\bibitem{KMR1} V.A.\,Kozlov, V.G.\,Maz'ya and J.\,Rossmann, {\it Elliptic 
Boundary Value Problems in Domains with Point Singularities}, AMS, 1997.

\bibitem{KMR2} V.A.\,Kozlov, V.G.\,Maz'ya and J.\,Rossmann, {\it Spectral 
Problems Associated with Corner Singularities of Solutions to Elliptic 
Equations}, AMS, 2001.

\bibitem{MPS} A.\,Maugeri, D.K.\,Palagachev and L.G.\,Softova, {\it Elliptic 
and Parabolic Equations with Discontinuous Coefficients}, Wiley-VCH, 2000.

\bibitem{Maz-Ctr} V.\,Maz'ya, {\it Examples of nonregular solutions of 
quasilinear elliptic equations with analytic coefficients}, Funkcional. 
Anal. i Prilozen., 2 (1968) 53-57; translated in  
Functional Anal. Appl., 2 (1968), 230-234.

\bibitem{Maz1} V.\,Maz'ya, {\it Sobolev Spaces}, Springer Series in Soviet 
Mathematics. Springer-Verlag, Berlin, 1985. 

\bibitem{Mazz-4} V.G.\,Maz'ya, {\it The Wiener test for higher order elliptic 
equations}, Duke Math. J., 115 (2002), no. 3, 479--512.

\bibitem{MS} V.\,Maz'ya and T.\,Shaposhnikova, {\it Theory of Multipliers in
Spaces of Differentiable Functions}, Monographs and Studies in Mathematics
Vol.\,23, Pitman, Boston, MA, 1985.

\bibitem{MS1} V.\,Maz'ya and T.\,Shaposhnikova, {\it On the regularity of 
the boundary in the $L_p$-theory of elliptic boundary value problems}. I. 
(Russian) Partial differential equations, pp.\,39--56, Trudy Sem. 
S.L.\,Soboleva, No.\,2, 80, Akad. Nauk SSSR Sibirsk. Otdel., Inst. Mat., 
Novosibirsk, 1980.   

\bibitem{MS2} V.\,Maz'ya and T.\,Shaposhnikova, {\it Higher regularity in 
the classical layer potential theory for Lipschitz domains},  
Indiana Univ. Math. J., 54 (2005), no.\,1, 99--142.

\bibitem{Mey} N.G.\,Meyers, {\it An estimate for the gradient of solutions of 
second order elliptic divergence equations}, Ann. Scuola Norm. Sup. Pisa (3), 
17 (1963), 189--206.

\bibitem{MMT} D.\,Mitrea, M.\,Mitrea and M.\,Taylor, {\it Layer Potentials, 
the Hodge Laplacian, and Global Boundary Problems in Nonsmooth Riemannian 
Manifolds}, Mem. Amer. Math. Soc., Vol.\,150, No.\,713, 2001.

\bibitem{MT1} M.\,Mitrea and M.\,Taylor, {\it Potential theory on Lipschitz
domains in Riemannian manifolds: Sobolev-Besov space results and the Poisson
problem}, J. Funct. Anal. 176 (2000), 1--79.

\bibitem{MT2} M.\,Mitrea and M.\,Taylor, {\it Sobolev and Besov space 
estimates for solutions to second order PDE on Lipschitz domains in manifolds 
with Dini or H{\"o}lder continuous metric tensors}, 
Comm. in PDE, 30 (2005), 1--37.

\bibitem{MT3} M.\,Mitrea and M.\,Taylor, {\it The Poisson problem in
weighted Sobolev spaces on Lipschitz domains}, Indiana Univ. Math. J., 
55 (2006), 1063--1089. 

\bibitem{Nec} J.\,Ne\v{c}as, {\it Les M\'ethodes Directes en Th\'eorie 
des \'Equations Elliptiques}, Masson et Cie, \'Editeurs, Paris, Academia, 
\'Editeurs, Prague, 1967.

\bibitem{Ni} S.M.\,Nikol'ski\u{\i}, {\it Approximation of Functions of Several 
Variables and Imbedding Theorems}, 2-nd edition, revised and supplemented, 
Nauka, Moscow, 1977. 

\bibitem{PV} J.\,Pipher and G.C.\,Verchota, {\it Dilation invariant estimates 
and the boundary G{\aa}rding inequality for higher order elliptic operators},  
Ann. of Math., 142 (1995), no. 1, 1--38.

\bibitem{Sar} D.\,Sarason, {\it Functions of vanishing mean oscillation}, 
Trans. Amer. Math. Soc., 207 (1975), 391--405.

\bibitem{Sh} Z.\,Shen, {\it Necessary and sufficient conditions for the 
solvability of the $L_p$ Dirichlet problem on Lipschitz domains}, 
Math. Ann., 336 (2006), 697--725.

\bibitem{Sn} I.Ya.\,\v{S}neiberg, {\it Spectral properties of linear
operators in interpolation families of Banach spaces}, Mat. Issled.,
9 (1974), 214--229.

\bibitem{Sob} S.L.\,Sobolev, {\it On a boundary value problem for the 
polyharmonic equations}, Mat. Sb., 2 (1937), 467--499.

\bibitem{Sob2} S.L.\,Sobolev, {\it Applications of Functional Analysis in 
Mathematical Physics}, Izd. LGU, Leningrad, 1950. 

\bibitem{Sol1} V.A.\,Solonnikov, {\it General boundary value problems
for systems elliptic in the sense of A. Douglis and L. Nirenberg. I},
(Russian) Izv. Akad. Nauk SSSR, Ser. Mat., 28 (1964), 665-706.

\bibitem{Sol2} V.A.\,Solonnikov, {\it General boundary value problems for
systems elliptic in the sense of A. Douglis and L. Nirenberg. II}, (Rusian)
Trudy Mat. Inst. Steklov, Vol.\,92 (1966), 233--297.

\bibitem{St} E.M.\,Stein, {\it Singular Integrals and Differentiability
Properties of Functions}, Princeton Mathematical Series, No.\,30,
Princeton University Press, Princeton, N.J., 1970.

\bibitem{Tr} H.\,Triebel, {\it Interpolation Theory, Function Spaces, 
Differential Operators}, North-Holland Publishing Co., Amsterdam-New York, 
1978. 

\bibitem{Usp} S.V.\,Uspenski\u{\i}, {\it Imbedding theorems for classes with
weights}, Trudy Mat. Inst. Steklov., 60 (1961), 282--303.

\bibitem{Ve2} G.C.\,Verchota, {\it The Dirichlet problem for the polyharmonic
equation in Lipschitz domains}, Indiana Univ. Math. J., 39 (1990), 
671--702.

\bibitem{Ve} G.C.\,Verchota, {\it The biharmonic Neumann problem in 
Lipschitz domains}, Acta Math., 194 (2005),  217--279.

\bibitem{Vi} M.\,I.Vi\v{s}ik, {\it On strongly elliptic systems of 
differential equations}, Mat. Sbornik N.S., 29 (1951), 615--676. 

\bibitem{Wh} H.\,Whitney, {\it Analytic extensions of differentiable 
functions defined in closed sets}, Trans. Amer. Math. Soc., 36 (1934), 63--89. 

\bibitem{Wh2} H.\,Whitney, {\it Functions differentiable on the boundaries 
of regions}, Ann. of Math., 35 (1934), 482--485.

\end{thebibliography}
\end{document}